\documentclass[a4paper,12pt]{article}
\usepackage{amsmath,amssymb,theorem}
\usepackage{amscd}
\usepackage{latexsym} 
\usepackage{amsfonts}
\usepackage{graphicx}
\usepackage{layout}

\newcommand{\bma}{\begin{math}}
\newcommand{\ema}{\end{math}}
\newcommand{\beq}{\begin{eqnarray}}
\newcommand{\eeq}{\end{eqnarray}}
\newcommand{\be}{\begin{eqnarray*}}
\newcommand{\ee}{\end{eqnarray*}}

\theorembodyfont{\itshape}
\newtheorem{thm}{Theorem}
\newtheorem{prop}[thm]{Proposition}
\newtheorem{lm}[thm]{Lemma}
\newtheorem{co}[thm]{Corollary}
\theorembodyfont{\rmfamily}
\newtheorem{df}[thm]{Definition}
\newtheorem{pf}{Proof.}
\newtheorem{re}{Remark}

\newtheorem{N}{Notation}

\newcommand{\qed}{\hfill $\Box$ \par}

\title{Hypergeometric Integrals associated with Hypersphere
Arrangements and Cayley-Menger Determinants}
\author{Kazuhiko AOMOTO and Yoshinori MACHIDA}
\date{}

\begin{document}

\maketitle

\begin{abstract}
The $n$ dimensional hypergeometric integrals associated with a hypersphere arrangement $S$
are formulated by the pairing of 
$n$ dimensional twisted cohomology $H_\nabla^n (X, \Omega^\cdot (*S))$
and its dual. Under the condition of general position there are stated some results which concern
an explicit representation of the standard form  by a special (NBC) basis of the twisted cohomology
(contiguity relation in positive direction),
the variational formula of the corresponding integral in terms of special invariant $1$-forms
$\theta_J$ written by  Calyley-Menger minor determinants,
and a connection relation of the unique twisted $n$-cycle identified with the unbounded chamber
to a special basis of  twisted $n$-cycles  identified with
bounded chambers. Gauss-Manin connection is formulated and  is  explicitly 
presented in two simplest cases.  In the appendix 
 contiguity relation in negative direction is presented in terms of
Cayley-Menger determinants.
\end{abstract}

%{\scriptsize 2010{ \bf Mathematical Subject Classification:} Primary 14F40, 33C70; Secondary 14H70.}

%{\scriptsize 2010 { \bf Key Words:} 
%hypergeometric integral, hypersphere arrangement, twisted rational 
%de Rham cohomology, 
%Cayley-Menger determinant, contiguity relation, Gauss-Manin connection. }

%{\scriptsize 2010 {\bf Running Title:} Hypergeometric Integral of Sphere Arrangement}

\renewcommand{\thefootnote}{\fnsymbol{footnote}}



\footnotetext{\scriptsize
{\it Running Title}: Hypergeometric Integral of Sphere Arrangement.}

\bigskip

\section{Introduction}

Let ${\alpha}_j=(\alpha_{j1},\alpha_{j2},\ldots,\alpha_{jn})\ (1\le j\le m)$ be
$m\,( n+1 \le m)$ pieces of $n$ dimensional real vectors, and $\alpha_{j0}
 \in {\bf R}$.
We take quadratic polynomials $f_j(x)$ of
$x=(x_1,x_2,\ldots,x_n)\in {\bf C}^n:$
\be
f_j (x)= Q(x)+2(\alpha_j,x)+\alpha_{j0},
\ee
where we put $Q(x) = (x,x)=\sum_{\nu=1}^n x_\nu^2,\; (\alpha_j,x)=\sum_{\nu=1}^n\alpha_{j\nu} x_\nu$.

\medskip

For  $\lambda=(\lambda_1,\ldots,\lambda_{m}) \in {\bf C}^{m}$,
we consider an analytic integral over a twisted $n$-cycle $\mathfrak{z}:$
\beq
&&{\cal J}_\lambda(\varphi) = \int_{\mathfrak{z}} \Phi(x)\varphi(x) \varpi,\\
&& \varpi = dx_1\wedge dx_2\wedge\cdots\wedge dx_n\nonumber,
\eeq
associated with the multiplicative meromorphic function
\be
\Phi(x)= f_1(x)^{\lambda_1} f_2(x)^{\lambda_2}\cdots  f_{m}(x)^{\lambda_{m}} \quad(\lambda_j\in {\bf C}).
\ee
Here $\varphi(x)$ is a rational function on ${\bf C}^n$
which is holomorphic in the complement
$X = {\bf C}^n - S\,(S=\bigcup_{j=1}^{m}{S_j})$,
 and 
$S_j$ denotes the $n-1$ dimensional complex hypersphere
\be
S_j :  f_j(x)=0.
\ee

The RHS of (1) will also be denoted by $\langle\varphi, \mathfrak{z}\rangle$.
We put $\lambda_\infty$ to be $\sum_{j=1}^{m} \lambda_j.$
From now on, we assume no resonance condition
\be
&&\lambda_j \ne 0,1,2,3,\ldots\quad(1\le j\le m),\\
&& 2\lambda_\infty \ne  0,-1,-2,-3,\cdots.
\ee

For the formulation of the integral (1), we use
the twisted de Rham cohomology on $X$
denoted by 
$H_\nabla^*(X,\Omega^\cdot(*S))$,
namely the cohomology related to
the twisted exterior differentiation $\nabla$ :
\be
\nabla \psi = d\psi + d\log\Phi \wedge \psi
\ee
over
the complex of rational differential forms
on the complex domain $X= {\bf C}^n - S,$ and
\be
\Omega^\cdot(*S) = \oplus_{\nu=0}^n \Omega^\nu(* S) .
\ee

In particular, ${\cal R}= \Omega^0(*S) $
is the set of all rational functions
which are holomorphic on $X$.
Let $\cal L$ denote the local system on $X$ defined by $\Phi(x)$
and ${\cal L}^*$ be its dual.  We can take $\mathfrak{z}$ as an element
of the $n$ dimensional homology $H_n(X,{\cal L}^*)$
with coefficients  in ${\cal L}^*$ (see [10], [13] for basic notions and details).

A theorem due to Grothendieck-Deligne (see [13])
shows that the following isomorphism is valid :
\be
H^*(X,{\cal L})\cong H_\nabla^*(X,\Omega^\cdot(*S)),
\ee
and that $H^*(X,{\cal L})$ and $H_*(X,{\cal L}^*)$
are dual to each other.

Since the integral (1) admits the group of isometric transformations in the Euclidean space ${\bf R}^n$
or its complexification ${\bf C}^n$, differential equations or contiguity relations can be represented 
by invariants related to $\Phi(x)$ under
this group .

\bigskip

Denote by ${\Re}S_j$ the real $n-1$ dimensional hypersphere
being the real part of $S_j$,
by $O_j$ its center $-\alpha_j=(-\alpha_{j1},-\alpha_{j2},\cdots,-\alpha_{jn})$. Then
the radius of ${\Re}S_j$
and the distance between $O_j$ and $O_k$
are given by
\be
&&r_j^2 =Q(\alpha_j)-\alpha_{j0},\\
&& \rho_{j,k}^2=Q(\alpha_j-\alpha_k)
\ee
respectively. These are the basic invariants with respect to the group of isometry.

%def1
\begin{df}

Let $B=(b_{jk})_{1\le j,k\le m+2}$ be the Cayley-Menger matrix of order $m + 2$
associated with the hypersphere arrangement $S$, i.e.,
the symmetric matrix of order $m + 2$ such that 
\be
&&b_{jj} = 0\,(1\le j\le m + 2),\ b_{1j} =1\,(2\le j\le m + 2),\ b_{2j} = r_{j-2}^2\,(3\le j\le m + 2),\\
&&b_{jk} = \rho_{j-2,k-2}^2\,(3\le j<k\le m + 2).
\ee
 For convenience in the sequel, the set of subscripts will de denoted by $0, \star, 1,2,\ldots, m$
 instead of $1,2,\ldots, m + 2$ respectively.
 
 {\rm Cayley-Menger} determinants
with the components $r_j^2$
and 
$\rho_{j,k}^2$ are defined as follows ([12],[14],[15]):

\be
B\left(
\begin{array}{cccc}
0  & i_1&  \ldots & i_p  \\
0 & j_1 &  \ldots & j_p   
\end{array}
\right)
=
\left|
\begin{array}{ccccc}
0 & 1 & 1 &\ldots & 1\\
1 & \rho_{i_1,j_1}^2  & \rho_{i_1,j_2}^2   & \ldots & \rho_{i_1,j_p}^2  \\
1 &  \rho_{i_2,j_1}^2 & \rho_{i_2,j_2}^2  & \ldots& \rho_{i_2,j_p}^2  \\
\vdots& \vdots  & \vdots  &\ddots & \vdots \\
1& \rho_{i_p,j_1}^2 & \rho_{i_p,j_2}^2&\ldots & \rho_{i_p,j_p}^2
\end{array}
\right|,
\ee
\be
B\left(
\begin{array}{cccc}
\star  & i_1&  \ldots & i_p  \\
0  & j_1, &  \ldots & j_p   
\end{array}
\right)
=
\left|
\begin{array}{ccccc}
1 & r_{j_1}^2 & r_{j_2}^2 &\ldots & r_{j_p}^2\\
1 & \rho_{i_1,j_1}^2  & \rho_{i_1,j_2}^2   & \ldots & \rho_{i_1,j_p}^2  \\
1 &  \rho_{i_2,j_1}^2 & \rho_{i_2,j_2}^2  & \ldots& \rho_{i_2,j_p}^2  \\
\vdots& \vdots  & \vdots  &\ddots & \vdots \\
1& \rho_{i_p,j_1}^2 & \rho_{i_p,j_2}^2&\ldots & \rho_{i_p,j_p}^2
\end{array}
\right|,
\ee
\be
B\left(
\begin{array}{cccc}
\star & i_1  & \ldots & i_p  \\
\star  & j_1  &  \ldots & j_p   
\end{array}
\right)
=
\left|
\begin{array}{ccccc}
0 & r_{j_1}^2 & r_{j_2}^2 &\ldots & r_{j_p}^2\\
r_{i_1}^2 & \rho_{i_1,j_1}^2  & \rho_{i_1,j_2}^2   & \ldots & \rho_{i_1,j_p}^2  \\
r_{i_2}^2 &  \rho_{i_2,i_1}^2 & \rho_{i_2,j_2}^2  & \ldots& \rho_{i_2,j_p}^2  \\
\vdots& \vdots  & \vdots  &\ddots & \vdots \\
r_{i_p}^2& \rho_{i_p,i_1}^2 & \rho_{i_p,j_2}^2&\ldots & \rho_{i_p,j_p}^2
\end{array}
\right|,
\ee
\be
B\left(
\begin{array}{ccccc}
0&\star & i_1  & \ldots& i_p  \\
0&\star  & j_1  &  \ldots & j_p   
\end{array}
\right)
=
\left|
\begin{array}{cccccc}
0& 1&1&1&\ldots&1\\
1&0 & r_{j_1}^2 & r_{j_2}^2 &\ldots & r_{j_p}^2\\
1&r_{i_1}^2 & \rho_{i_1,j_1}^2  & \rho_{i_1,j_2}^2   & \ldots & \rho_{i_1,j_p}^2  \\
1&r_{i_2}^2 &  \rho_{i_2,i_1}^2 & \rho_{i_2,j_2}^2  & \ldots& \rho_{i_2,j_p}^2  \\
\vdots& \vdots  & \vdots  &\ddots & \vdots \\
1&r_{i_p}^2& \rho_{i_p,i_1}^2 & \rho_{i_p,j_2}^2&\ldots & \rho_{i_p,j_p}^2
\end{array}
\right|.
\ee
%p4
When $i_\nu=j_\nu\;(1\le\nu\le p)$,
we simply write
$B(0i_1\ldots i_p), B( i_1\ldots i_p),
B(0\star i_1\ldots i_p)$
instead of
$
B\left(
\begin{array}{ccccc}
0& i_1  & \ldots& i_p  \\
0& i_1  &  \ldots & i_p   
\end{array}
\right),
B\left(
\begin{array}{ccccc}
 i_1  & \ldots& i_p  \\
 i_1  &  \ldots & i_p   
\end{array}
\right),
B\left(
\begin{array}{ccccc}
0&\star& i_1  & \ldots& i_p  \\
0&\star&  i_1  &  \ldots & i_p   
\end{array}
\right)
$ etc.

\end{df}

For example,
\be
&&B(0\star j)=2r_j^2,\quad B(0ij)= 2\rho_{ij}^2,\\
&& 
B\!
\left(
\begin{array}{ccc}
 0 &i   & k  \\
 0 & j  & k  
\end{array}
\right) 
=\rho_{ik}^2+\rho_{jk}^2 - \rho_{ij}^2,\quad 
B\!
\left(
\begin{array}{ccc}
 0 & \star  & i  \\
 0 & \star  & j  
\end{array}
\right)
=r_i^2+r_{j}^2 - \rho_{ij}^2,
\\
&& 
B\!
\left(
\begin{array}{ccc}
 0 &j   & k  \\
 0 & \star & k 
\end{array}
\right)
=\rho_{jk}^2+r_{k}^2 - r_{j}^2,
\\
&&B(0jkl)= \rho_{kl}^4 + \rho_{jl}^4 + \rho_{jk}^4
- 2\rho_{jk}^2\rho_{jl}^2 - 2 \rho_{jk}^2\rho_{kl}^2
-2\rho_{jl}^2\rho_{kl}^2.
\ee

\bigskip

As is seen below, the arrangement of hyperspheres $\{S_j\}_{1\le j\le  m}$
is equivalent to
an arrangement of hyperplane sections of the fundamental
unit hypersphere $\tilde{S}_0^n$ below
through the stereographic projection.

Let 
\be
&&\iota: x_j = \frac{\xi_j}{\xi_{n+1}+1}\quad(1\le j\le n),\\
&&\iota^{-1} : \xi_j = \frac{2x_j}{Q(x)+1}\;(1\le j\le n),\; \xi_{n+1}=
\frac{1-Q(x)}{1+Q(x)}\quad(\xi= (\xi_1,\xi_2,\ldots,\xi_{n+1})\in \tilde{S}_0^n)
\ee
be the stereographic projection onto ${\bf C}^n$
from the fundamental unit hypersphere
\be
\tilde{S}_0^n :\quad \xi_1^2+\xi_2^2+\cdots +\xi_{n+1}^2=1.
\ee

We have the conformal isomorphism
\be
\iota: \tilde{S}_0^n - \{\xi_{n+1}+1=0\} \cong  {\bf C}^n - \{Q(x)+1=0\},
\ee
such that the quadratic functions $f_j(x)$ on ${\bf C}^n$ are transformed into
the  linear functions $\tilde{f}(\xi)$ on $\tilde{S}_0^n$:

\be
&&\tilde{f}_j(\xi)= \frac{(\xi_{n+1}+1)}{2r_j}\,f_j\!\left(\frac{\xi_1}{\xi_{n+1}+1},\frac{\xi_2}{\xi_{n+1}+1},\ldots,\frac{\xi_n}{\xi_{n+1}+1}\right)\\
&&{}\quad \quad= u_{j0}+\sum_{\nu=1}^{n+1}
u_{j\nu}\xi_\nu\,\,{\mbox{(\,{\rm linear\, function}\,)}},\\
&&u_{j0}=\frac{1+\alpha_{j0}}{2r_j},\quad
u_{j\nu}=\frac{2\alpha_{j\nu}}{2r_j}\;(1\le j\le m),\quad
u_{j\,n+1}=\frac{\alpha_{j0}-1}{2r_j}.
\ee
We now normalize $\tilde{f}_j$
such that they are invariant under the standard
Lorentz transformations :
\be
-u_{j0}^2 + \sum_{\nu=1}^{n+1} u_{j\nu}^2=1.
\ee

It is convenient to put
\be
\tilde{f}_{m+1}(\xi)= \xi_{n+1}+1.
\ee

Then $\tilde{S}_j  : \{\tilde{f}_j = 0\}\cap \tilde{S}_0^n\,(1\le j\le m)$ define
hyperplane sections of $\tilde{S}_0^n$
whose image by $\iota$ is $S_j.$

The Lorentz inner product
$a_{jk}=(\tilde{f}_j,\tilde{f}_k)$ between $\tilde{f}_j,\tilde{f}_k$
can be defined as
\be
&&a_{jk}= -u_{j0}u_{k0}+\sum_{\nu=1}^{n+1} u_{j\nu}u_{k\nu}\\
&&{}\quad \quad= \frac{r_j^2+r_k^2-\rho_{jk}^2}{2r_jr_k}
\quad(1\le j,k\le m+1)
\ee
such that $a_{jj}=1\,(1\le j\le m), a_{m+1\,m+1}=0$.
We put further
\be
a_{j0}=u_{j0}\quad(1\le j\le m+1),\, a_{00}=-1
\ee
and obtain the $(m+2)\times (m+2)$ symmetric configuration matrix
$A=(a_{jk})_{0\le j,k\le m+1}$.

Denote the minor determinant 
of $i_1,\ldots,i_p$th row, $j_1,\ldots,j_p$th column
by
$A\!\left(
\begin{array}{ccc}
 i_1 & \ldots  & i_p  \\
 j_1 & \ldots  & j_p    
\end{array}
\right).$

As is well-known, all the singularities appearing in the integral
 (1) are contained in the zero set
 of at least one of the principal minor determinants (so called Landau
 singularity, see [22]):
\be
 A(i_1,\ldots,i_p) = A\!\left(
\begin{array}{ccc}
 i_1 & \ldots  & i_p  \\
 i_1 & \ldots  & i_p    
\end{array}
\right) =0.
\ee

\bigskip

We denote by $\varepsilon_j\,(1\le j\le  m)$ the standard basis of ${\bf C}^{m}$,
the space of exponents $\lambda$:\ $\lambda = \sum_{j=1}^{m} \lambda_j\varepsilon_j.$
We denote by $N_m$ the set of all indices $\{1,2,\ldots, m\}.$
For $J =\{j_1,\ldots,j_p\}\subset N_m,$ denote by $|J| = p$ the size of $J.$
$J^c$ denotes the complement of $J$ in $N_m $.
$\partial_{j_\nu}J$ denotes the ordered subset $\{j_1,\ldots,j_{\nu-1},j_{\nu+1},\ldots,j_p\}$ 
the $\nu$ th element $j_\nu$ being deleted.
Denote by $\mathfrak{S}_{m}$ the symmetric group of degree $m$
generated by  $\sigma_{ij}$ the transposition between the indices $i,j\,(i\ne j)$.

A subset of indices $J\subset N_m$ such that $1\le |J|\le n+1$
will be called \lq\lq admissible". 
We assume that none of the  principal 
Cayley-Menger determinants vanish. Namely  
\be
({\cal H}1) \quad  B(0J)\ne 0,\; B(0\star J)\ne0
\ee
for all $J\subset N_m$ and $1 \le |J| \le n+1$.

For two ordered subsets of the same size $J =\{j_1,\ldots,j_p\}\subset N_m,\,K=\{k_1,\ldots,k_p\} \subset N_m,$
we say that 
``$J$ and $K$ are equivalent"
 if the sequence $k_1,\ldots,k_p$ are
 a permutation of the sequence $j_1,\ldots,j_p,$ and denote 
 $J \equiv K.$ We may also denote $J \equiv K,$ when $J$ and $K$ are ordered 
and  unordered subsets
respectively and any sequence of elements of $K$ is equivalent to $J.$ 

Through $\iota,$ the arrangement of (complex) spheres in ${\bf C}^n$ is identified by 
an arrangement of hyperplane sections in $\tilde{S}_0^n.$
Since the Euler number $Eu$(X) of $X$ and  of the complement $\tilde{S}_0^n - \bigcup_{j=1}^m\tilde{S}_j$
are related as
\be
Eu(X) = Eu\{\tilde{S}_0^n - \bigcup_{j=1}^m\tilde{S}_j\}-1= (-1)^n \,\bigl(\sum_{\nu=1}^n \, {m\choose\nu}
+  {m-1\choose n}\bigr),
\ee
we have the equality (see [8] Proposition 2.2)
\be
\dim H_\nabla^n(X,\Omega^\cdot(*S)) = \sum_{\nu=1}^n \, {m\choose\nu}
+  {m-1\choose n}.
\ee

Now we take the more restrictive hypothesis which is crucial 
in the sequel:
\be
({\rm \cal H}2)\quad (-1)^{p-1}\, B(0\star J) > 0 \quad(J\subset N_m, |J|=p \le n+1).
\ee

Let $D_j,\,\overline{D}_j$
denote the $n$-open ball and its closure bounded by $\Re S_j$
respectively.
Take an arbitrary $J$ such that $J \subset N_m$ and $|J| = n+1.$
Then the common part $\bigcap_{j\in J}\,D_j$
is non-empty and each $n$-tuple  $\bigcap_{j\in \partial_\nu J}\,\Re S_j \,(\nu\in J)$,
$\nu$ being fixed,  consists of  two different points.
All bounded chambers in ${\bf R}^n - \bigcup_{j\in N_m}\mathfrak{ R}S_j$ make as representative
a basis of the twisted $n$-homology $H_n(X,{\cal L}^*)$ (see [7],[8],[10]).
Their total number is equal to the dimension of $H_\nabla^n(X,\Omega^\cdot(*S)).$

\medskip

${\cal J}_\lambda(\varphi)$ can be regarded as an analytic function 
of the parameters $\alpha_{j\nu}$.
Denote by $d_B$ the total differentiation with respect to the parameters,
then 
\beq
d_B {\cal J}_\lambda(\varphi) = \int_{\mathfrak{z}}\Phi\,\nabla_B(\varphi\varpi),
\eeq
where $\nabla_B$ denotes the covariant differentiation (Gauss-Manin connection)
operating
on $H_\nabla^n(X,\Omega^\cdot(*S))$:
\be
\nabla_B(\varphi\varpi) = d_B\,(\varphi\varpi) + 
d_B\log\Phi \ \varphi\,\varpi,
\ee
where
\be
d_B = \sum_{j=1}^{m}\sum_{\nu=0}^{n+1}d\alpha_{j\nu}\frac{\partial}{\partial \alpha_{j\nu}}.
\ee

Basis of the cohomology $H_\nabla^n(X,\Omega^\cdot(*S))$ cannot be represented any more by
logarithmic forms  contrary to the case of hyperplanes (see [18]) . Nevertheless we can construct
in a canonical way a basis extending logarithmic forms.  

\bigskip

The plan of this article is as follows. 

 At  first  in \S3, we present Main Results (Theorems 1 and 2)  in the simpler case $m = n+1$
and prove them in \S5 by using contiguity relation  among elements of $H_\nabla^n(X,\Omega^\cdot(*S)).$
An extension of these theorems to general $m \ge n+1$ (Theorems 3 and Theorem 4) will be given in $\S6$.
Gauss-Manin connection is formulated in terms of Cayley-Menger determinants (Theorem 5) 
and its simplest cases
 are concretely presented in $\S7$.
In the appendix, we give contiguity relation in negative direction  (Theorem A), 
which plays an important role
to represent Gauss-Manin connection in terms of the admissible $NBC$ basis $F_J\,(J\in {\cal B})$.

Partial results related to them  and some conjectures generalizing them have been announced in [11].

\bigskip

%sect2 basic facts
\section{Basic Facts}

The following three lemmas are elementary
and easily checked.

%lm2
\begin{lm}
\be
A(j\,k)= \frac{- B(0\star jk)}{4 r_j^2r_k^2}\quad (1\leq j,k\leq m).
\ee

$A(j\,k)=0$ holds if and only if $S_j$ and $S_k$
have contact with each other.
\end{lm}

Remark that
\be
&&- B(0\star jk) =-r_j^4-r_k^4-\rho_{jk}^4
+2r_j^2r_k^2+2r_j^2\rho_{jk}^2+2r_k^2\rho_{jk}^2\\
&&=(r_j+r_k+\rho_{jk})(-r_j+r_k+\rho_{jk})(r_j-r_k+\rho_{jk})
(r_j+r_k-\rho_{jk}).
\ee

More generally,
\begin{lm}
For an admissible $J\,(|J| = p),$
\be
&&2^p \prod_{j\in J}r_j^2\,A(J)= (-1)^{p-1}\,B(0\star J)\quad(J\subset N_m).
\ee

$A(J)=0$ if and only if $S_j\,(j \in J)$
have a common tangent affine subspace of codimension $p-1$.
\end{lm}

Remark that
\be
&&a_{j\,m+1}=-\frac{1}{r_j},\ 
A(j\,m+1)= -\frac{1}{r_j^2}\quad(1\le j\le m),\\
&&A(jk\,m+1)= \frac{\rho_{jk}^2}{r_j^2r_k^2}\quad(1\le j<k\le m),
\ee
so that $A(j\,m+1)=0$, or $A(jk\,m+1)=0$ are equivalent to
$r_j=\infty$, or $
\rho_{j,k}=0$.

%lem4

\begin{lm}
Take an arbitrary $J \subset N_m$ such that $|J| = n+1.$
Then $S_j \, (j\in J)$ have a common point if and only if

\be
 B(0\star J)=0,
\ee
which is also equivalent to the equality $A(J) = 0.$
%as is seen from {\rm Lemma 2}.
The centers of
$S_k$ denoted by $O_k\,(k\in J)$
lie in the same hyperplane
if and only if
\be
B(0\, J)=0.
\ee
\end{lm}

The proof is elementary. For the proof, see [8].

Now we assume that $S_j\;(j\in J)$
are in general position, i.e.,
\be
B(0\,J) \ne 0,\quad B(0\star \,J) \ne 0\quad (J \subset N_m,\, |J| \le n+1).
\ee
Remark that for $|J|=p,$ $(-1)^p B(0J)>0\ (2\le p\le n+1)$
and $(-1)^{p-1}\,B(0\star J) > 0\ (1\le p\le n+1),$
when $\alpha_{j\nu}, \, \alpha_{j0}$ are all real.

\bigskip

%def5
\begin{df}

For $\varphi\varpi \in \Omega^n(*S)$,
$\mathfrak{M}_{f_j^{\pm1}}\,\varphi\varpi$ means the multiplication
of $\varphi\varpi$ by $f_j^{\pm1}$. More generally $\mathfrak{M}_{f_J^{\pm1}} = \prod_{j\in J}\,\mathfrak{M}_{f_j^{\pm1}}$
for $J \subset N_m$.
We shall  simply abbreviate by $J\,K$ the ordered join $J\,\cup\, K$ of
two disjoint sets of indices $J,\,K$.

For $\varphi\,\varpi \in \Omega^n(*S)$ being independent of $\lambda$,
the multiplication by $f_j^{\pm1},\, (f_jf_k)^{\pm1}$ denoted by $\mathfrak{M}_{f_j^{\pm1}}\varphi\varpi$
and $\mathfrak{M}_{(f_jf_k)^{\pm1}}\varphi\varpi\;(j\ne k)$
correspond to the shift operators for the exponents
respectively:
\be
&&T_{\pm \varepsilon_j} :  \lambda\to \lambda\pm \varepsilon_j,\\
&&T_{\pm \varepsilon_j \pm\varepsilon_k} : \lambda \to \lambda \pm \varepsilon_j \pm \varepsilon_k.
\ee
In other words,
\be
&&T_{\pm\varepsilon_j}(\Phi\,\varphi\varpi) 
= \Phi\,\mathfrak{M}_{f_j^{\pm1}}\, T_{\pm \epsilon_j}\varphi \varpi=\Phi\, f_j^{\pm1}\,T_{\pm \epsilon_j}\varphi\varpi,\\
&&T_{\pm(\varepsilon_j+\varepsilon_k))}\,\Phi\,\varphi\varpi =\Phi\, \mathfrak{M}_{(f_jf_k)^{\pm1}}
\,T_{\pm(\varepsilon_j + \varepsilon_k)}\varphi\varpi = 
\Phi\,(f_jf_k)^{\pm1}T_{\pm (\varepsilon_j + \varepsilon_k)}\varphi\varpi.
\ee
\end{df}
These multiplications operate  on $H_\nabla^n (X,\Omega^\cdot(*S))$:
\be
 H_\nabla^n (X,\Omega^\cdot(*S))\overset{
 \mathfrak{M}_{{f_j}^{\pm1}}
 }
\longrightarrow H_\nabla^n (X,\Omega^\cdot(*S)).
\ee

Generally 
$T_{-\varepsilon_J} = \prod_{j\in J} T_{-\varepsilon_j}$
are well-defined for $\varepsilon_J = \sum_{j\in J}\varepsilon_j$.

The fundamental identity for the integral (1)
is given by
\be
T_{-\varepsilon_J}{\cal J}_\lambda(\varphi) = {\cal J}_{\lambda- \varepsilon_J}(T_{-\varepsilon_J}\varphi)
= {\cal J}_\lambda\bigl(\mathfrak{M}_{ f_J^{-1}} \,T_{-\varepsilon_J}\varphi\bigr).
\ee

The following Lemma plays a crucial role  in the sequel,
although it is almost obvious from the definition of the integrals (1).
It will be frequently  used without mentioning it.
%lm 6
\begin{lm}
Suppose $\varphi \varpi \in \Omega^n(*S)$ is cohomologously trivial, i.e.,
$
\varphi\,\varpi \sim 0.
$
Then also
$
\mathfrak{M}_{f_j^{\pm1}}\, T_{\pm \varepsilon_j}\varphi\,\varpi \sim 0\  
(with\, the \,same\, sign)
$
respectively, and vice versa.
\end{lm}
 
\begin{pf}
In fact, for $\varphi\varpi \in \Omega^n(*S)$ such that
$
{\cal J}_\lambda (\varphi) = 0,\ 
{\cal J}_{\lambda\pm\varepsilon_j} (\varphi) = 0,
$
and vice-versa. This implies Lemma 6. \qed
\end{pf}

%main thms $3
\section{Main Theorems (Case where $m = n+1)$}

\hspace{1cm}

   From this section until $\S5,$ we assume that $m = n+1$
so that $\dim H_\nabla^n(X, \Omega^\cdot(*S))$ 
$= 2^{n+1} - 1$.
$N = N_{n+1}$ will denote the set of indices $\{1,\ldots, n+1\}$. 
Denote by $\cal B$ the family of all admissible sets of 
indices which corresponds to the basis of
$H_\nabla^n(X, \Omega^\cdot(*S)).$

As for the basis of $H_n(X, {\cal L}^*),$
the following  Proposition follows from
the preceding results in [7], [8] and [10].

%prop8
\begin{prop}
As a basis of $H_n(X,{\cal L}^*),$
we can choose the homology classes of twisted n-cycles $\mathfrak{z}_J$
regularized from all relatively compact connected components 
of ${\bf R}^n - S$. Namely,
\be
&&{\rm(i)}\, \mathfrak{z}_{N} = \bigcap_{j\in N} \overline{D}_j,\\
&&{\rm(ii)}\, \mathfrak{z}_{J} = \bigcap_{j\in J} \overline{D}_j- \bigcup_{k\in J^c}\,D_k, \,(J \subset N, 1\le |J|\le n),\\
\ee
where $D_j,\,\overline{D}_j$ denotes the domain in ${\bf R}^n$
defined by $f_j(x) < 0, \, f_j(x) \le 0,$ respectively.
\end{prop}

For the proof and its background, see [7],[8],[10].

%re
\begin{re}\quad 
In Proposition 7, $\mathfrak{z}_\infty = {\bf R}^n -  \bigcup_{j\in N}\, D_j$,
the complement  of  $\bigcup_{j\in N}D_j,$ can also be regarded
as a twisted cycle. This cycle can be represented by a linear combination of
$\mathfrak{z}_J\,(J \in {\cal B})$ in homological sense. 
%We give it without proof   :
 
 In the case $n$ is odd,
\be
 \mathfrak{z}_{\infty} \sim
 - \sum_{J \in {\cal B}, |J|\le n} 
 \frac{\sin\,\pi\,\lambda_{J^c}}{\sin\,\pi\lambda_\infty}\,\mathfrak{z}_J .
 \ee
 In the case $n$ is even,
 \be
 \mathfrak{z}_{\infty} \sim
 - \sum_{J \in {\cal B}}
  \frac{\cos\,\pi\,\lambda_{J^c}}{\cos\,\pi\lambda_\infty}\,\mathfrak{z}_J .
\ee
Here $J^c$ denotes the complement $N - J$
and $\lambda_{J^c} = \sum_{j\in J^c}\,\lambda_j.$

The above formula can be proved by using the Cauchy  Integral Theorem
on a family of lines through a point in $\bigcap_{j\in N}\, D_j $.
\end{re}

\bigskip

We can introduce two kinds of systems of bases of $H_\nabla^n(X, \Omega^\cdot(*S))$
containing logarithmic forms
as in the following Proposition (see [7],[8],[11],[17]).

%prop8
\begin{prop}
The  $n$-forms
\be
 F_J = \frac{\varpi}{f_J} \quad(J \in {\cal B}).
\ee
represent a basis of $H_\nabla^n(X,\Omega^\cdot(*S)),$
where the product $\prod_{j\in J}\,f_j$ is abrreviated by $f_J$.

Likewise the following $n$-forms 
\be
 W_0(J)\varpi \,\,(J \in {\cal B})
\ee
represent another basis of $H_\nabla^n(X,\Omega^\cdot(*S))$,
where $W_0(J)\varpi$ are defined as
\be
&&W_0(j)\varpi = B(0\star j)\,F_j \quad ( j \in N ),\\
&&W_0(J)\varpi = - \sum_{\nu\in J}\,B\!
\left(
\begin{array}{ccc}
0  &\star   & \partial_\nu J  \\
 0 & \nu  & \partial_\nu  J\\
\end{array}
\right)\,F_{\partial_\nu J}
+ B(0\star J) \,F_J\quad(J \subset N, |J|\ge 2),
\ee
where $\partial_\nu J$ denotes the deletion of the element $\nu \in J$
from $J.$
\end{prop}

We shall call the former basis $F_J$ \lq\lq of  first kind", and the latter basis
$W_0(J)\varpi$ \lq\lq of  second kind".

%\vspace{.5cm}
%def9
\begin{df}
We can define the \lq\lq degree \rq\rq
of a rational differential $n$-form $\varphi\,\varpi \in \Omega^n(*S)
$ by the postulate 
\be
&&{\rm deg}\, x_\nu= {\rm deg}\, dx_\nu  = {\rm deg}\,\alpha_{j\nu} = {\rm deg}\, d\alpha_{j\nu} = 1,\\
&&{\rm deg}\, f_j = {\rm deg}\, df_j  = {\rm deg}\,\alpha_{j0}  = {\rm deg}\, d\alpha_{j0} = 2.
\ee

On the other hand, suppose that $\varphi\,\varpi$ is a linear combination of admissible $F_J:$
\be
\varphi\,\varpi = \sum_{J\in {\cal B}} [\varphi\varpi : F_J]\, F_J\quad([\varphi\varpi : F_J]
\,{\rm denotes\, the\, coefficient\,of}\,F_J).
\ee

We can define the \lq\lq weight\rq\rq of $\varphi\,\varpi$ as follows.
We say that ${\rm weight}\,(\varphi\varpi) \le q$ (or  ${\rm weight}\,(\varphi\varpi) \ge q$)
if all $F_J$ such that 
$ [\varphi\varpi : F_J] \ne 0$ satisfy $|J| \le q$ (or $|J| \ge q$).
We say that ${\rm weight}\,(\varphi\varpi) = q$ if all $F_J$ such that 
$ [\varphi\varpi : F_J] \ne 0$ satisfy $|J| = q$.
In particular, we have
 ${\rm weight}\,F_j = 1,\,{\rm weight}\,F_J =  p\,(|J| = p)$.
 \end{df}

The following Lemma immediately follows from the definition.
%lem10
\begin{lm}
For  an admissible $J,$ $F_J$ can be described as a linear combination of
$W_0(K)\varpi\,(K\subset J)$, {\rm i.e}., there exists a triangular matrix
$(\beta_{K,J})$ such that
\beq
B(0\star J)\,F_J = \sum_{K\subset J} \beta_{K,J} \,W_0(K)\varpi,
\eeq
where $\beta_{K,J}$ are uniquely determined by the relations
\be
&&\beta_{J,J}=1, \; \beta_{J,J\cup\{l\}}= \frac{
B\!
\left(
\begin{array}{ccc}
 0 &\star   & J  \\
 0 & l  &  J \\
\end{array}
\right)
}{B(0\star J)}\;\;(l\notin J),\\
&&
\beta_{K,J} = \sum_{l\in J-K} \beta_{K,K\cup\{l\}} \beta_{K\cup\{l\},J}.
\ee
More explicitly,
\be
&& \beta_{K,J} = \sum_{K=L_0 \subset L_1\subset \cdots \subset L_p = J}
\prod_{\nu = 1}^p\,\frac{
B\!
\left(
\begin{array}{ccc}
 0 &\star   &L_{\nu-1}   \\
 0 & l_\nu  &L_{\nu-1}   \\
\end{array}
\right)
}
{B(0\star L_{\nu-1})},
\ee
where $L_\nu = L_{\nu-1}\cup \{l_\nu\} \,(1\le \nu \le p)$ 
and the summation in the {\rm RHS} ranges  over the set of all sequences
$\{l_1,\l_2,\ldots,l_p\} = J - K$.

\end{lm}

We shall denote by $\tilde{\beta} = (\tilde{\beta}_{K,J})$
the upper triangular matrix consisting of the elements
\be
\tilde{\beta}_{K,J} = 
\left\{
\begin{array}{ccc}
\frac{\beta_{K,J}}{B(0\star J)} & (K \subset J), \\
 0 &\,{\rm (otherwise)}.\\
\end{array}
\right.
\ee

\bigskip
 
An arbitrary element $\varphi\,\varpi \in \Omega^n(*S)$ can be 
uniquely represented as a cohomological class in  
$H_\nabla^n(X,\Omega^\cdot(*S))$
by a linear combination of either of the above bases.

We are now in a position to give an explicit expression of the standard $n$-form $\varpi$
by the basis of second kind
in the following way.

%thm1
%\vspace{.5cm}
\bigskip

[{\bf Theorem 1}]\quad
We have
\beq
(2\lambda_\infty + n) \varpi \sim \sum_{\nu=1}^{n+1} (-1)^\nu
\sum_{J\subset N, |J| = \nu}\frac{\prod_{j\in J}\lambda_j}{\prod_{s=1}^{\nu-1} (\lambda_\infty + n - s)}
W_0(J)\,\varpi.
\eeq
{\it
Here $J$ ranges over the family of all {\rm(}unordered{\rm)} subsets of $N$ such that $1\le |J| \le n+1.$}

\bigskip

The second aim of this article is to give an explicit expression of $\nabla_B\varpi$
in terms of  invariant special 1-forms $\theta_J$\,(Theorem 2).
To do that, it is necessary to compute the contiguity relation for the shift
$\lambda \to \lambda \pm \varepsilon_j$. 
Theorem 2 plays a basic part in giving 
the covariant derivation $\nabla_B F_J$ for $F_J$
with an arbitrary admissible $J$.

%\vspace{.5cm}
\bigskip
%thm2
[{\bf Theorem 2}]\quad
We have
{\it
\beq
\nabla_B(\varpi) \sim
 \sum_{p=1}^{n+1} 
\sum_{J\subset N, |J| = p}\frac{\prod_{j\in J}\lambda_j}{\prod_{s=1}^{p-1} (\lambda_\infty + n - s)}
\,\theta_J\,W_0(J) \,\varpi,
\eeq
where $J \subset N$ move over the family of sets of indices 
such that $|J| = p$. $\theta_J$ represents one degree form 
defined by 
\beq
&&\theta_j = - \frac{1}{2}\,d\log \,B(0\star j)
 = - \frac{1}{2}
  \,d\log\,r_j^2,
  \quad(J = \{j\},\,{\rm i.e.,}\, |J| = 1)
  \nonumber\\
&&\theta_{jk} = \frac{1}{2}\,d\log\,B(0jk) = \frac{1}{2}\, d\log\rho_{jk}^2,
\quad(J = \{j,k\},\,{\rm i.e.,}\,|J| = 2)\nonumber\\
&&\theta_J = (-1)^p \sum_{j,k\in J, j<k }\,\frac{1}{2}\, d\log B(0jk)\cdot \nonumber\\
&& \sum_{\{\mu_1,\ldots,\mu_{p-2}\}= \partial_j\partial_k J  }\prod_{s=1}^{p-2}
 \frac{
B\!
\left(
\begin{array}{cccccccc}
 0 &\star   &\mu_{s-1}&\cdots&\mu_1&j&k   \\
 0 & \mu_s  &\mu_{s-1}& \cdots&\mu_1&j&k  
\end{array}
\right)
}
{B(0\mu_s\,\mu_{s-1}\cdots\mu_1\,j\,k)},\quad(|J|= p \ge 2)\nonumber\\
\eeq
where $\{\mu_1,\ldots,\mu_{p-2}\}$ move over the family of all ordered
sequences consisting of  $p-2$ different elements of $\partial_j\partial_k J$.
}
%endthm2

%\vspace{.5cm}
\bigskip

[{\rm \bf Question}]

Generally the one degree forms $\theta_J$ in {\rm Theorem 2} seem 
\lq\lq  generically  logarithmic" along the singularity, i.e., logarithmic  in the sense of K.Saito {\rm [25]}.

\begin{re}
Theorem 1 and Theorem 2 still hold true for $1\le m \le n$ by taking 
$\lambda_j = 0\,(m+1\le j\le n+1)$.
The origin of Theorem 2  goes back to
the variational formula for the volume of a spherical  geodesic simplex
due to L.Schl\"afli  (See [1],[3],[16],[20], [26],[27].
See also [14] related to the bellows conjecture.)
The point is that the coefficients involved in the differential $\nabla_B$
can be described in terms of Cayley-Menger invariants relative to
the isometry group ({\rm [24],[28]}).
\end{re}

\bigskip

Under the assumptions $({\cal H}1)$ and $({\cal H}2),$
 $f_j$ are can be described  in the following  normalized 
form by choosing  suitable coordinates $(x)$:
\beq
&&f_{n+1}(x) = Q(x) + \alpha_{n+1\,0},\nonumber\\
&& f_j(x) = Q(x)+2\sum_{\nu=1}^{n+1- j}\alpha_{j\nu}x_\nu+ \alpha_{j0}\quad(1\le j\le n),\nonumber\\
\eeq
for $\alpha_{j\nu},\,\alpha_{j0}\in {\bf R}$.
Then
\be
&&r_j^2= \sum_{\nu = 1}^{n+1-j}\alpha_{j\,\nu}^2-\alpha_{j\,0},\quad(1\le j\le n)\\
&&r_{n+1}^2 = -\alpha_{n+1\,0}.\\
\ee

We may take $\alpha_{n1}>0, \alpha_{n-1\,2}>0, \ldots, \alpha_{1\,n}>0.$

%lm11
\begin{lm}
Under the condition $({\cal H}2),$
we have
\be
&& \alpha_{n+1\,0} = - \frac{1}{2}\,B(0\star n+1),\quad \alpha_{n1}^2 = \frac{1}{2}B(0\ n\,n+1),
\quad \alpha_{n1}=\rho_{n\,n+1},\\
&&\alpha_{j\,n+1-j}= \sqrt{\frac{B(0\, j\ldots\, n+1)}{-2B(0\,j+1\,\ldots\, n+1 )}}\quad(1\le j\le n),\\
&& \prod_{\nu = 1}^{p}\alpha_{n - \nu + 1, \nu}^2 = (-1)^{p-1} \frac{B(0 \,n-p+1\ldots  n+1) }{2^p}\quad(1\le p\le n) ,\\
&&\sum_{\nu=1}^ {n+1-k} \alpha_{j\nu}\alpha_{k\nu}=\frac{1}{2} 
B
\!\left(
\begin{array}{ccc}
 0 & j & n+1 \\
 0 &  k&  n+1  
\end{array}
\right) = \frac{1}{2}\,(\rho_{j\,n+1}^2 + \rho_{k\,n+1}^2 - \rho_{j\,k}^2)\\
&&\qquad \qquad \qquad \qquad \qquad \qquad \qquad \qquad \qquad \qquad 
(1\le j\le k\le n),\\
&&\alpha_{j2}\alpha_{n-1\,2}= - \frac{
B\!
\left(
\begin{array}{cccc}
 0&j &n   & n+1  \\
 0&n-1 & n  & n+1  \\
\end{array}
\right)
}{4\rho_{n\,n+1}^2}\quad(1\le j\le n-1),\\
&&\alpha_{j0}-\alpha_{n+1\,0}= \rho_{j\,n+1}^2+r_{n+1}^2-r_j^2 = 
B
\!\left(
\begin{array}{ccc}
 0 & j  & n+1 \\
 0 &  \star &  n+1 
\end{array}
\right),\\
&&\alpha_{j0}+\alpha_{n+1\,0}= - 
B\!\left(
\begin{array}{ccc}
 0 & \star  & j  \\
 0 &  \star &  n+1  
\end{array}
\right).\\
\ee
\end{lm}

Note that $B(0 12\ldots p)=0$ holds if and only if 
the $p-1$ dimensional volume of the $p-1$-simplex  $\Delta O_{1}O_2\ldots O_{p}$ equal to $0,$
i.e., degenerates into a $p-2$ dimensional affine subspace.

\bigskip

%notation1
\begin{N}

We define the logarithmic forms for an each ordered subset
$J = \{j_1,\ldots, j_p\}\ (1\le p \le n)$ by
\be
&&e_j=\frac{df_j}{f_j}, \\
&&e_J=\frac{df_{j_1}}{f_{j_1}}\wedge \cdots \wedge \frac{df_{j_p}}{f_{j_p}},
\ee
and the function $W$
such that the alternating sum is expressed as
\beq
e_{2\ldots n+1} - e_{13\ldots n+1} + \cdots  + (-1)^{n}\,e_{12\ldots \, n}= 2^n W\varpi\quad(W\in \Omega^0(*(S)).
\eeq
\end{N}

The following identity can be obtained by an elementary
calculation.
%lem12
\begin{lm}
\beq
&&W=(-1)^{\frac{n(n-1)}{2} + 1}
\frac{1}{\sqrt{(-1)^{n+1}\, 2^n\, B(0\,N)}}\,W_0(N).
\eeq

\end{lm}

%lm13
\begin{lm}
In $H^n_{\nabla}(X,\Omega^\cdot(*S))$
%as cohomology class,  
we have the identities
\beq
e_{12\ldots\,\widehat{\nu}\,\ldots \, n+1}\sim 
(-1)^{\nu-1}\,\frac{2^n\lambda_\nu}{\lambda_\infty}W\varpi\quad(1\le \nu\le n+1).
\eeq
\end{lm}

\begin{pf}
Indeed, by the Stokes formula it holds
\be
\nabla(e_{3\ldots n+1}) = \lambda_1e_{13\ldots n+1}+\lambda_2 e_{23\ldots n+1}\sim 0.\\
\ee
More generally, for any ordered $K \subset N$ such that $|K| = n-1$ and 
two different $\mu,\nu \in K^c,$
\be
\lambda_\mu e_{\mu K} + \lambda_\nu e_{\nu K} \sim 0.
\ee

These relations plus {\rm(8)} show {\rm Lemma 13}.\qed
\end{pf}

\begin{re}
 In course of proving Theorem 1 and 2, we will see that the following three properties concerning 
Cayley-Menger determinants play an essential role:
\it{the cofactor expansion, Sylvester(Jacobi) formula
and  Pl\"ucker's relation}.
\end{re}

\bigskip

%sect4
\section{Hierarchy of Contiguity Relations}

%\hspace{1cm}
\bigskip

To prove Theorems 1 and 2,
we make use of representation of  the shift operators
in terms of the basis stated in Proposition 8.
These are given in Propositions 21, 23 and 24.

Let  $\boldsymbol \alpha$ be the $n \times n $ matrix
with $j,\nu$ th elements $\alpha_{j\nu}\,(1\le j\le n,\,1\le \nu \le n).$
The Gram matrix 
 associated with $\boldsymbol \alpha$ is equal to
 the $n\times n$ matrix
\be
{\boldsymbol \alpha}\cdot {^t}\!{\boldsymbol \alpha} = \bigl (\frac{1}{2}\,
B\!
\left(
\begin{array}{ccc}
 0 & j  & n+1  \\
 0 & k  & n+1  \\   
\end{array}
\right)\bigr)_{
1\le j\le n, 1\le k\le n
}.
\ee
The minor determinant consisting of $j_1,\ldots,j_p$th rows
and $\nu_1,\ldots,\nu_p$th columns of ${\boldsymbol \alpha}$ is denoted by
\be
{\boldsymbol \alpha}\!
\left(
\begin{array}{ccc}
 j_1 &\ldots   &j_p   \\
 \nu_1 &\ldots   &\nu_p   \\
\end{array}
\right).
\ee

From Lemma 11 and Gram determinant formula, it follows that
%llm14
\begin{lm}
For $1\le p\le n,$
\be
&&\alpha_{j\,n-p+1} = (-1)^{n-p} 
\frac{
B\!
\left(
\begin{array}{ccccc}
 0 &j   & p+1&\ldots &n+1 \\
 0 & p  & p+1&\ldots&n+1  
\end{array}
\right)
}
{
\sqrt{-2 B(0\, p\ldots n+1)\,B(0\,p+1\ldots n+1)}
}\quad(1\le j\le p),\\
&&
{\boldsymbol\alpha}\!
\left(
\begin{array}{ccccc}
p   & p+1&\ldots\widehat{\;j\;}\ldots &n \\
 1  & 2&\ldots&n-p 
\end{array}
\right) = (-1)^{\frac{(n-p)(n-p-1)}{2}+ (n- j)}\\
&&
{}\qquad \qquad \cdot
\frac{
B\!
\left(
\begin{array}{ccccc}
 0 &p  & p+1&\ldots\widehat{\;j\;}\ldots &n+1 \\
 0 & j & p+1&\ldots\widehat{\;j\;}\ldots&n+1  
\end{array}
\right)}
{\sqrt{2^{n-p}(-1)^{n+1-p}\,B(0\,p+1\ldots n+1)}}\quad(p+1\le j\le n),\\
&&{\boldsymbol\alpha}\!
\left(
\begin{array}{ccccc}
 p+1&p+2&\ldots &n \\
 1  & 2&\ldots&n-p 
\end{array}
\right)
= (-1)^{\frac{(n-p)(n-p-1)}{2}}\,\\
&&
{}\qquad \qquad \cdot\sqrt{
\frac{(-1)^{n+1-p}\,B(0\,p+1\,\ldots\, n+1)}
{2^{n-p}
}}\quad(j = p).
\ee
\end{lm}

%lem15
\begin{lm}
Fix $p\; (1\le p\le n)$.
Then the following equality holds:
\be
&&\sum_{\nu=1}^p (-1)^{\nu-1} \alpha_{\nu\,n-p+1}df_1\wedge \ldots\widehat{\;df_\nu\;}\ldots
\wedge df_p\wedge\ldots\wedge df_{n+1} \\
&&{}\qquad = 2^n (-1)^{\frac{n(n+1)}{2}+1} x_{n-p+1} \prod_{\mu=1}^n \alpha_{n-\mu+1\,\mu}\,\varpi.
\ee
In particular, when $p=1$ or $n$ we have
\be
&&\alpha_{1n}\,df_2\wedge \cdots \wedge df_{n+1} 
= 2^n (-1)^{\frac{n(n+1)}{2}+ 1}x_n \prod_{\mu=1}^n \alpha_{n-\mu+1\,\mu}\varpi,\\
&& \sum_{\nu=1}^n (-1)^{\nu-1} \alpha_{\nu\,1}df_1\wedge\cdots\wedge df_{\nu-1}\wedge df_{\nu+1}\wedge\cdots\wedge
df_{n+1} \\
&&{}\qquad = 2^n  (-1)^{\frac{n(n+1)}{2}+1}x_1 \prod_{\mu=1}^n \alpha_{n-\mu+1\,\mu}\varpi,
\ee
respectively.
\end{lm}

From Lemma 13 and  Lemma 15, the following holds.

%co16
\begin{co}
For $J = \{1,2,\ldots, p\}\;(1\le p\le n),$
\beq
&&\frac{x_{n-p+1}}{f_{p+1}\,\ldots\,f_{n+1}}\varpi
 \sim
\frac{(-1)^n}{\sqrt{2^n (-1)^{n+1} B(0 N)}}\cdot \nonumber\\
&&
{}\qquad \frac{1}{\prod_{\nu=1}^n \alpha_{n-\nu+1\,\,\nu}}\;
\sum_{\nu=1}^p \alpha_{\nu\,n-p+1}\frac{\lambda_\nu}{\lambda_\infty+p-1}
\prod_{j\in \partial_{\nu}J}\mathfrak{M}_{f_j}(W_0(N)\varpi)\nonumber.\\
\eeq
\end{co}

%pf
\begin{pf}
Indeed, {\rm Lemma 15} implies the following identity
\be
2^n (-1)^{\frac{n(n+1)}{2}}\prod_{\mu=1}^n \alpha_{n-\mu+1\,\mu}
\frac{x_{n-p+1}}{f_{p+1}\cdots f_{n+1}} \,\varpi = 
\sum_{\nu=1}^p
(-1)^{\nu-1} \alpha_{\nu\,n-p+1}
 \prod_{\mu\in \partial_\nu J}\,\mathfrak{M}_{f_\mu}
 \, (e_{1\ldots\widehat{\nu}\ldots n+1}).
\ee
We can apply {\rm Lemma 13} to this identity and obtain {\rm (11)}.
\qed
\end{pf}

On the other hand, as an immediate consequence
of the definition of $W_0(N)\varpi$, we have

%lem17
\begin{lm}
For $J \subset N,$
\be
&&\mathfrak{M}_{f_J}(W_0(N)\varpi)
= -
\sum_{j\in J}
B\!\left(
\begin{array}{ccc}
0  &\star   &  \partial_j N \\
0  & j  &\partial _j N  \\
\end{array}
\right)
\mathfrak{M}_{f_j}F_{J^c}\\
&&
{}\qquad \qquad \qquad \quad -
\sum_{k\in J^c}
B\!\left(
\begin{array}{ccc}
0  &\star   &  \partial_k N \\
0  & k  &\partial _k N  \\
\end{array}
\right)
F_{\partial_k J^c}
+
B(0 \star N)F_{J^c}.
\ee
\end{lm}

Denote by  $\overline{f}_j$  the homogeneous part of ${f}_j:
\overline{f}_j = f_j - \alpha_{j\,0}\quad(1\le j\le n+1).
$
Then
%lem18
\begin{lm}
For $1\le p\le n,$ the following equality holds
\be
&&(\overline{f}_p - \overline{f}_{n+1})\, \prod_{\nu=1}^{n-p}\alpha_{n-\nu+1\,\nu}
= 2 x_{n-p+1}\prod_{\nu=1}^{n-p+1}\alpha_{n-\nu+1\,\nu}\\
&&{}\ + (-1)^{\frac{(n-p)(n-p-1)}{2}}
\sum_{\mu=p+1}^{n} (\overline{f}_\mu - \overline{f}_{n+1}) (-1)^{\mu-p-1}
{\boldsymbol \alpha}\!\!
\left(
\begin{array}{ccccc}
 p & p+1  & \ldots\widehat{\,\mu\,}\ldots& n  \\
 1 & 2  & \ldots& n-p 
\end{array}
\right).
\ee
\end{lm}

\begin{pf}
Indeed, from (7)
%$\reflectbox{\ddots}$
\be
\left|
\begin{array}{cccccccc}
 -x_{n-p+1} & 0  &  0& \ldots&\ldots&0&1 \\
 -\overline{f}_p + \overline{f}_{n+1} & 2\alpha_{p\,1}  & \alpha_{p\,2}&\ldots&\ldots&2\alpha_{p\,n-p}&2\alpha_{p\,n-p+1}  \\
  -\overline{f}_{p+1} + \overline{f}_{n+1} & 2\alpha_{p+1\,1}  &  \alpha_{p+1\,2} &\ldots&\ldots&2\alpha_{p+1\,n-p}&0\\
  %
 % \vdots&\vdots&\vdots&& &\mbox{\reflectbox{$\ddots$}}&\vdots\\
  %
  \vdots&\vdots&\vdots&&
  \mbox{\reflectbox{$\ddots$}}&&\vdots\\
  -\overline{f}_{n-1}+ \overline{f}_{n+1}&2\alpha_{n-1\,1}&2\alpha_{n-1\,2}&0&\ldots&\ldots&0\\
  - \overline{f}_n + \overline{f}_{n+1}& 2\alpha_{n\,1}&0 &\ldots&\ldots&\ldots&0\\
\end{array}
\right|= 0.
\ee
We can deduce {\rm Lemma 18} by expanding the cofactors along the first column.

\qed
\end{pf}

Lemma 14 and Lemma 18 imply the following.
 
%lem19
\begin{lm}
For $1\le p \le n,$
\beq
&&\frac{\overline{f}_p - \overline{f}_{n+1}}{f_{p+1}\ldots f_{n+1}}\varpi
=
\sqrt{-\frac{B(0\,p\ldots\,n+1)}{2\,B(0\,p+1\;\ldots\;n+1)}}
\frac{2x_{n-p+1}}{f_{p+1}\ldots f_{n+1}}\,\varpi\nonumber\\
&&{}\ + \sum_{\nu=p+1}^n 
\frac{
B\!
\left(
\begin{array}{cccccc}
 0 &p   & p+1&\ldots\widehat{\nu}\ldots &n+1  \\
  0&\nu   &p+1&\ldots\widehat{\nu}\ldots & n+1  
\end{array}
\right)
}
{B(0\,p+1\,\ldots\,n+1)}
\frac{\overline{f}_\nu - \overline{f}_{n+1}}{f_{p+1}\ldots f_{n+1}}\varpi.
\eeq
\end{lm}

Due to Jacobi identities of determinants and the cofactor expansion, we 
see the following identities:
\be
&&\alpha_{p 0} - \alpha_{n+1\,0} - \sum_{\nu=p+1}^n 
\frac{(\alpha_{\nu0}\,- \alpha_{n+1\,0})\,
B\!
\left(
\begin{array}{cccccccc}
 0 &p   &p+1&\ldots&\widehat{\nu}&\ldots&n+1   \\
 0 & \nu  &p+1&\ldots&\widehat{\nu}&\ldots&n+1   
\end{array}
\right)
}{B(0\,p+1\ldots n+1)}\\
&&
{}\qquad \qquad \qquad = - \frac{B\!
\left(
\begin{array}{cccccccc}
 0 &\star   &p+1&\ldots&n+1   \\
 0 & p  &p+1&\ldots&n+1   
\end{array}
\right)}{B(0\,p+1\ldots n+1)},\\
&& 1 - \sum_{\nu= p+1}^n \frac{B\!
\left(
\begin{array}{cccccccc}
 0 &p   &p+1&\ldots&\widehat{\nu}&\ldots&n+1   \\
 0 & \nu  &p+1&\ldots&\widehat{\nu}&\ldots&n+1   
\end{array}
\right)
}{B(0\,p+1\ldots n+1)}\\
&&
{}\qquad \qquad \qquad = \frac{B\left(
\begin{array}{cccccccc}
 0 &p   &p+1&\ldots&n+1   \\
 0 &n+1 &p+1&\ldots&n+1   
\end{array}
\right)}{B(0\,p+1\ldots n+1)}.
\ee

Seeing that
\be
&&\overline{f}_j - \overline{f}_{n+1} = f_j - f_{n+1} + \alpha_{n+1\,0} - \alpha_{j\,0}\quad(1\le j\le n),
\ee

(12) can be rewritten in the following form
%(21p)
\beq
&&\frac{f_p}{f_{p+1}\cdots f_{n+1}} \,\varpi =
\frac{2\,x_{n-p+1}}
{f_{p+1}\,\cdots\, f_{n+1}}
\sqrt{-\frac{B(0p\ldots n+1)}{2\,B(0\,p+1\ldots \,n+1)}
}\,\,\varpi\nonumber\\
&&{}\quad - 
\frac{B\!
\left(
\begin{array}{ccccccc}
 0 & \star  &p+1&\ldots& n+1   \\
 0 &p   &p+1&\ldots&n+1   \\
\end{array}
\right)
}{B(0\,p+1\,\ldots\,n+1)}
\frac{\varpi}{f_{p+1}\,\cdots f_{n+1}}\nonumber\\
&&{}\quad +
\sum_{\nu= p+1}^{n+1}
\frac{B\!
\left(
\begin{array}{cccccccc}
 0 & p  &p+1&\ldots&\widehat{\nu}&\ldots& n+1   \\
 0 &\nu   &p+1&\ldots&\widehat{\nu}&\ldots&n+1   \\
\end{array}
\right)
}{B(0\,p+1\;\ldots\,n+1)}\,\frac{\varpi}{f_{p+1}\cdots\widehat{f}_\nu\cdots f_{n+1}}.\nonumber\\
\eeq

 (11) and (13) leads us to the following.

%lm20
\begin{lm}
For $1\le p\le n,$
\beq
&&\mathfrak{M}_{f_p }F_{p+1\,\ldots\,n+1} \sim \frac{1}{B(0\,N) B(0\,p+1\,\ldots n+1)}\cdot \nonumber\\
&&{}\quad \sum_{\nu=1}^p 
B\!
\left(
\begin{array}{ccccc}
 0 & \nu  & p+1\ldots n+1 \\
 0 & p  &p+1\ldots n+1 \\
\end{array}
\right)\frac{\lambda_\nu}{\lambda_\infty + p-1} 
\mathfrak{M}_{f_1\ldots\widehat{\;f_\nu\;}\ldots f_p} (W_0(N)\varpi)\nonumber\\
&&{}\quad - \frac{
B\!
\left(
\begin{array}{ccccc}
 0 & \star  & p+1\ldots n+1 \\
 0 & p  &p+1\ldots n+1 \\
\end{array}
\right)}{B(0\, p+1\ldots n+1)} F_{p+1\,\ldots\, n+1}\nonumber\\
&&{}\quad + \sum_{\nu = p+1}^{n+1} 
\frac{
B\!
\left(
\begin{array}{ccccc}
 0 & \nu  & p+1\ldots\widehat{\,\nu\,}\ldots n+1 \\
 0 & p  &p+1\ldots\widehat{\,\nu\,}\ldots n+1 \\
\end{array}
\right)}
{B(0\,p+1\ldots n+1)}\,F_{p+1\ldots\widehat{\nu}\ldots n+1},
\eeq
where the last term in the {\rm RHS} of {\rm (14)} reduces to $\varpi$ for $p= n.$
\end{lm}

%Prop21
%\vspace{.5cm}

\bigskip

\begin{prop}
For $J \subset N\, (|J| = p,\, 2\le p\le n+1)$ and $j\in J,$
\beq
&&\mathfrak{M}_{f_j}\,F_J = F_{\partial_j J}\quad(j\in J),\\
&&\mathfrak{M}_{f_j}\,F_j = \varpi.
\eeq
On the contrary,
for $j\in J^c,$
the following identity holds{\rm :}
\beq
&&\mathfrak{M}_{f_j} \, F_J \sim
\{\sum_{\nu \in J^c}\frac{B\!
\left(
\begin{array}{cccc}
 0 &\nu   &J   \\
 0 & j  & J   
\end{array}
\right)}{B(0\,N)\,B(0\, J)}
\frac{\lambda_\nu}{\lambda_\infty + n - p}
\prod_{l\in \partial_\nu J^c }\,\mathfrak{M}_{f_l}(W_0(N)\varpi)\}\nonumber\\
&&
{}\qquad \qquad -\frac{B\!
\left(
\begin{array}{cccc}
 0 &\star   &J  \\
 0 & j  & J   
\end{array}
\right)}{B(0 J)}\, F_J
+ \sum_{\nu\in J}\,\frac{B\!
\left(
\begin{array}{cccc}
 0 &\nu  &\partial_\nu J \\
 0 & j  & \partial_\nu J   
\end{array}
\right)}{B(0 J)}\,F_{\partial_{\nu} J}.
\eeq
Here $F_\emptyset$ means $\varpi.$
\end{prop}

\bigskip

%\vspace{.5cm}
%pf of Prop21
[{\bf proof of Proposition 21 }]

%\begin{pf}
Because of symmetry we have (17) from (14) replacing $f_p$ with $f_j$ 
and the set $\{p+1,\ldots,n+1\}$ 
with $J\, (|J| = p)$ such that $j\in J^c$ repectively.
Hence Proposition 21 holds. 
\qed
%\end{pf}

\bigskip

Now Proposition 21 and Lemma17 give a recurrence relation
for $\mathfrak{M}_{f_J}W_0(N)\varpi$
in the following way.

%lm22
\begin{lm}
For $J \subset N,\,|J| = p\le n,$
\beq
&&\mathfrak{M}_{f_J}W_0(N)\varpi \sim - \sum_{\mu \in J}
\frac{
B\!
\left(
\begin{array}{ccc}
 0 &\star   &J^c   \\
 0 &\mu   &J^c   
\end{array}
\right)
}{B(0\,J^c)} \frac{\lambda_\mu}{\lambda_\infty + p-1}\,\mathfrak{M}_{f_{\partial_\mu\! J}} W_0(N)\varpi
\nonumber\\
&&{}\qquad \qquad \qquad + \frac{B(0\,N)}{B(0\,J^c)}\,W_0(J^c)\varpi + 
(-1)^n \, \delta_{p\,n}\, B(0\,N)\,\varpi,
\eeq
where $\delta_{p\,n}$ means  Kroneckers's delta.
\end{lm}

%proof
\begin{pf}
In fact, we can derive (18) from (17) and {\rm Lemma 17} by making  use of the following identities 
(Sylvester determinant formula plus cofactor expansion):
\be
&&\sum_{j\in J}\, B\!
\left(
\begin{array}{ccc}
 0 &\star   &\partial_j N   \\
  0&j   &\partial_j N    \\
\end{array}
\right)
\,B\!
\left(
\begin{array}{ccc}
 0 &k   &J^c   \\
  0&j   &J^c   \\
\end{array}
\right)
= B(0\,N)\,B\!
\left(
\begin{array}{ccc}
 0 &\star   &J^c   \\
  0&k   &J^c   \\ 
\end{array}
\right)\quad(k\in J),\\
&&\sum_{j\in J}\, B\!
\left(
\begin{array}{ccc}
 0 &\star   &\partial_j N   \\
  0&j   &\partial_j N    \\
\end{array}
\right)
\,B\!
\left(
\begin{array}{ccc}
 0 &\star  &J^c   \\
  0&j   &J^c   \\
\end{array}
\right) + B(0\star N)\,B(0\,J^c)
= B(0\,N)\,B(0\star J^c),\\
&&\sum_{j\in J}\, B\!
\left(
\begin{array}{ccc}
 0 &\star   &\partial_j N   \\
  0&j   &\partial_j N    \\
\end{array}
\right)
\,B\!
\left(
\begin{array}{ccc}
 0 &\nu   &\partial_\nu J^c   \\
  0&j   &\partial_\nu J^c   \\
\end{array}
\right)
+ B\!
\left(
\begin{array}{ccc}
 0 &\star   &\partial_\nu N   \\
  0&\nu   &\partial_\nu N   \\
\end{array}
\right)\,B(0\,J^c)\\
&&
{}\qquad \qquad = B(0\,N)\,
 B\!
\left(
\begin{array}{ccc}
 0 &\star   &\partial_\nu J^c   \\
  0&\nu   &\partial_\nu J^c   \\
\end{array}
\right)\quad(\nu \in J^c).
\ee
\qed
%\hspace{1cm}$\Box$
\end{pf}

As a result of Lemma 22, an explicit formula for $\mathfrak{M}_{f_J}(W_0(N)\varpi)$
is deduced as follows:

%Proposition23
%\vspace{.5cm}

\begin{prop}
For $J \subset N\,(|J| = p, 0\le p\le n),$
\beq
&&\mathfrak{M}_{f_J}\,\bigl(W_0(N) \varpi \bigr)  = 
 \prod_{j\in J} \mathfrak{M}_{f_j} \bigl(W_0(N) \varpi \bigr)\nonumber\\
&& \sim 
 (-1)^n \,\delta_{p\,n}\,B(0\,N)\,\varpi  + \sum_{q = 0}^p (-1)^q 
 \sum_{\{\mu_1,\ldots,\mu_q\}\subset J}
\frac{\prod_{\nu=1}^q\lambda_{\mu_\nu}}{\prod_{\nu=1}^q (\lambda_\infty + p - \nu)}
\cdot
\nonumber\\
&&
{}\quad \prod_{\nu=1}^q
\frac{B\!
\left(
\begin{array}{cccccc}
 0 &\star  &\mu_{\nu-1}\ldots &\nu_1& J^c  \\
 0 &\mu_\nu  &\mu_{\nu-1}\ldots&\mu_1& J^c  \\
\end{array}
\right)
}
{B(0\mu_1\ldots\mu_{\nu} \;J^c)}\, \frac{B(0\,N)}{B(0\,J^c)}
\;W_0(\mu_1\ldots\mu_q J^c)\varpi,
\eeq
where $\{\mu_1,\ldots,\mu_q\}$ represents an ordered subset of $J$
consisting of all different elements of $J$ and  it means the empty set
if $q=0$.
\end{prop}

\bigskip

%\vspace{.5cm}
%pf of Prop23
[{\bf proof of Proposition 23 }]

Proposition 23 can be proved by the
repeated use of  the recurrence relations (18)
concerning the size $|J|$.
\qed

\bigskip

%\vspace{.5cm}
(14) can be rewritten through {\rm Proposition 23} by a linear
combination of the terms
$F_{p+1\ldots n+1}, \,F_{p+1\ldots\widehat{\nu}\ldots n+1}$
and $W_0(K)\varpi\ (J^c \subset K)
$.

In fact, for $1\le p\le n,$ we have
\be
\mathfrak{M}_{f_p}F_{p+1\ldots n+1}\sim
U_0 + 
\sum_{\nu=1}^ p \sum_{J\subset N, |J| = \nu} \frac{\prod_{j\in J} \lambda_j}{\prod_{s=1}^\nu(\lambda_\infty+p-s)} \;U_J.
\ee
Both of $U_0, \,U_J$ do not depend on $\lambda$.
$U_0$ is expressed as
\be
&&U_0 = - \frac{
B\!
\left(
\begin{array}{ccccc}
 0 & \star  & p+1\ldots n+1 \\
 0 & p  &p+1\ldots n+1 \\
\end{array}
\right)}{B(0\, p+1\ldots n+1)} F_{p+1\,\ldots\, n+1}\\
&&{}\qquad + \sum_{\nu = p+1}^{n+1} 
\frac{
B\!
\left(
\begin{array}{ccccc}
 0 & \nu  & p+1\ldots\widehat{\,\nu\,}\ldots n+1 \\
 0 & p  &p+1\ldots\widehat{\,\nu\,}\ldots n+1 \\
\end{array}
\right)}
{B(0\,p+1\ldots n+1)}\,F_{p+1\ldots\widehat{\nu}\ldots n+1}.
\ee
 The weights appear 
 with the condition
\be
\mbox{ weight}\, U_0  =  n- p ,\quad \mbox{weight}\, U_J \ge  n + \nu - p \ \ (|J| = \nu ).
\ee
.

The simplest case $p=1$ is written as
\be
&&\mathfrak{M}_{f_1}F_{2\ldots\,n+1} \sim  U_0 + \frac{\lambda_1}{\lambda_\infty} \;U_1,\\
&& U_0 = 
- \frac{B\!
\left(
\begin{array}{ccc}
 0 &\star   &2\ldots n+1   \\
  0&1   & 2\ldots n+1  \\
\end{array}
\right)
}{B(02\ldots n+1)}\,F_{2\ldots n+1}
\,+\,
\sum_{\nu= 2}^{n+1} 
\frac{B\!
\left(
\begin{array}{ccc}
 0 & \nu  &2 \ldots \widehat{\nu}\ldots n+1 \\
  0& 1  & 2\ldots\widehat{\nu}\ldots n+1  \\
\end{array}
\right)
}{B(02\ldots n+1)} \;F_{2\ldots\widehat{\nu}\ldots n+1},\\
&&U_1 = \frac{1}{B(0 2\ldots n+1))}\,W_0(N)\varpi.
\ee

%re
\begin{re}
In (19),
the weight in the LHS is at least $\,n-p\,$
 since
${\rm weight}\, \mathfrak{M}_{f_J}\bigl(W_0(N) \varpi \bigr) \ge   n - p,$ 
and
in the RHS ${\rm weight}\, W_0(\mu_1\ldots \mu_s\,J^c)\,\varpi\ge   n  - p + s$.
Hence in (17), ${\rm weight}\,(\mathfrak{M}_{f_j}F_J) \ge p - 1$ for every $J,$
if $p \ge 2$.
\end{re}

In the simplest case where
$p=1,$ we have
\be
&&\mathfrak{M}_{f_1}\bigl(W_0(N)\varpi\bigr)
= - \frac{B\!
\left(
\begin{array}{ccccc}
 0&\star &2   &\cdots& n+1   \\
  0&1   &2&\cdots& n+1   \\
\end{array}
\right)
}
{B(02\cdots\,n+1)}\;
\frac{\lambda_1}{\lambda_\infty}\;W_0(12\ldots n+1)\varpi\\
&&{}\qquad \qquad \qquad \quad \ + \frac{B(012\ldots\,n+1)}{B(02\ldots\,n+1)} 
\, W_0(23\ldots\,n+1)\varpi.
\ee

On the other hand, in the case where $p = n,$ we have
%prop24
\begin{prop}
For $j,k\in N,\; j\ne k$
\beq
&&\mathfrak{M}_{f_j - f_k} F_{k} \sim
\frac{B\!
\left(
\begin{array}{cccccc}
 0 &j & k \\
 0 &\star & k  \\
\end{array}
\right)
}
{B(0\star k)
}\,W_0(k)\varpi
 + \sum_{\nu=1}^n (-1)^\nu
 \sum_{M = \{\mu_1,\ldots,\mu_\nu \}\subset \partial_k N }\frac{\prod_{s=1}^\nu \lambda_{\mu_s}}{\prod_{s=1}^\nu(\lambda_\infty + n - s)}\nonumber\\
&&
{}\qquad \quad  \cdot\frac{B\!
\left(
\begin{array}{cccccc}
 0 &j& k  \\
 0 &\mu_1 &k  \\
\end{array}
\right)
}
{B(0\mu_1 k)}
 \prod_{s=2}^\nu\frac{B\!
\left(
\begin{array}{cccccc}
 0 &\star&\mu_{s-1}&\ldots&\mu_{1} &k  \\
 0 &\mu_s&\mu_{s-1}&\ldots&\mu_{1} &k  \\
\end{array}
\right)
}
{B(0\mu_s \mu_{s-1}\ldots \mu_1 k)}
W_0(\mu_1\ldots \mu_\nu k)\varpi,
\eeq
where $M = (\mu_1,\ldots,\mu_\nu)$
ranges over the family of ordered sets such that $M\subset \partial_k N$ and 
$|M| = \nu.$
\end{prop}
%end prop

\bigskip

%pf of Prop25
%\vspace{.5cm}
[{\bf proof of Proposition  24}]

In the case where $p=n,$
(19) is equivalent to
\beq
&&\mathfrak{M}_{f_n - f_{n+1}}\,F_{n+1}\sim
- \sum_{\nu=1}^{n} \frac{
B\!\left(
\begin{array}{ccc}
 0 & \nu  &n+1   \\
 0 & n  &n+1   
\end{array}
\right)
}{B(0\,N)} \frac{\lambda_\nu}{\lambda_\infty + n -1}\mathfrak{M}_{f_1\ldots\widehat{\;f_\nu\;}\ldots f_n}(W_0(N)\varpi)\nonumber\\
&&{}\qquad \qquad \quad \quad \quad+ \frac{
B\!\left(
\begin{array}{ccc}
 0 & \star  &n+1   \\
 0 & n  &n+1   
\end{array}
\right)}{B(0\star\,n+1)}W_0(n+1)\varpi.
\eeq
(Notice that the last term in the {\rm RHS} of (19) reduces to $\varpi$.)
We may apply Proposition 23 to the RHS of (21), then we sees the formula (20) is
indeed valid in case where $j= n$ and $k = n+1$. 

The identity (20) follows by symmetry
replacing $n,n+1$ with arbitrary pair $j, k$ respectively.
\qed

\bigskip

%sec5 pf of Th 1-2
\section{Proof of Theorems 1 and  2}

 %\hspace{.5cm}
 
Assume that all the quadratics $f_j$ are normalized as in (7).
Let $*dQ$  be the $n-1$-form
\be
*dQ = \sum_{\nu = 1}^n \,(-1)^{\nu-1}\,x_\nu\,dx_1\wedge\cdots\widehat{dx}_\nu\cdots\wedge dx_n.
\ee

Since
\be
df_j\wedge *dQ = (2f_j - 2\,\sum_{\nu=1}^n \alpha_{j\nu}\,x_\nu - 2\alpha_{j0})\,\varpi ,
\ee
the  Stokes identity
\be
0 \sim \nabla(\sum_{\nu=1}^{n} (-1) ^{\nu-1} x_\nu dx_1\wedge\cdots\widehat{dx_\nu}\cdots
\wedge dx_n) 
\ee
implies

\beq
(2\lambda_\infty + n)\, \varpi \sim
\sum_{j=1}^n \lambda_j \mathfrak{M}_{f_j - f_{n+1}} F_j + \sum_{j=1}^{n+1} \lambda_j (\alpha_{j\,0} + \alpha_{n+1\,0})
F_j.
\eeq

In view of
\be
\alpha_{j0} + \alpha_{n+1\,0} = - 
B\!\left(
\begin{array}{ccc}
 0 &n+1   &j   \\
  0& \star  &j   
\end{array}
\right),
\ee
and
\be
\mathfrak{M}_{f_j - f_{n+1}} F_j = - \mathfrak{M}_{f_{n+1}- f_j}F_j,
\ee
Proposition 24 implies
\be
&&\mathfrak{M}_{f_j - f_{n+1}} F_j \sim
- \frac{
B\!\left(
\begin{array}{ccc}
 0 &n+1   &j   \\
  0& \star  &j   
\end{array}
\right)}{B(0\star j)}\,W_0(j)\varpi\\
&&{}\qquad \qquad \qquad  - \sum_{\nu=1}^n  (-1)^\nu
\sum_{\{\mu_1,\ldots,\mu_\nu\} \subset \partial_j N}\frac{\prod_{s=1}^\nu \lambda_{\mu_s}}{\prod_{s=1}^\nu (\lambda_\infty + n- s)}
\,\frac{
B\!\left(
\begin{array}{ccc}
 0 &n+1   &j   \\
  0& \mu_1  &j   
\end{array}
\right)}{B(0\mu_1\,j)}\\
&&{}\qquad \qquad \qquad \cdot\prod_{s=2}^\nu \frac{
B\!
\left(
\begin{array}{ccccccc}
 0 &\star   &\mu_{s-1}&\ldots& \mu_1& j  \\
  0& \mu_s  &\mu_{s-1} &\ldots&\mu_1&j  
\end{array}
\right)
}{B(0\mu_s\ldots\mu_1\,j)}
W_0(\mu_1\ldots\mu_\nu\,j)\varpi.
\ee
Substituting these identities into the RHS of (22),
we get the following
\beq
(2\lambda_\infty + n)\,\varpi \sim
\sum_{\nu=1}^{n+1} (-1)^\nu
\sum_{J\subset N, |J| = \nu}\frac{\prod_{j\in J}\lambda_j}{\prod_{s=1}^{\nu-1} (\lambda_\infty + n - s)}
\,\eta_J \,W_0(J)\,\varpi,
\eeq
where $J$ ranges over the family of {\it unordered} subsets
 such that $J\subset N=N_{n+1}$ and $|J| = p.$
$\eta_j, \,\eta_J$ denote
\beq
&&\eta_j = \frac{
B\!
\left(
\begin{array}{ccc}
  0&n+1   &j   \\
 0 &\star   &j    
\end{array}
\right)
+ 
B\!
\left(
\begin{array}{ccc}
  0&\star   &j   \\
 0 &\star   &n+1   
\end{array}
\right)
}{B(0\star \,j)}\quad \quad (J = \{j\}),\\
&&\eta_{J} = 
\sum_{j\in J}\sum_{\{\mu_1\,\ldots\, \mu_{\nu-1}\} \equiv \partial_j J}
\frac{
B\!
\left(
\begin{array}{ccc}
0  &n+1   &j   \\
0  & \mu_1  & j  
\end{array}
\right)
}{B(0\mu_1\,j)}
\,\frac{\prod_{s=1}^{\nu-1} 
B\!
\left(
\begin{array}{ccccccc}
 0 & \star  &  \mu_{s-1}&\ldots&\mu_1&j \\
  0&  \mu_s &\mu_{s-1} &\ldots&\mu_1&j  
\end{array}
\right)
}{\prod_{s=1}^{\nu-1}B(0\mu_s\mu_{s-1}\ldots\mu_1\,j)}\nonumber\\
&&\qquad \qquad \qquad \qquad  \qquad \qquad \qquad \qquad \qquad  (2\le \nu\le n+1,\, |J| = \nu),
\eeq
and$\{\mu_1,\ldots,\mu_{\nu-1}\}$ move over the family of all ordered sets
consisting of different elements of $\partial_jJ$.

Now as to $\eta_j,\, \eta_J$ we have the following lemma.

%lm25
\begin{lm}
\be
\eta_j = 1,\quad 
\eta_J = 1.
\ee

\end{lm}

\begin{pf}
It follows that 
\beq
&&B\!
\left(
\begin{array}{ccc}
  0&n+1   &j \\
 0 &\star   &j    
\end{array}
\right)
+ 
B\!
\left(
\begin{array}{ccc}
 0 &\star   & j  \\
 0 &  \star & n+1  \\
\end{array}
\right)
= B(0\,\star\,j),\\
&&
B\!
\left(
\begin{array}{ccc}
  0&n+1&j \\
 0 &k &j    
\end{array}
\right)
+ 
B\!
\left(
\begin{array}{ccc}
  0&n+1   &k   \\
  0& j  & k  \\
\end{array}
\right) = B(0\,k\,j).
\eeq

More generally,
for $1\le k\le n$
due to the cofactor expansion, it follows that  
%(8)
\beq
&&B(0\,\mu_1\ldots \mu_k) = \sum_{s=1}^k
B\!
\left(
\begin{array}{ccccccc}
0  &\star   & \mu_1\ldots\widehat{\mu_s} \ldots\mu_k \\
 0 &  \mu_s & \mu_1\ldots\widehat{\mu_s}\ldots\mu_k 
\end{array}
\right),\\
&&B(0\, \star\mu_1\ldots \mu_k) = 
\sum_{s=1}^k B\!
\left(
\begin{array}{ccccccc}
0  &n+1   & \star&\mu_1\ldots\widehat{\mu_s} \ldots\mu_k \\
 0 &  \mu_s &\star& \mu_1\ldots\widehat{\mu_s}\ldots\mu_k 
\end{array}
\right) \nonumber\\
&&{}\qquad \qquad \qquad \qquad +
B\!
\left(
\begin{array}{ccccccc}
0  &n+1  & \mu_1\ldots \ldots\mu_k \\
 0 &  \star & \mu_1\ldots\ldots\mu_k  
\end{array}
\right).\nonumber\\
\eeq

Suppose that the ordered set
$\{\mu_0,\mu_1,\ldots,\mu_{k-1}\}\,(k\le \nu)$ coincides with 
$L = \{l_1,\ldots,l_{k}\}\,(1< l_1<\cdots < l_{k})$
a subset of $J$. Then  from (28),
\beq
\frac{1}{(k-1)!}\frac{\sum_{\{\mu_0,\mu_1,\ldots,\mu_{k-1}\} \equiv L}
B\!
\left(
\begin{array}{ccccccc}
0  &\star   & \mu_{k-2}&\ldots&\ldots&\mu_0 \\
 0 &  \mu_{k-1} & \mu_{k-2}&\ldots&\ldots&\mu_0  
\end{array}
\right)}{B(0\,L)} = 1.
\eeq

By applying successively the identities (28) and (29), 
\be
&&\eta_J = \frac{1}{2} \sum_{\{\mu_0,\mu_1,\ldots,{\mu_\nu -1}\} \equiv J}
\prod_{s=2}^{\nu-1}
\frac{
B\!
\left(
\begin{array}{ccccccc}
 0 & \star  & \mu_{s-1} &\ldots&\mu_1&\mu_0 \\
 0 & \mu_s  &\mu_{s-1} &\ldots&\mu_1&\mu_0  
\end{array}
\right)
}{B(0\mu_s\mu_{s-1}\ldots\mu_1\mu_0)}\\
&&=
\frac{1}{k!} \sum_{\{\mu_0,\mu_1,\ldots,\mu_{\nu-1}\} \equiv J}
\prod_{s=k}^{\nu-1}\frac{
B\!
\left(
\begin{array}{ccccccc}
 0 & \star  & \mu_{s-1} &\ldots&\mu_1&\mu_0 \\
 0 & \mu_s  &\mu_{s-1} &\ldots&\mu_1&\mu_0  
\end{array}
\right)
}{B(0\mu_s\mu_{s-1}\ldots\mu_1\mu_0)}\quad(2\le k\le \nu-1)\\
&&=
\frac{1}{(\nu-1)!} \sum_{\{\mu_0,\mu_1,\ldots,\mu_{\nu-1}\} \equiv J}
\frac{
B\!
\left(
\begin{array}{ccccccc}
 0 & \star  & \mu_{\nu-2} &\ldots&\mu_1&\mu_0 \\
 0 & \mu_{\nu-1}  &\mu_{\nu-2} &\ldots&\mu_1&\mu_0  
\end{array}
\right)
}
{B(0\mu_{\nu-1}\mu_{\nu-2}\ldots\mu_1\mu_0)}
= 1.
\ee
\qed
%$\Box$
\end{pf}

\bigskip

%pf of them
%\vspace{.5cm}
[{\bf proof of Theorem 1}]

(23) and Lemma 25 imply Theorem 1. 
\qed

\bigskip
%\vspace{5cm}

We prove Theorem 2.

By definition, the covariant derivation with respect to the parameters
$\alpha_{j\nu}, \alpha_{j0}$ are given by
\beq
&&\nabla_B(\varpi) = 
d_B{\rm log}\Phi \varpi=\sum_{j=1}^{n+1}\lambda_j\frac{d_Bf_j}{f_j}\varpi\nonumber\\
%\sum_{j=1}^{n+1}\sum_{\nu=0}^{n+1-j} \,d\alpha_{j\nu}\,\nabla_{B, \frac{\partial}{\partial\alpha_{j\nu}}}\varpi\nonumber\\
&& {}\qquad \quad= \sum_{j=1}^{n+1}\lambda_j\bigl( \sum_{\nu=1}^{n-j + 1}\,2d\alpha_{j\nu}\,x_\nu + d\alpha_{j0} \bigr)\,F_j.
\eeq

%df 26
\begin{df}
Differential $1$-forms $\theta_k^j\,(1\le j\le k\le n)$
can be defined uniquely such that
\beq
d\alpha_{j\nu} = \sum_{s=j}^ {n-\nu+1}\alpha_{s\,\nu} \theta_{s}^j \quad(1\le j\le n,\,1\le \nu\le n,\, j+\nu\le n+1).
\eeq
In particular,
\beq
\theta_j^j = d\log \alpha_{j\,n+1-j}\quad(1\le j\le n).
\eeq
\end{df}
Namely, $\theta_k^j$ are related with the right invariant differential matrix $1$-form
$\omega = (\omega_{jk})_{1\le j,k\le n+2}:$
\beq
\omega = dg_\alpha\,g_\alpha^{-1},
\eeq
associated with the  $(n+2)\times (n+2)$ lower triangular matrix
$g_\alpha = (g_{jk})_{1\le j,k\le n+2}$
depending on $\alpha_{j\nu},$
where
\be
&&g_{jk} = 0\quad(j<k),\\
&&g_{11} = g_{22} = 1, \, g_{21}= \alpha_{n+1\,0},\\
&&g_{j1} = \alpha_{n+3-j\,0} ,\ g_{j2} = 1,\ g_{jk} = \alpha_{n+3-j\,k-2}\quad(3\le j \le n+2,\,3\le k\le j).
\ee

Indeed, the following Lemma holds.
%lem27
\begin{lm}
\be
&&\omega_{jk} = 0\quad(j<k),\ \omega_{11} = \omega_{22} = 0,\ \omega_{21} = d\alpha_{n+10},\\
&&\omega_{j2} = - \sum_{\nu=3}^j \omega_{j\nu}\quad(3\le j\le \,n+2),\\
&&\omega_{jk} = \theta_{n-k+3}^{n-j+3}\quad(3\le k\le j\le n+2).
\ee
\end{lm}

(32) gives the recurrence relations among $\theta_k^j:$
\beq
\theta_{k}^j = \frac{1}{\alpha_{k\,\,n+1-k}}
\{d\alpha_{j\,n+1-k} - \sum_{s=j}^{k-1} \alpha_{s\,\,n+1-k}\,\theta_{s}^j \}\quad(1\le k\le j\le n),
\eeq
which enable us to evaluate $\theta_k^j$ recurrently
using the parameters $r_j^2, \, \rho_{kl}^2$
starting from 
\be
&&\theta_j^j = d\log\alpha_{j\,n+1-j} =
\frac{1}{2}\,d\log \bigl(-\frac{B(0j\ldots \,n+1)}{B(0\,j+1\,\ldots\,n+1)}\bigr)\quad(1\le j\le n).
\ee

Furthermore,
\be
&&\alpha_{j0} =  \frac{1}{2}\bigl\{ \,B(0\,j\,n+1) - \,B(0\,\star\, j)\bigr\}\quad(1\le j\le n+1).
\ee

The RHS  of (31) can be rewritten only in terms of $f_j, \,\theta_k^j$:
\beq
\nabla_B(\varpi) = \sum_{j=1}^{n+1} \lambda_j\,V_j,
\eeq
where
\be
&&V_j = \sum_{k=j}^n \,\theta_{k}^j \, \frac{f_{k}- f_{n+1}}{f_j} \,\varpi + 
\bigl(d\alpha_{j0} - \sum_{k=j}^n (\alpha_{k0} - \alpha_{n+1\,0})\,\theta_k^j \bigr)\,F_j\\
&&{}\quad  = \sum_{k=j}^n \,\theta_k^j \, \frac{(f_k - f_j) - (f_{n+1} - f_j)}{f_j} \,\varpi 
+ \bigl(d\alpha_{j0} - \sum_{k=j}^{n} (\alpha_{k0} - \alpha_{n+1\,0})\,\theta_k^j \bigr)\,F_j,\\
&&V_{n+1} = d\alpha_{n+1\,0}\,F_{n+1}.
\ee
Hence,
\beq
&&\sum_{j=1}^{n+1} \,\lambda_j\,V_j =
\sum_{j=1}^n \lambda_j \bigl(\sum_{k= j+1}^n  \, \theta_k^j  \,\frac{f_k - f_j}{f_j}\,\varpi
- \sum_{k=j}^n \,\theta_k^j \,\frac{f_{n+1}- f_j}{f_j}\,\varpi\bigr)\nonumber\\
&&{}\qquad \qquad \quad+ \sum_{j=1}^{n+1} \lambda_j 
\bigl(d\alpha_{j0} - \sum_{k=j}^{n+1} (\alpha_{k0} - \alpha_{n+1\,0})\,\theta_k^j \bigr)\,F_j.
\eeq

Applying  Proposition 23 to  (36), we can represent 
$\nabla_B(\varpi)$ as  a linear combination of $W_0(J)\varpi$:
\beq
\nabla_B(\varpi) \sim  \sum_{\nu=1}^ {n+1}\sum_{J\subset N,\, |J| = \nu} 
\frac{\prod_{j\in J} \lambda_j}
{\prod_{s=0}^{\nu-1}(\lambda_\infty + n-s)}\, \,\tilde{\theta}_J\, W_0(J)\,\varpi.
\eeq

%re
\begin{re}
$\tilde{\theta}_J$ does not depend on $\lambda$. It is a differential $1$-form 
depending only on $r_j^2,\,dr_j^2, \,\rho_{jk}^2,\,d\rho_{jk}^2\,(j,k\in J) $.
The expression does not depend on $n$, namely  for a fixed $J\subset N,$
$\tilde{\theta}_J$ depends only on $J$.
The set of all $\tilde{\theta}_J$
has  the symmetry under
the symmetric group
$\mathfrak{S}_{n+1}$ which acts on $N$.
The set of all $\theta_J$ has the same property.
Therefore to prove {\rm Theorem 2}
we have only to prove that the explicit formula for $\tilde{\theta}_{12\ldots\,n+1}$
coincides with $\theta_{12\ldots \,n+1}$.
\end{re}

It is convenient to write a differential $1$-form $\psi$
which is a linear combination of $dr_j^2, d\rho_{jk}^2$
in the following form:
\be
\psi= \sum_{j=1}^{n+1}[\psi : dr_j^2] \,dr_j^2 +
\sum_{\mu<\nu} [\psi : d\rho_{\mu\nu}^2 ] \,d\rho_{\mu\nu}^2.
\ee

$\theta_k^j, \tilde{\theta}_J\,(|J|\ge 1)$ are written by linear combinations of $dr_j^2,\,d\rho_{jk}^2$.

The following Lemma is an immediate consequence from the definition:
%lem28
\begin{lm}
\be
&&[\theta_k^j : dr_l^2] = 0\quad(1\le l\le n+1)\,(1\le j\le k\le n),\\
&&[\theta_k^j : d\rho_{12}^2] = 0\quad(2\le j\le n).
\ee
\end{lm}

Further, a direct calculation gives
\be
&&\theta_{j+1}^j = \frac{1}{B(0\,j+1\ldots\,n+1)}
\{d\,B\!\left(
\begin{array}{cccccc}
 0 &j   &j+2\ldots& n+1   \\
  0& j+1  &j+2\ldots&n+1   
\end{array}
\right)\\
&&  - \frac{1}{2}\, B\!
\left(
\begin{array}{ccccc}
  0& j  & j+2\ldots n+1  \\
  0& j+1  & j+2\ldots n+1    
\end{array}
\right)\,
d\log\,\bigl(B(0j\,j+1\ldots n+1)\,B(0\,j+2\ldots n+1)\bigr)\}.
\ee

When $K\subset N$ and $1,2 \in K^c$, we see
\beq
&&[d\,B\!
\left(
\begin{array}{cccccc}
 0&1&K   \\
  0&2&K
\end{array}
\right): d\rho_{12}^2] = B(0\,K),\\
&&[ d\,B(012\,K) : d\rho_{12}^2] = - 2\,B\!
 \left(
\begin{array}{cccccc}
 0&1&K   \\
  0&2&K 
\end{array}
\right).
\eeq

From these relations and (30), we can deduce by recurrence
the following result.

%lm29
\begin{lm}
\be
&&[\theta_1^1 : d\rho_{12}^2] = - \frac{
B\!
 \left(
\begin{array}{cccccc}
 0&1&3&\ldots&\,n+1   \\
  0&2&3&\ldots&\,n+1  
\end{array}
\right)
}{B(0N)},\\
&&[\theta_2^1 : d\rho_{12}^2] = 
\frac{B(0134\ldots\,n+1)}{B(0 N)},\\
&&
[\theta_k^1: d\rho_{12}^2] = -
\frac{
B\!
\left(
\begin{array}{ccccccc}
 0 & 2  & 1&3&4\ldots\widehat{k}\ldots&\,n+1  \\
 0 & k  & 1 &3&4\ldots\widehat{k}\ldots&\,n+1  
\end{array}
\right)}
{B(0N)}\quad(3\le k\le n),\\
&&[\sum_{k=1}^n \theta_k^1 : d\rho_{12}^2]
= \frac{
B\!
\left(
\begin{array}{ccccccc}
 0 & 2  & 1&3&4&\ldots&\,n  \\
 0 & n+1  & 1 &3&4&\ldots&\,n  
\end{array}
\right)}{B(0N)}.
\ee
\end{lm}

In the case where $J =\{j\}\,(1\le j\le n+1),$
seeing that  $W_0(j)\,\varpi = B(0\star j) F_j$,
the following identity is derived:
\be
&&\tilde{\theta}_j = - \,
\frac{
B\!
\left(
\begin{array}{ccc}
 0 &n+1   &j   \\
  0&\star   &j   
\end{array}
\right)
}{B(0\star\,j)}\sum_{k=j}^n \theta_k^j 
+ \sum_{k=j+1}^n\frac{
B\!
\left(
\begin{array}{ccc}
 0 &k  &j   \\
  0&\star   &j   
\end{array}
\right)
}{B(0\star\,j)}\, \theta_{k}^j\\
&&{}\qquad + \frac{1}{B(0\star j)}\bigl(d\alpha_{j0} - \sum_{k=j}^{n+1} (\alpha_{k0} - \alpha_{n+1\,0}) \theta_k^j\bigr) .
\ee
In particular, in the case where $j= n+1,$
\be
\tilde{\theta}_{n+1} = - \frac{1}{2} \,d\log B(0\star\,\,n+1 ).
\ee
By symmetry also for any $j\in N,$ 
\be
\tilde{\theta}_{j} = - \frac{1}{2} \,d\log B(0\star\,j ).
\ee
In the case where $J = N,$ we have the following.

%prop30
\begin{prop}
Fix $n \ge 1$, then for $N = \{1\,2\ldots\,n+1\}$
the following equality holds:
\beq
\tilde{\theta}_{N} = \theta_{N}.
\eeq
\end{prop}

Since $\tilde{\theta}_N,\,\theta_N$ are symmetric under any permutation among
the elements of $N$, we have only to prove
the following equality
\beq
[\tilde{\theta}_N : d\rho_{12}^2] = [\theta_N : d\rho_{12}^2].
\eeq

Indeed  if (42) holds true, then
\be
[\tilde{\theta}_N : d\rho_{jk}^2 ] = [\theta_N : d\rho_{jk}^2]
\ee
holds true by symmetry for  arbitrary pair $j,k\in N$ by  the same argument.

Owing to Proposition 24,  Lemmas 28 and 29 and the identities (36), (37),
\be
&&[\tilde{\theta}_{N}: d\rho_{12}^2] 
 = 
\sum_{k=2}^n \,[\theta_{k}^1 : d\rho_{12}^2] \\
&&(-1)^n
\sum_{\{\mu_1,\ldots,\mu_{n}\} = \partial_1N}
\frac{
B\!
\left(
\begin{array}{cccccc}
  0&k   & 1 \\
  0& \mu_1  & 1 
\end{array}
\right)
}{B(0\mu_1\,1)}
%&&
\cdot
\prod_{s=2}^n \frac{
B\!
\left(
\begin{array}{cccccc}
 0 & \star  &\mu_{s-1}  &\ldots&\mu_1& 1 \\
 0 & \mu_s  & \mu_{s-1}&\ldots&\mu_1&1
\end{array}
\right)
}
{B(0\mu_s\ldots\mu_1\,1)} \\
&&{}\qquad - \sum_{k=1}^n [\theta_k^1 : d\rho_{12}^2] (-1)^n\cdot \\
&&{}\qquad
\sum_{\{\mu_1\ldots\mu_n\}= \partial_1 N }
\frac{
B\!
\left(
\begin{array}{cccccc}
  0&n+1   & 1  \\
  0& \mu_1  & 1  
\end{array}
\right)
}{B(0\mu_1\,1)}
\,
\prod_{s=1}^n \frac{
B\!
\left(
\begin{array}{cccccc}
 0 & \star  &\mu_{s-1}  &\ldots&\mu_1& 1\\
 0 & \mu_s  & \mu_{s-1}&\ldots&\mu_1&1  
\end{array}
\right)
}{B(0\mu_s\ldots\mu_1\,1)}\\
\ee
where $\{\mu_1\ldots\mu_n\}$ move over the family of all ordered sequences
consisting of different elements of $\partial_1 N$.

Namely,
\be
[\tilde{\theta}_N : d\rho_{12}^2] = U_1 + U_2 + U_3,
\ee
\be 
&&U_1 = (-1)^n\,\frac{B(0\,\partial_2 N)}{B(0\,N)}
\,
\sum_{\{\mu_1,\ldots,\mu_{n}\} = \partial_1 N}
\frac{
B\!
\left(
\begin{array}{cccccc}
  0&2   & 1 \\
  0& \mu_1  & 1 
\end{array}
\right)
}{B(0\mu_1\,1)}\\
&&{}\qquad \cdot
\prod_{s=2}^n \frac{
B\!
\left(
\begin{array}{cccccc}
 0 & \star  &\mu_{s-1}  &\ldots&\mu_1& 1 \\
 0 & \mu_s  & \mu_{s-1}&\ldots&\mu_1&1
\end{array}
\right)}{B(0\mu_s\ldots \mu_1\, 1)},\\
&&U_2 = - (-1)^n
\,\frac{B\!
\left(
\begin{array}{cccccc}
  0&2   &\partial_2\partial_{n+1}N \\
 0 &n+1   &\partial_1\partial_{n+1}N   \\ 
\end{array}
\right)
}{B(0\,N)}\, \sum_{\{\mu_1,\ldots,\mu_{n}\} = \partial_1 N}
\frac{
B\!
\left(
\begin{array}{cccccc}
  0&n+1   & 1 \\
  0& \mu_1  & 1 
\end{array}
\right)
}{B(0\mu_1\,1)}\\
&&{}\qquad \cdot
\prod_{s=2}^n \frac{
B\!
\left(
\begin{array}{cccccc}
 0 & \star  &\mu_{s-1}  &\ldots&\mu_1& 1 \\
 0 & \mu_s  & \mu_{s-1}&\ldots&\mu_1&1
\end{array}
\right)}{B(0\mu_s\ldots \mu_1\, 1)},
\ee
\be&& U_3 = - (-1)^n 
\sum_{k=3}^n
\frac{B\!
\left(
\begin{array}{ccc}
 0 &2&\partial_2\partial_k N     \\
  0&k   & \partial_2\partial_k N    
\end{array}
\right)
}{B(0\,N)}
\sum_{\{\mu_1,\ldots,\mu_{n}\} = \partial_1 N}
\frac{
B\!
\left(
\begin{array}{cccccc}
  0&k   & 1 \\
  0& \mu_1  & 1 
\end{array}
\right)
}{B(0\mu_1\,1)}\\
&&{}\qquad \cdot
\prod_{s=2}^n \frac{
B\!
\left(
\begin{array}{cccccc}
 0 & \star  &\mu_{s-1}  &\ldots&\mu_1& 1 \\
 0 & \mu_s  & \mu_{s-1}&\ldots&\mu_1&1
\end{array}
\right)}{B(0\mu_s\ldots \mu_1\, 1)}.
\ee

By cofactor expansion, the following holds.
%lem31
\begin{lm}
\beq
&& B\!
\left(
\begin{array}{ccc}
 0 &2   &1   \\
 0 &l  &1    
\end{array}
\right)
\,B(0\,\partial_2\,N)
-
\sum_{k=3}^n 
B\!
\left(
\begin{array}{ccc}
 0 &k   &1   \\
 0 &l   &1    
\end{array}
\right)
B\!
\left(
\begin{array}{ccc}
 0 &2   &\partial_2\partial_k N   \\
  0&k   & \partial_2\partial_k N  
\end{array}
\right) \nonumber\\
&&
- B\!
\left(
\begin{array}{ccc}
 0 &n+1   &1   \\
 0 &l   &1    
\end{array}
\right)
\,B\!
\left(
\begin{array}{cccccc}
  0&2   &\partial_2\partial_{n+1} N   \\
 0 &n+1   &\partial_2\partial_{n+1} N   
\end{array}
\right)
=
\left\{
\begin{array}{ccc}
- B(0\,N) & (l = 2)   \\
0  & (l\in \partial_2 N)  &  
\end{array}.
\right.\nonumber\\
\eeq
\end{lm}

Indeed, 
let $Y$ be the $n\times n$ matrix
whose $j-1,k-1$th component 
$(j,k\ge 2)$ is
\be
y_{jk} = B\!
\left(
\begin{array}{ccc}
 0 & 1  &j   \\
 0 &  1 & k  
\end{array}
\right)\quad(j,k \in \partial_1 N),
\ee
and let the $p-1$ order minor consisting of 
$j_1-1,\ldots,j_p-1$th rows, and  $k_1-1,\ldots, k_p-1$th  columns
be denoted by
\be
Y\!
\left(
\begin{array}{ccc}
 j_1 &\ldots   &j_p   \\
 k_1 & \ldots  &k_p   
\end{array}
\right)
= 
\left|
\begin{array}{cccccc}
 y_{j_1k_1} &\ldots   &y_{j_1k_p}   \\
  \vdots&  \ddots & \vdots  \\
 y_{j_pk_1} & \ldots  &y_{j_pk_p}   
\end{array}
\right|
\ee
(in particular when $j_\nu = k_\nu,$ it is denoted by $Y(j_1\ldots j_p)$).
Then the LHS of (43) equals
\beq
(-1)^n\, Y\!\left(
\begin{array}{cccccc}
 l &3&\ldots   &n+1   \\
  2&3 &\ldots  &n+1   
\end{array}
\right) =
- B\!
\left(
\begin{array}{cccccc}
 0 &1   &l&3&\ldots&n+1   \\
  0& 1  &2&3&\ldots&n+1   
\end{array}
\right). 
\eeq
Hence this equals $- B(0\,N)$ for $l = 2$
and $0$ for $3\le l\le n+1$.
\qed

\bigskip

%pf of prop30
%\vspace{.5cm}
[{\bf proof of Proposition 30}]

According to Lemma 31,
\be
&&U_1 + U_2 + U_3\\
&& = (-1)^{n+1} \frac{1}{B(0\,2\,1))}
\sum_{\{\mu_2,\ldots,\mu_n\} = \partial_1\partial_2 N}
\prod_{s=2}^n
\,\frac{B\!
\left(
\begin{array}{ccccccc}
 0 &\star   &\mu_{s-1}&\ldots&\mu_2&2&1   \\
 0 &\mu_s&\mu_{s-1}&\ldots&\mu_2&2&1     
\end{array}
\right)
}{B(0\,\mu_s\,\ldots\,\mu_2\,2\,1)}.
\ee
This means (42).
Since $\{\theta_J\}_J$ and $\{ \tilde{\theta}_J\}_J$ are symmetric 
under the action of ${\mathfrak S}_{n+1},$
Proposition 30 has been proved.
%$\Box$
\qed

\medskip

Likewise the following Proposition holds.

%prop32
\begin{prop}
For $J\,\subset N\, (|J| \ge 1),$ we have
\beq
\tilde{\theta}_J = \theta_J.
\eeq
\end{prop}

\begin{pf}
Because of symmetry, we have only to prove it when
$J = \{ k+1,\ldots, n+1\}\,(1\le k\le n)$.
We can apply the preceding Proposition 30
in the special case $\lambda_l = 0 \,(1\le l\le k),$
where the multiplicative meromorphic function
\be
\Phi_k(x) = \prod_{j= k+1}^{n+1} \,f_j(x)^{\lambda_j}
\ee
is taken.

On the other hand,
under the condition  $x_1= \ldots = x_k = 0,$
we put
\be
y = (x_{k+1},\ldots, x_{n}),\,
\varpi_k = dx_{k+1}\wedge\cdots\wedge dx_{n}
\ee
and consider the multiplicative function of $y$
\be
&&\Phi_k^*(y) = \prod_{j=k+1}^{n+1} \,f_j^*(y)^{\lambda_j},\\
&&f_j^*(y) = \sum_{j = k+1}^{n+1}\,x_{j}^2 + \sum_{\nu=k+1}^{n+1}
\,2\alpha_{j\nu} x_\nu + \alpha_{j0}\quad(k+1\le j\le n+1).
\ee

In the complement  $X_k : = {\bf C}^{n-k} - \bigcup_{j=k+1}^{n+1} \{f_j^*(y) = 0\},$
we may apply Proposition 30 to
the $n-k$ dimensional twisted de Rham cohomology $H_\nabla^{n-k}(X_k, \Omega^\cdot(* S))$
associated with $\Phi_k^*(y)$.
From the definition, the basic invariants of $\Phi_k, \Phi_k^*$ coincide with
$r_j^2, \rho_{jk}^2\,(k+1\le j,k\le n+1)$. The only one difference is that
$\lambda_\infty + n - \nu $ must be replaced by $\lambda_\infty + n - k -  \nu$ 
in the case of $\Phi_k^*$.
Let the differential $1$-form $\tilde{\theta}_J\,(J \subset N_k^c = \{k+1,\ldots, n+1\})$
corresponding to $\Phi_k^*(y)$ be denoted by $\tilde{\theta}_J^*$.
Then the following holds:
\be
&&\nabla_B(\varpi) \sim \sum_{\nu = 1}^{n+1-k} \sum_{J\subset N_k^c, |J| = \nu}
\frac{\prod_{j\in J} \lambda_j}{\prod_{s=1}^{\nu-1} (\lambda_\infty + n - s)} \,\tilde{\theta}_J\,W_0(J)\,\varpi,\\
&&\nabla_B(\varpi_k) \sim \sum_{\nu = 1}^{n+1-k} \sum_{J\subset N_k^c, |J| = \nu}
\frac{\prod_{j\in J} \lambda_j}{\prod_{s=1}^{\nu-1} (\lambda_\infty + n - k- s)} \,\tilde{\theta}_J^*\,W_0(J)\,\varpi_k,
\ee
Furthermore,
\be
\tilde{\theta}_{J}^* = \tilde{\theta}_{J}\quad(J \subset N_k^c).
\ee
In particular,
\beq
&&\tilde{\theta}_{k+1\,\ldots\,n+1}^* = \tilde{\theta}_{k+1\,\ldots\,n+1}.
\eeq
As was proved in {\rm Proposition 30},
\beq
\tilde{\theta}_{k+1\,\ldots\, n+1}^*
=\theta_{k+1\,\ldots\,n+1}.
\eeq
{\rm(46) and (47)} means
\be
\tilde{\theta}_{k+1\,\ldots\,n+1} = \theta_{k+1\,\ldots\,n+1}.
\ee
In this way, {\rm(45)} has been proved in case  $J = \{k+1,\ldots,n+1\}$.
Because of symmetry on $N,$
{\rm (45)} holds true for an arbitrary $J \subset N \,(|J| = n- k + 1)$. $\Box$
\end{pf}

\bigskip

%pf of th 2
%\vspace{.5cm}
[{\bf proof of Theorem 2}]

$\nabla_B(\varpi)$ is expressed in the form of (31).
Owing to  Proposition 32,
$\theta_J$ and $\tilde{\theta}_J$ coincide with each other
for all $J \subset N \,(|J|\ge 1)$, 
and hence the RHS of (38) coincides with the RHS of (5).
\qed

\bigskip

%sec6,p35
\section{Case  for general $m$}

%\hspace{.5cm}

We now want to show that Theorems 1 and 2 still hold for general $m$ such that
$n+1 \le m.$
Since the proof can be done almost in the same procedure, we roughly sketch the proof,
focusing on different points of argument.

Denote by $N= N_{n+1}$ the set of indices $\{1,2,\ldots, n+1\}.$
We may assume $f_j \,(1\le j\le n+1)$ have the same expression as in (7).

For an arbitrary (unordered) subset $J \subset N_m,$
denote by $W_0(J)\varpi$ the form
\be
&&W_0(J)\,\varpi = - \sum_{j\in J} B\!
\left(
\begin{array}{cccc}
 0 &j&\partial_j J    \\
  0&\star&\partial_j  J \\
\end{array}
\right)\,F_{\partial_j J}
+ B(0\star\,J) F_J.
\ee

Since both $B\!
\left(
\begin{array}{cccc}
 0 &j&\partial_j J    \\
  0&\star&\partial_j  J \\
\end{array}
\right)
$ and $B(0\star J)$ vanish provided $|J| \ge n+3,$
$W_0(J)\varpi$ vanishes too.
 
 In the case $|J| = n+2,$ both $B(0\,J)$ and $W_0(J)\varpi$ vanish: 
 \be
 B(0\,J) = 0,  \quad W_0(J)\varpi = 0, 
\ee
so that the following partial fraction expansion holds:
\beq
 B(0\star J) \,F_J = \sum_{\nu =1}^{n+1} B\!
\left(
\begin{array}{ccc}
 0 &\star   &\partial_{j_\nu} J  \\
 0 &j_\nu   & \partial_{j_\nu}J  \\
\end{array}
\right)\,F_{\partial_{j_\nu}J} .
 \eeq

In this way, $F_J\,(|J| \ge n+2) $ can be represented by a linear combination
of $F_K\, (|K| = n+1),$ owing to $B(0\star J) \not= 0\,(|J| = n+2)$ by assumption.

Furthermore, similarly to Lemma 12
, for $|J| = n+1, \,J = \{j_1,\ldots, j_{n+1}\}\subset N_m,$ 
$W(J)\varpi$ denotes
\be
 2^{-n}\sum_{\nu=1}^{n+1}\,(-1)^{\nu-1}
\,e_{
j_1\ldots\widehat{j_\nu}\ldots\,j_{n+1}
} = W(J)\varpi .
\ee

If $j_1 < \ldots < j_{n+1},$ then
\beq
W(J)\varpi = (-1)^{\frac{n(n-1)}{2}+1}
 \frac{1}{\sqrt{(-1)^{n+1}\,2^n\,B(0\,J)}}\,W_0(J)\varpi.
 \eeq

Since,
for an arbitrary ordered set $J = \{j_1,\ldots, j_{n+2}\}  \subset N_m, \,(|J| = n+2)$
such that $j_1<\ldots < j_{n+2},$
 \beq
\sum_{\nu=1}^{n+2} W(\partial_{j_\nu} J)\varpi = 0,
\eeq
 we get 
 \be
 W_0(\partial_{j_{n+2}}J)\varpi - \sum_{\nu=1}^{n+1}
 \frac{B\!
\left(
\begin{array}{ccccc}
 0 &j_{n+2}   &\partial_{j_{n+2}}\partial_{j_\nu}J   \\
  0&j_\nu   &\partial_{j_{n+2}}\partial_{j_\nu}J
\end{array}
\right)}{B(0\,\partial_{j_\nu}J)}W_0(\partial_{j_\nu} J)\varpi
 = 0,
 \ee
 in view of the Jacobi identity
\be
 B(0 \,\partial_{j_\mu}J)\,
 B(0 \,\partial_{j_\nu}J) 
 - B^2\!
\left(
\begin{array}{ccccc}
  0&j_\mu   &\partial_{j_\mu}\partial_{j_\nu}J   \\
  0&j_\nu   & \partial_{j_\mu}\partial_{j_\nu}J  \\
\end{array}
\right) = 0.
\ee

\medskip

%re
\begin{re}
$H_\nabla^n(X,\Omega^\cdot(*S))$ is generated by the forms $F_J\, \,(1\le |J| \le n+1)$
called \lq\lq admissible".
However they are no more linearly independent. 

Owing to {\rm Proposition 21, Proposition 24} and {\rm Theorem A.1}  in the {\rm Appendix}
(50) means that we can take as a basis of
$H_\nabla^n(X,\Omega^\cdot(*S))$ the collection of forms 
$F_J\,(1\le |J| \le n)$ and $F_J\, (|J| = n+1)$ such that  $ n+1\in J,$
or equivalently $W_0(J)\varpi\,(1\le |J| \le n)$ and $W_0(J)\varpi\, (|J| = n+1)$
such that $n+1 \in J$.  
The set of indices $J$ satisfying this property will be denoted by $\cal B.$
$\cal B$ can be regarded as \lq\lq NBC (non broken circuit) " of $N_m$
with respect to the total ordering $\cal O$ of
$N_m$:
\be
{\cal O}  : \quad  n+1 \prec 1\prec \cdots \prec n \prec n+2\prec \cdots\prec m,
\ee
where $J$ is dependent if and only if $|J| \ge n+2$ {\rm(}see {\rm [17],[18]} for NBC basis{\rm)}.
\end{re}

Then (50) shows the similar identity modified from (10) holds true:
%lem33
\begin{lm}
For $J \subset N_m, \,|J| = n,$
\beq
e_J \sim 2^n \sum_{k \in  J^c} \,\frac{\lambda_k}{\lambda_\infty}\,W(k\,J)\varpi,
\eeq
In particular,
\beq
&&(-1)^{\nu-1} \,e_{1\ldots\widehat{\nu}\ldots n+1} \sim 
2^n 
\frac{(-1)^{
\frac{n(n-1)}{2}+1
}}
{\sqrt{(-1)^{n+1}2^n\,B(0\,N)}}
\{
\frac{\lambda_\nu}{\lambda_\infty}W_0(N)\varpi \nonumber\\
&&+ \sum_{k=n+2}^m 
\frac{\lambda_k}{\lambda_\infty}\frac{B\!
\left(
\begin{array}{ccc}
 0 &\nu&\partial_\nu N      \\
  0&k   &\partial_\nu N   \\
\end{array}
\right)
}
{B(0\,k\,\partial_\nu N)}\,W_0(k\,\partial_\nu N)\varpi\}\quad(1\le \nu \le n+1).
\eeq
\end{lm}

%pf
\begin{pf}
 Indeed, the Stokes formula shows that for any $K\subset N_m, \,|K| = n-1,$
  \beq
 0 \sim \nabla(e_K) = \sum_{k\notin K} \lambda_k\, e_{k\, K}.
 \eeq
We can uniquely solve (53) with respect to  $e_J\,(|J| = n)$
in terms of  $W(L)\varpi $\\$(|L| = n+1)$ and gets (51).
Furthermore, (49) and (51) imply (52) in view of the Jacobi identity
\be
 B(0 N)\,B(0\,k\,\partial_\nu N)
  - B^2\!
\left(
\begin{array}{cccc}
 0 &\nu   &\partial_\nu N   \\
 0 &k&   \partial_\nu  N \\
\end{array}
\right) = 0.
\ee
\qed
\end{pf}
%*
%*
%*

We can see now from (52) that a similar identity to (17)  still holds true  
by the substitutions
\be
&&J^c \longrightarrow  J^c\cap N,\\
&&\lambda_\nu \,W_0(N)\varpi 
\longrightarrow
\lambda_\nu \,W_0(N)\varpi \nonumber
 + \sum_{k=n+2}^m 
\lambda_k\,\frac{B\!
\left(
\begin{array}{ccc}
 0 &\nu&\partial_\nu N      \\
  0&k   &\partial_\nu N   \\
\end{array}
\right)
}
{B(0\,k\partial_\nu N)}\,W_0(k\,\partial_\nu N)\varpi .
\ee

Namely according to Lemma 13, we have
the analog of Proposition 21.

%prop34
\begin{prop}
For admissible $J \subset N_m$ and $j\in J^c,$ the following identity holds{\rm :}
\be
&&\mathfrak{M}_{f_j} F_J \sim \{\sum_{\nu\in N\cap J^c}
\frac{B\!
\left(
\begin{array}{ccc}
  0&\nu   &J   \\
 0 &j   &J   \\
\end{array}
\right)
}{B(0\,N)\,B(0\,J)}\,\frac{1}{\lambda_\infty + n - p}\\
&&\cdot \Bigl(
\lambda_\nu 
\prod_{l\in \partial_\nu J^c\cap N}
\mathfrak{M}_{f_l}
\,W_0(N)\varpi + 
\sum_{k=n+2}^m 
\lambda_k
\frac{B\!
\left(
\begin{array}{ccc}
 0 &\nu&\partial_\nu N      \\
  0&k   &\partial_\nu N   \\
\end{array}
\right)
}
{B(0\,k\,\partial_\nu N)}\,
\prod_{l\in \partial_\nu J^c\cap N}
\mathfrak{M}_{f_l}W_0(k\,\partial_\nu N)\varpi\Bigr)\}\\
&&- \frac{B\!
\left(
\begin{array}{ccc}
 0 &\star   &J   \\
  0&j   &J   \\
\end{array}
\right)
}{B(0\,J)}\,F_J + \sum_{\nu\in J}
\frac{B\!
\left(
\begin{array}{ccc}
 0 &\nu   &\partial_\nu J   \\
  0&j   &\partial_\nu J   \\
\end{array}
\right)
}{B(0\,J)}\,F_{\partial_\nu J},
\ee
where $F_{\partial_\nu J}$ denotes $\varpi$ if $\partial_\nu J = \emptyset$ {\rm (the empty set)}.
\end{prop}

Corresponding to Lemma 22, we get the following 
recurrence formula.

%lm35,p39
\begin{lm}
For $J\subset N,\,1\le |J| = p \le n,$
\be
&&\mathfrak{M}_{f_J}W_0(N)\varpi \sim
- \sum_{\nu\in J} \frac{B\!
\left(
\begin{array}{ccc}
  0&\star   &J^c\cap N   \\
  0&\nu   &J^c\cap N   \\
\end{array}
\right)
}{B(0\,J^c\cap N)} \frac{1}{\lambda_\infty+p-1}\\
&&{}\qquad \cdot \mathfrak{M}_{f_{\partial_\nu J}}
\{\lambda_\nu\,W_0(N)\varpi
+ \sum_{k=n+2}^m 
\lambda_k\frac{B\!
\left(
\begin{array}{ccc}
 0 &\nu&\partial_\nu N      \\
  0&k   &\partial_\nu N   \\
\end{array}
\right)
}
{B(0\,k\,\partial_\nu N)}\,W_0(k\,\partial_\nu N)\varpi\}\\
&&{}\qquad+ \frac{B(0\,N)}{B(0\,J^c\cap N)}\,W_0(J^c\cap N)\varpi 
+ (-1)^n \delta_{pn}\,B(0\,N)\,\varpi.
\ee
\end{lm}

Repeated application of the preceding Lemma, we have
an extension of Proposition 23 stated as follows.

%prop36
\begin{prop}
Let $J\subset N$ be given and fixed.
Under the same condition as above,
\be
&&\mathfrak{M}_{f_J}W_0(N)\varpi \sim
(-1)^n\,\delta_{pn}\,B(0\,N)\varpi +\\
&&
\sum_{q=0}^p \,(-1)^q \sum_{K = \{k_1,\ldots,k_q\}\subset J\cup N^c}\frac{\prod_{\nu=1}^q\lambda_{k_\nu}}
{\prod_{\nu=1}^q (\lambda_\infty+p-\nu)}\, \frac{B(0\,N)}{B(0\,J^c\cap N)}\\
&&\sum_{L = \{\mu_1,\ldots,\mu_q\}\subset J}
\frac{
B\!
\left(
\begin{array}{ccc}
 0 &\star   &J^c\cap N   \\
  0&\mu_1   &J^c\cap N   \\
\end{array}
\right)
\prod_{s=2}^q 
B\!
\left(
\begin{array}{cccccc}
 0 &\star   &k_{s-1}&\ldots&k_1&J^c\cap N   \\
 0 &\mu_s   &k_{s-1} &\ldots&k_1&J^c\cap N  \\
\end{array}
\right)
}
{\prod_{s=1}^q B(0\,k_s\ldots k_1\,J^c\cap N)}\\
&&
\cdot \frac{\prod_{s=1}^q \,
B\!
\left(
\begin{array}{cccccc}
 0 &k_s   &k_{s-1}&\ldots&k_1&\partial_{\mu_s}\cdots\partial_{\mu_1}N   \\
 0 & \mu_s  &k_{s-1} &\ldots&k_1&\partial_{\mu_s}\cdots\partial_{\mu_1}N  \\
\end{array}
\right)
}
{\prod_{s=1}^q\,B(0\,k_s\ldots k_1\,\partial_{\mu_s}\ldots\partial_{\mu_1}N) }
W_0(k_1\ldots\,k_q\,J^c\cap N)\varpi.
\ee
where $M = \{\mu_1,\ldots, \mu_q\}$ ranges over the family of all ordered subsets
of $J$ consisting of $q$ different elements.
$K = \{k_1,\ldots, k_q\}$ ranges over the family of all ordered subsets of $J\cup N^c$ consisting
of different elements.
\end{prop}

Finally, take the special case $J = \{n+1\},\, j = n.$ Then Propositions 34 and 36 show
the following generalization of (20).

%lem37
\begin{lm}

\beq
&& \mathfrak{M}_{f_n - f_{n+1}}F_{n+1} \sim 
\sum_{q=0}^{n} (-1)^q \,V_q,
\eeq
where

\beq
&&V_0 = 
\frac{B\!
\left(
\begin{array}{ccc}
 0 & \star  &n+1   \\
 0 & n  &  n+1 \\
\end{array}
\right)}{B(0\star\, n+1)}\,W_0(n+1)\varpi,\\
%V_q
&&V_q = 
\sum_{\{k_1,\ldots,k_q\}\subset \partial_{n+1}N_m }
\frac{\prod_{s=1}^q \lambda_{k_s}}{\prod_{s=1}^q (\lambda_\infty + n - s)}\nonumber\\
&&\sum_{\{\mu_1,\ldots,\mu_q\}\subset \partial_{n+1}N}
\frac{B\!
\left(
\begin{array}{ccc}
 0 & n  &n+1   \\
 0 & \mu_1  &  n+1 \\
\end{array}
\right)\,\prod_{s= 2}^q 
 B\!
\left(
\begin{array}{cccccc}
  0&\star   & k_{s-1}&\ldots&k_1&n+1  \\
  0&\mu_s   &k_{s-1}&\ldots&k_1&n+1   \\
\end{array}
\right)
}
{\prod_{s=1}^q
B(0\,k_s\,\ldots\,k_1\,n+1)}\nonumber
\eeq
\beq
&&{}\quad \cdot
\frac{\prod_{s=1}^q  
B\!
\left(
\begin{array}{cccccc}
 0 &k_s   &k_{s-1}&\ldots&k_1&\partial_{\mu_s}\ldots\partial_{\mu_1}N   \\
  0& \mu_s  &k_{s-1}&\ldots&k_1&\partial_{\mu_s}\ldots\partial_{\mu_1}N   \\
\end{array}
\right)
}
{\prod_{s=1}^q B(0\,k_s,\ldots,k_1\,\partial_{\mu_s}\ldots\partial_{\mu_1}N)}\nonumber\\
&&\nonumber\\
&&{}\quad \cdot
\frac{B(0\,k_q\ldots k_1\,\partial_{\mu_q}\cdots\partial_{\mu_1}N)}{B(0\,N)}\,
W_0(k_q\ldots k_1\,n+1)\varpi\quad(1\le q\le n).
\eeq
Here $\{k_1,\ldots,k_q\}$ ranges over the family of all ordered subsets of $\partial_{n+1}N_m $
cosisting of different elements, while $\{\mu_1,\ldots, \mu_q\}$ ranges over all
ordered subsets of $N$ satisfying the following property
\be
N\cap \{k_1,\ldots, k_s\} \subset \{\mu_1,\ldots,\mu_s\}\quad(1\le s\le q).
\ee
Every other term in the {\rm RHS} of {\rm(56)} vanishes with no contribution to the summation.
\end{lm}
 
 As a result,
 %co38
\begin{co}
The form $\frac{f_n - f_{n+1}}{f_{n+1}}\,F_{n+1}$ can be represented by 
a linear combination of $W_0(J)\varpi\,(1\le |J| \le n+1)$ such that $n+1 \in J,$
i.e., $J\in {\cal B}$.
In particular, the coefficients of $W_0(J)\varpi$ such that $J \subset N$
coincide with the ones in {\rm (21)}.
\end{co}
 
 The formula (51) can be simplified into the final formula analogous to (21)
 for $j= n$ and $k= n+1$. Namely,
a generalization of Proposition 24 holds true for general $m$
replacing $N = N_{n+1}$ with $N_m$:

%prop39
\begin{prop}
Take $j,k \in N_m$ such that $j\ne k.$ Then
we get
\beq
&&\mathfrak{M}_{f_j-f_k}\,F_{k} \sim
\frac{B\!
\left(
\begin{array}{cccccc}
 0 &j& k \\
 0 &\star & k  \\
\end{array}
\right)
}
{B(0\star \,k)
}\,W_0(k)\varpi\nonumber\\
&& + \sum_{q=1}^n (-1)^q
 \sum_{K = \{k_1,\ldots,k_q \}\subset \partial_{n+1} N_m }\frac{\prod_{s=1}^q \lambda_{k_s}}{\prod_{s=1}^q(\lambda_\infty + n - s)}\nonumber\\
&&
 \cdot
 \frac{B\!
\left(
\begin{array}{cccccc}
 0 &j& k  \\
 0 &k_1 &k  \\
\end{array}
\right)
}
{B(0\,k_1 \,k)}
 \prod_{s=2}^q\frac{B\!
\left(
\begin{array}{cccccc}
 0 &\star&k_{s-1}&\ldots&k_{1} &k  \\
 0 &k_s&k_{s-1}&\ldots&k_{1} &k  \\
\end{array}
\right)
}
{B(0\,k_s\ldots k_1\,k)}
 W_0(k_1\ldots k_q\, k)\varpi.\nonumber\\
\eeq
where $K = \{k_1,\ldots,k_q\}$ ranges over the family of all ordered subsets of
$\partial_{k}N_m$ consisting of $q$ different elements of $\partial_{k}N_m$.
\end{prop}

To prove Proposition 39, we may assume that $j = n,\,k= n+1.$ 
As is seen from {\rm(54)}, we have an expression
\beq
\frac{f_n - f_{n+1}}{f_{n+1}}\,\varpi \sim
 \sum_{q=0}^n \sum_{K\subset \partial_{n+1}N_m,\,|K| = q}
\frac{\prod_{k\in K} \lambda_k}{\prod_{s=1}^q (\lambda_\infty + n - s)}
\eta_{K}
\,W_0(K\,n+1)\varpi,
\eeq
where $K = \{k_1,\ldots,k_q\}$ ranges over the family of all ordered subsets
of $\partial_{n+1}N_m$ consisting of different elements
and $\eta_K$ is independent of the exponents $\lambda$.
If $K \subset \partial_{n+1}N$, then Corollary 38 shows that
$\eta_K$ coincides with the coefficients of $W_0(K\,n+1)\varpi$
written in the RHS of (57).
Since $W_0(K)\varpi\,(1\le |K|\le n, K\subset \partial_{n+1}N_m)$
and  $W_0(K\,n+1)\varpi\,(0\le |K|\le n, K\in \partial_{n+1}N_m)$
make as representative a basis of $H_\nabla^n(X,\Omega^\cdot(*S)),$
the expression is unique. Hence, because of symmetry for expression, 
the coefficients $\eta_K$ must be described as in the RHS of (57)
provided $|K| \le n-1.$
Therefore, it is sufficient to prove that $\eta_K$ coincides with the one
in the RHS of (57) in the case where $|K| = n .$

Namely, we must prove the following identity:
\beq
&&V_n =  \sum_{K = \{k_1,\ldots,k_n \}\subset \partial_{n+1} N_m }\frac{\prod_{s=1}^n \lambda_{k_s}}{\prod_{s=1}^n(\lambda_\infty + n - s)}\nonumber\\
&&
 \cdot\frac{B\!
\left(
\begin{array}{cccccc}
 0 &n& n+1  \\
 0 &k_1 &n+1  \\
\end{array}
\right)
}
{B(0\,k_1 \,n+1)}
 \prod_{s=2}^q\frac{B\!
\left(
\begin{array}{cccccc}
 0 &\star&k_{s-1}&\ldots&k_{1} &n+1  \\
 0 &k_s&k_{s-1}&\ldots&k_{1} &n+1  \\
\end{array}
\right)
}
{B(0\,k_s\ldots k_1\,n+1)}
 W_0(k_1\ldots k_n\, n+1)\varpi, \nonumber\\
\eeq
where $K$ ranges over the family of ordered subsets of $\partial_{n+1}N_m.$

To do it, we need three preliminary lemmas.
The first one is stated as follows:
%lem40
\begin{lm}
For an ordered subset $M = \{\mu_1,\ldots,\mu_q\}\subset N\,(1\le q\le n),$
\beq
&&\frac{\prod_{s=1}^q 
B\!
\left(
\begin{array}{cccccc}
 0 &k_s   &k_{s-1}&\ldots&k_1&\partial_{\mu_s}\ldots\partial_{\mu_1}N   \\
  0& \mu_s  &k_{s-1}&\ldots&k_1&\partial_{\mu_s}\ldots\partial_{\mu_1}N   \\
\end{array}
\right)
}
{\prod_{s=1}^{q-1}\,
 B(0\,k_s,\ldots,k_1\,\partial_{\mu_s}\ldots\partial_{\mu_1}N)}\nonumber\\
&&=
B\!
\left(
\begin{array}{cccccc}
 0 &k_q   &k_{q-1}&\ldots&k_1&\partial_{\mu_q} \cdots\partial_{\mu_1}N  \\
  0& \mu_q  &\mu_{q-1}&\ldots&\mu_1&\partial_{\mu_q}\cdots\partial_{\mu_1}N   \\
\end{array}
\right)\quad(1\le q\le n).
\eeq
In particular,
\beq
&&\frac{\prod_{s=1}^n
B\!
\left(
\begin{array}{cccccc}
 0 &k_s   &k_{s-1}&\ldots&k_1&\partial_{\mu_s}\ldots\partial_{\mu_1}N   \\
  0& \mu_s  &k_{s-1}&\ldots&k_1&\partial_{\mu_s}\ldots\partial_{\mu_1}N   \\
\end{array}
\right)
}
{\prod_{s=1}^{n-1}\,
 B(0\,k_s,\ldots,k_1\,\partial_{\mu_s}\ldots\partial_{\mu_1}N)}\nonumber\\
&&=
B\!
\left(
\begin{array}{cccccc}
 0 &k_n   &k_{n-1}&\ldots&k_1&n+1  \\
  0& \mu_n  &\mu_{n-1}&\ldots&\mu_1&n+1   \\
\end{array}
\right).
\eeq
\end{lm}

%proof
\begin{pf}
{\rm(60)} can be derived  by induction
using  the following Jacobi identity
\be
B\!
\left(
\begin{array}{ccc}
 j_1   &J \\
 k_1   &K   \\
\end{array}
\right)
B\!
\left(
\begin{array}{ccc}
 j_2   &J \\
 k_2    &K   \\
\end{array}
\right)
=B\!
\left(
\begin{array}{ccc}
 j_1   &J \\
 k_2   &K   \\
\end{array}
\right)
B\!
\left(
\begin{array}{ccc}
 j_2   &J \\
 k_1   &K   \\
\end{array}
\right)
\ee
for $|J| = |K| = n,$
because of
\be
B\!
\left(
\begin{array}{ccc}
 j_1   &j_2&J \\
 k_1   &k_2&K   \\
\end{array}
\right) = 0.
\ee
\qed
\end{pf}

The following identity is an immediate consequence of {\rm(60) and (61)}:
\beq
&&B\!
\left(
\begin{array}{ccccc}
 0 &\mu_1 &\ldots  &\mu_n& n+1   \\
 0 &k_1   &\ldots&k_n&n+1   \\
\end{array}
\right)
B\!
\left(
\begin{array}{cccccc}
 0 &\mu_1 &\ldots  &\mu_{n-1}&\mu_n& n+1   \\
 0 &k_1   &\ldots&k_{n-1}&n&n+1   \\
\end{array}
\right)\nonumber\\
&&=
B(0\,N)\,
 B\!
\left(
\begin{array}{cccccc}
 0 &k_1 &\ldots  &k_{n-1}&n& n+1   \\
 0 &k_1   &\ldots&k_{n-1}&k_n&n+1   \\
\end{array}
\right)
\eeq
for $\{\mu_1,\mu_2\ldots,\mu_n\} \equiv\{1,2,\ldots,n\}$
as a set.

\medskip

The second one follows from Laplace expansion of determinant.

%lm41
\begin{lm}
Denote  $K_\nu = \{k_1,\ldots,k_\nu\} \subset N_m, M_\nu=\{\mu_1,\ldots,\mu_\nu\}\subset N_{n+1}\,(1\le \nu\le n)$
respectively {\rm (where }$K_0$ means the empty set{\rm )}.
Then
\beq
&&\sum_{\nu=1}^{p+1}
(-1)^{p+1-\nu}\, B\!
\left(
\begin{array}{cccccccccc}
 0 &\partial_{\mu_\nu} M_{p+1}&&n+1  \\
  0& K_{p-1}&n&n+1
\end{array}
\right)
B\!
\left(
\begin{array}{ccccc}
 0 &\star   &K_p   &n+1\\
 0 &\mu_\nu&K_p &n+1  \\
\end{array}
\right)\nonumber\\
&&=
B\!
\left(
\begin{array}{cccccc}
 0 &M_{p+1} &&n+1  \\
0  & K_p&n&n+1   \\
\end{array}
\right)
B\!
\left(
\begin{array}{cccccc}
0 &K_{p-1}&\star&n+1  \\
0 &K_p&&n+1   \\
\end{array}
\right)\nonumber\\
&&{}\quad + 
B\!
\left(
\begin{array}{cccccccc}
0 &M_{p+1} &&&n+1  \\
0 & K_{p-1}&n&\star&n+1   \\
\end{array}
\right)
\,B(0\,K_p\,n+1).
\eeq
\end{lm}

\begin{pf}
Define the matrix $Y= (y_{\mu\,k})$ with the $\mu,\,k$th elements
$y_{\mu\, k}=$ 
$B\!
\left(
\begin{array}{ccc}
 0 &n+1   &\mu   \\
  0&n+1   & k  \\
\end{array}
\right).$
Here $\mu$ or $k$ should be replaced with $\star$
if necessary.

Consider the determinant of the $(2p+1)\times (2p+1)$
matrix $Z= (z_{\nu\,s})$ defined by putting
\be
&&z_{\nu1} = y_{\nu\,n}\,(1\le \nu\le p+1),\,z_{\nu1} = 0\ (p+2\le \nu\le 2p+1)\;{\rm for }\, s=1;\\
&&z_{\nu s} = y_{\nu\,k_{s}}\,(1\le \nu\le p+1),z_{\nu s} = 0\ (p+2\le \nu\le 2p+1)\;{\rm for}\, 2\le s\le p;\\
&&z_{\nu\, p+1 } = y_{\nu \star}\,(1\le \nu\le p+1), z _{\nu p+1}= y_{\nu- p-1\,\star}\ (p+2\le \nu\le 2p+1);\\
&&z_{\nu\,s} = y_{\nu\,k_{s-p-1}}\,(1\le \nu\le p+1), z_{\nu s}= y_{\nu-p-1\,k_{s-p-1}}\ (p+2\le \nu\le 2p+1),\\
&&{}\qquad \qquad \qquad \qquad \qquad \qquad \qquad \qquad \qquad \qquad 
{\rm for}\ (p+2\le s\le 2p+1).
\ee

Denote by $Y\!
\left(
\begin{array}{ccccc}
 \nu_1 &\ldots   &\nu_\alpha   \\
 s_1 &\ldots   &s_\alpha   \\
\end{array}
\right)
$ the subdeterminat corresponding to the $\nu_1,\ldots,\nu_\alpha$th row
and $s_1,\ldots,s_\beta$th columns of $Y.$
Laplace expansion as to the first $p$ columns times the last $p+1$ ones of the determinant
of $Z$
equals the expansion as to the first $p+1$ rows times the last $p$ rows. 
This equality shows the following equality
\be
&&\sum_{\nu=1}^{p+1} (-1)^{p+1-\nu}\,Y\!
\left(
\begin{array}{ccccc}
\mu_1  &\ldots&\widehat{\mu_\nu}\ldots   &\mu_{p+1}   \\
 n &k_1& \ldots  & k_{p-1}  \\
\end{array}
\right)\,Y\!
\left(
\begin{array}{ccccc}
 \star &k_1   &\ldots&k_p   \\
 \nu & k_1  &\ldots&k_p   \\
\end{array}
\right)\\
&&= Y\!
\left(
\begin{array}{cccccc}
  \mu_1&\mu_2&\ldots   &\mu_p&\mu_{p+1}   \\
 n &   k_1&\ldots&k_{p-1}&\star   \\  
\end{array}
\right)\,Y(k_1\ldots k_p)\\
&&{}\quad + (-1)^p\,
 Y\!
\left(
\begin{array}{cccccc}
  \mu_1&\mu_2&\ldots   &\mu_p&\mu_{p+1}   \\
 n &   k_1&\ldots&k_{p-1}&k_p   \\  
\end{array}
\right)\,
Y\!
\left(
\begin{array}{ccccc}
  k_1&k_2&\ldots   &k_p   \\
 \star &k_1&\ldots   &k_{p-1}   \\
\end{array}
\right),
\ee
which is nothing else than {\rm(63)}.
\qed
\end{pf}

\medskip

The following Lemma is elementary.
%lm42
\begin{lm}
For $L \subset N_m,\,|L| \le n$ and $j,k ,n \in L^c$ such that $j\ne n, k\ne n,\,j\ne k,$
\beq
&&\frac{B\!
\left(
\begin{array}{cccc}
 0 &n   &k& L   \\
  0& j  &k& L   \\
\end{array}
\right)
B\!
\left(
\begin{array}{ccccc}
 0 &\star   & L   \\
  0& k& L   \\
\end{array}
\right)}{B(0\,k\,L)}
+
\frac{B\!
\left(
\begin{array}{ccccc}
 0 &n   &j& L   \\
  0& k  &j& L   \\
\end{array}
\right)
B\!
\left(
\begin{array}{ccccc}
 0 &\star   & L   \\
  0& j& L   \\
\end{array}
\right)
}{B(0\,j\,L)}\nonumber\\
&&= \frac{B\!
\left(
\begin{array}{cccc}
 0 &\star   &k& L   \\
  0& j  &k& L   \\
\end{array}
\right)
B\!
\left(
\begin{array}{ccccc}
 0 &n   & L   \\
  0& k& L   \\
\end{array}
\right)}{B(0\,k\,L)}
+
\frac{B\!
\left(
\begin{array}{ccccc}
 0 &\star   &j& L   \\
  0& k  &j& L   \\
\end{array}
\right)
B\!
\left(
\begin{array}{ccccc}
 0 &n   & L   \\
  0& j& L   \\
\end{array}
\right)
}{B(0\,j\,L)},\nonumber\\
\eeq
namely, the sum is invariant  under the transposition between
the symbols $n$ and $\star.$
\end{lm}

\bigskip

%proof of Prop39
%\vspace{.5cm}
[{\bf proof of Proposition 39}] 

Let $K_\nu \,(1\le \nu\le p)$ be fixed.
Define $\xi_p(M_p)$ and $\xi_{p}^*(M_p)$
for $1\le p\le n$ by
\be
&&\xi_p(M_p) = B\!
\left(
\begin{array}{cccccc}
  0&M_p&&n+1   \\
 0 &K_{p-1}&n&n+1   \\
\end{array}
\right)
\prod_{\nu=1}^{p-1} B\!
\left(
\begin{array}{cccccc}
  0&K_{\nu-1} &\star&n+1  \\
  0&K_\nu&&n+1   \\
\end{array}
\right),\\
&&
\xi_p^*(M_p) =  B\!
\left(
\begin{array}{cccccccc}
  0&M_p&&&n+1   \\
 0 &K_{p-2}&n&\star&n+1   \\
\end{array}
\right)
\prod_{\nu=1}^{p-2} B\!
\left(
\begin{array}{cccccc}
  0&K_{\nu-1} &\star&n+1  \\
  0&K_{\nu}&&n+1   \\
\end{array}
\right),\\
\ee
and the alternating sum
\be
&&\zeta_{p+1}(M_{p+1})
= \sum_{\nu=1}^{p+1} (-1)^{p+1-\nu} \xi_{p}(\partial_{\mu_\nu} M_{p+1})
B\!
\left(
\begin{array}{cccccc}
  0&K_p &\star&n+1  \\
  0&K_p&\mu_{\nu}&n+1   \\
\end{array}
\right)\\
&&= \frac{1}{p!} \sum_{\sigma\in \mathfrak{S}_{p+1}}\,\mbox{sign}(\sigma) \,\xi_p\bigl(\sigma(\partial_{\mu_1}M_{p+1} )\bigr)
\,B\!
\left(
\begin{array}{cccccc}
  0&K_p &\star&n+1  \\
  0&K_p&\sigma(\mu_1)&n+1   \\
\end{array}
\right).
\ee

Notice that $\xi_{p+1}(M_{p+1}),\xi_{p+1}^*(M_{p+1}) $ and $\zeta_{p+1}(M_{p+1})$ all are skew symmetric 
with respect to
$\mu_1,\ldots,\mu_{p+1}.$

According to {\rm Lemma 41}, the following recurrence relation holds true:
\beq
\zeta_{p+1}(M_{p+1}) = 
\xi_{p+1}(M_{p+1}) + B(0 K_p\,n+1)\,\xi_{p+1}^*(M_{p+1}).
\eeq

Moreover, $\xi_{p+1}^*(M_{p+1})$ is symmetric with respect to
the transposition $\sigma_{k_{p}k_{p+1}}$ 
between $k_p$ and $k_{p+1}:$ 
\beq
\sigma_{k_{p}k_{p+1}}\xi_{p+1}^*(M_{p+1}) = \xi_{p+1}^*(M_{p+1}) .
\eeq

From (56) and (61), $V_n$ can be rewritten as alternating sum:
\be
&&V_n = \sum_{\{k_1,\ldots,k_n\}\subset \partial_{n+1}N_m} 
\frac{\prod_{\nu=1}^n  \lambda_{k_\nu}}{\prod_{\nu=1}^n (\lambda_\infty + n - \nu)}
\frac{B\!
\left(
\begin{array}{ccccc}
 0 &K_n   &n+1   \\
 0 &\partial_{n+1}N  &n+1   \\ 
\end{array}
\right)
}{B(0\,N)\,B(0\,K_n\,n+1)}\\
&&\cdot \sum_{\sigma\in \mathfrak{S}_n} \mbox{sgn}(\sigma)\,\frac{
\xi_1(\sigma(1))\,\prod_{\nu=2}^{n} B\!
\left(
\begin{array}{ccccc}
 0 & \star  &K_{\nu-1}  &n+1 \\
 0 &\sigma(\nu)&K_{\nu-1}   &n+1   \\
\end{array}
\right)
}{\prod_{\nu=1}^{n-1} B(0\,K_\nu\,n+1)}
\,W_0(K_n\,n+1)\varpi,
\ee
where
$\sigma$ ranges over all permutations of  the symbols $1,2,\ldots,n$.

According to (65) and (66), we have the following recurrence relation
\be
&&\sum_{K_n=\{k_1,\ldots,k_n\}\subset N_m} 
\frac{1}{p!}\,\sum_{\sigma\in \mathfrak{S}_n}
\mbox{sgn}(\sigma)\frac{\xi_p(\sigma(M_p)) \prod_{\nu = p+1}^n B\!
\left(
\begin{array}{ccccc}
  0&\star   &K_{\nu-1}  &n+1 \\
 0 &\sigma(\mu_\nu)   &K_{\nu-1}  &n+1 \\
\end{array}
\right)
}{\prod_{\nu=1}^{n-1} B(0\,K_\nu\,n+1)}\\
&&=
\sum_{\{k_1,\ldots,k_n\}\subset N_m} \sum_{\sigma\in \mathfrak{S}_n}
\frac{1}{p!\,(p+1)!}\sum_{\sigma^\prime\in \mathfrak{S}_{p+1}}
\mbox{sgn}(\sigma\sigma^\prime)\\
&&
{}\quad \cdot \frac{\xi_p(\sigma\sigma^\prime(M_p) \prod_{\nu = p+1}^n B\!
\left(
\begin{array}{ccccc}
  0&\star   &K_{\nu-1}  &n+1 \\
 0 &\sigma\sigma^\prime(\mu_\nu)   &K_{\nu-1}  &n+1 \\
\end{array}
\right)
}{\prod_{\nu=1}^{n-1} B(0\,K_\nu\,n+1)}\\
&&=\sum_{K_n=\{k_1,\ldots,k_n\}\subset N_m} 
\,\sum_{\sigma\in \mathfrak{S}_n}
\frac{1}{(p+1)!}\,
\frac{\xi_{p+1}(\sigma(M_{p+1}) \prod_{\nu = p+2}^n B\!
\left(
\begin{array}{ccccc}
  0&\star   &K_{\nu-1}  &n+1 \\
 0 &\sigma(\mu_\nu)   &K_{\nu-1}  &n+1 \\
\end{array}
\right)
}{\prod_{\nu=1}^{n-1} B(0\,K_\nu\,n+1)}\\
&&{}\qquad \qquad \qquad \qquad \qquad \qquad \qquad \qquad \qquad \qquad 
\qquad \qquad \qquad (1\le p\le n-1).
\ee

Hence in view of $\xi_n(\sigma(M_n)) = \mbox{sgn}(\sigma)\,\xi_n(M_n)\,(\sigma\in \mathfrak{S}_n),$
\beq
&&V_n = \sum_{K_n=\{k_1,\ldots,k_n\}\subset \partial_{n+1}N_m} 
\frac{\prod_{\nu=1}^n \lambda_{k_\nu}}{\prod_{\nu=1}^n (\lambda_\infty + n - \nu)}
\,\frac{B\!
\left(
\begin{array}{ccccc}
 0 &K_n   &n+1   \\
 0 &\partial_{n+1}N  &n+1   \\ 
\end{array}
\right)
}{B(0\,N)\,B(0\,K_n\,n+1)}\nonumber\\
&&{}\qquad \cdot\frac{\xi_n(M_n)}
{\prod_{\nu=1}^{n-1} B(0\,K_\nu\,n+1)}.\nonumber\\
\eeq

Furthermore, $J\subset \partial_{n+1}N_m\,(|J|=n)$ being fixed,
repeated application of {\rm(62)} to the RHS of (67) gives
\be
&&\sum_{K_n\equiv J,\,K\subset \partial_{n+1}N_m} 
\frac{B\!
\left(
\begin{array}{ccccc}
 0 &K_n   &n+1   \\
 0 &\partial_{n+1}N  &n+1   \\ 
\end{array}
\right)
}{B(0\,N)\,B(0\,K_n\,n+1)}
\frac{\xi_n(M_n)}{\prod_{\nu=1}^n B(0\,K_\nu\,n+1)}\\
&&= \sum_{K_n\equiv J} \frac{
B\!
\left(
\begin{array}{cccccc}
 0 &K_{n-1}   &n&n+1   \\
  0&K_n &&n+1   \\
\end{array}
\right)
\prod_{\nu=1}^{n-1}B\!
\left(
\begin{array}{cccc}
  0&K_{\nu-1}   &\star&   n+1\\
 0 &K_\nu   &&n+1  \\   
\end{array}
\right)
}
{\prod_{\nu=1}^n\,B(0\,K_\nu\, n+1)}\\
&&=
\sum_{K_n\equiv J}
\frac{\prod_{\nu=2}^{n}B\!
\left(
\begin{array}{ccccc}
  0&K_{\nu-1}   &\star  &n+1 \\
 0 &K_\nu   &&n+1  \\   
\end{array}
\right)
B\!
\left(
\begin{array}{ccccc}
 0 &n& n+1   \\
 0 &k_1   &n+1   \\
\end{array}
\right)
}
{\prod_{\nu=1}^n\,B(0\,K_\nu\, n+1)},
\ee
owing to {\rm  Lemma 40} and {\rm  Lemma 42}.
In this way, (59) and therefore Proposition 39 have been completely proved.
\qed

\bigskip

By using Propositionn 39,
we can conclude the following generalization of Theorem 1, 
which can be deduced in the same procedure
as Theorem 1.

\bigskip
%\vspace{.5cm}

%th3
[{\bf Theorem 3}] \quad We have
\beq
(2\lambda_\infty + n) \varpi \sim \sum_{\nu=1}^{n+1} (-1)^\nu
\sum_{J\subset N_m, |J| = \nu}\frac{\prod_{j\in J}\lambda_j}{\prod_{s=1}^{\nu-1} (\lambda_\infty + n - s)}
W_0(J)\,\varpi,
\eeq
{\it where $J$ ranges over the set of {\rm(}unordered{\rm )} subsets of $N_m$ such that $1\le |J|\le n+1.$}

\begin{pf}
Indeed, the analog of (22) holds true:
\be
(2\lambda_\infty + n) \varpi \sim \sum_{j=1}^n \lambda_j \mathfrak{M}_{f_j-f_{n+1}} F_j + \sum_{j=1}^{m} \lambda_j
(\alpha_{j0}+\alpha_{n+1\,0}) F_j.
\ee
By substituting the formula {\rm(57)} into the above $\mathfrak{M}_{f_j-f_{n+1}} F_j,$  
we obtain (68).\qed
\end{pf}

%\vspace{1cm}
\bigskip

Before proving Theorem 4, it is useful to notify the following Lemma.

%le43
\begin{lm}
Suppose that, for some constants $u_J \in {\bf C}\,(J\,\mbox{admissible}, p = |J|)$
independent of $\lambda,$
\be
\sum_{J \, admissible} \frac{\prod_{j\in J}\lambda_j}{\prod_{s=1}^p\,(\lambda_\infty + n -s)}\,u_J\,W_0(J)\varpi \,\sim \,0.
\ee
Then all $u_J = 0.$
\end{lm}

\begin{pf}
As is seen from {\rm(50)},
$W_0(K)\varpi\;(K \in {\cal B})$ make as representative a basis of $H_\nabla^n(X,\Omega^\cdot(*S)),$
i.e., for any admissible $J,$
\be
W_0(J)\varpi = \sum_{K\in {\cal B}} v_{KJ} W_0(K)\varpi,
\ee
where $v_{KJ}$ are all independent of $\lambda.$
Hence for each $K\in {\cal B},$ we have identically with respect to $\lambda$
\be
\sum_{J \, admissible} \frac{\prod_{j\in J}\lambda_j}{\prod_{s=1}^p\,(\lambda_\infty + n -s)}\,u_J\,v_{KJ} = 0.
\ee
\end{pf}
Hence $u_J\,v_{KJ} = 0.$ For each admissible $J$ there exists 
$K \in {\cal B}$ such that $v_{KJ} \ne 0.$
This means $u_J = 0.$\qed

A generalization of Theorem 2 for general $m$
can be stated as follows:

%\vspace{.5cm}
\bigskip

%th4
[{\bf Theorem 4}]  \quad We have
{\it
\beq
\nabla_B(\varpi) \sim
 \sum_{p=1}^{n+1} 
\sum_{J\subset N_m, |J| = p}\frac{\prod_{j\in J}\lambda_j}{\prod_{s=1}^{p-1} (\lambda_\infty + n - s)}
\,\theta_J\,W_0(J) \,\varpi,
\eeq
where $J \subset N_m$ move over the family of {\rm (unordered)} subsets of indices 
such that $|J| = p$. $\theta_J$ represents a $1$-form 
defined by 
\be
&&\theta_j = - \frac{1}{2}\,d\log \,B(0\star j)
 = - \frac{1}{2}
  \,d\log\,r_j^2,
  \quad(J = \{j\},\,{\rm i.e.,}\, |J| = 1)
  \nonumber\\
&&\theta_{jk} = \frac{1}{2}\,d\log\,B(0jk) = \frac{1}{2}\, d\log\rho_{jk}^2,
\quad(J = \{j,k\},\,{\rm i.e.,}\,|J| = 2)\nonumber\\
&&\theta_J = (-1)^p \sum_{j,k\in J, j<k }\,\frac{1}{2}\, d\log B(0jk) \nonumber\\
&& \sum_{\{\mu_1,\ldots,\mu_{p-2}\}\equiv \partial_j\partial_k J  }\prod_{s=1}^{p-2}
 \frac{
B\!
\left(
\begin{array}{cccccccc}
 0 &\star   &\mu_{s-1}&\cdots&\mu_1&j&k   \\
 0 & \mu_s  &\mu_{s-1}& \cdots&\mu_1&j&k  
\end{array}
\right)
}
{B(0\mu_s\,\mu_{s-1}\cdots\mu_1\,j\,k)},\quad(|J|= p \ge 2)\nonumber\\
\ee
where $M = \{\mu_1,\ldots,\mu_{p-2}\}$ move over the family of all ordered
sets consisting of  $p-2$ different elements of $\partial_j\partial_k J$.
}
%endthm4

\begin{pf}
Lemma 43 shows that (69) is a unique expression.
We have the analog of (31) as follows:
\be
&&\nabla_B(\varpi) \sim \sum_{j=1}^{n+1}
 \sum_{\nu=0}^{n+j-\nu}\,d\alpha_{j\nu}
\frac{\partial}{\partial \alpha_{j\nu}}{\rm log}\Phi \varpi
%\nabla_{B, \frac{\partial}{\partial \alpha_{j\nu}}} \varpi
+ \sum_{j=n+2}^{m} \sum_{\nu=0}^{n+1}
d\alpha_{j\nu}
\frac{\partial}{\partial \alpha_{j\nu}}{\rm log}\Phi \varpi.
%\nabla_{B, \frac{\partial}{\partial \alpha_{j\nu}}
%}
%\varpi\\
 \ee

 In addition to the differential $1$-forms $\theta_k^j\,(1\le j\le k\le n)$ defined in 
 {\bf{\rm Definition 26}}, define further
 the differential $1$-forms $\theta_k^j \,(n+2\le j \le m,\,1\le k\le n)$
 by the following reccurrence relations: 
 \be
&& \theta_{k}^j = \frac{1}{\alpha_{k\,n+1-k}}\{d\alpha_{j,\,n+1-k} - \sum_{s=1}^{k-1} \alpha_{s\,n+1-k}\,\theta_s^j\},\quad(n+2\le j\le m,\,1\le j\le k\le n)\\
&&\theta_{1}^j  = \frac{1}{\alpha_{1\,n}}\,d\alpha_{j\,n}.
 \ee
 
 Then according to Proposition 40, we get 
  \be
 &&\nabla_B(\varpi) = \sum_{j=1}^n \lambda_j 
 \bigl( \sum_{k= j+1}^n \,\theta_k^j \,\frac{f_k - f_j}{f_j}\varpi
 - \sum_{k=1}^{n}\, \theta_k^j\,\frac{f_{n+1}- f_j}{f_j}\varpi\bigr) \\
 &&{}\qquad \quad \qquad+ 
 \sum_{j=n+2}^m \lambda_j 
 \bigl( \sum_{k= 1}^n \,\theta_k^j \,\frac{f_k - f_j}{f_j}\varpi
 - \sum_{k=1}^{n}\, \theta_k^j\,\frac{f_{n+1}- f_j}{f_j}\varpi\bigr) \\
&&{}\qquad \quad  = \sum_{\nu=1}^{n+1} \sum_{J\subset N_m, |J|=\nu} \frac{\prod_{j\in J}\lambda_j}
{\prod_{s=0}^{\nu-1} (\lambda_\infty + n-s)} \,\tilde{\theta}_J \,W_0(J)\varpi,
\ee
where the differential $1$-form $\tilde{\theta}_J$ does not depend on $\lambda.$ 
When we specify $J$ such that $J\subset N,$ then $\tilde{\theta}_J$ coincides with
$\theta_J$ owing to {\rm Theorem 2}. Hence, because of symmetry and uniqueness, we can conclude that $\tilde{\theta}_J$
also coincides with $\theta_J$ for every $J$ such that $J\subset N_m.$
\qed
\end{pf}
%*
%*

\bigskip

%sec7
\section{Gauss-Manin connection}

Take a non empty (unordered ) subset of indices
$J$ such that $J \subset N_m = \{1,2,\ldots, m\}, |J| =p \le n+1$
(called \lq\lq admissible" element).

$T_{-\varepsilon_j}F_J$ and $\nabla_B F_J\,$
can be described in integrable matrix  form:
\beq
&&T_{-\varepsilon_j}F_J 
\sim \sum_{K ; admissible} \tilde{\gamma}_{K,J}^j\, F_K,\\
&&\nabla_B F_J  \sim \sum_{K : admissible} F_K\,\Theta_{K,J},
\eeq
where $F_J$ and $W_0(K)\varpi$ are connected by (3)
with the inverse of the transition matrix $(\tilde{\beta}_{K,J})$
between the two system of admissible elements $(F_J)$ and $(W_0(J)\varpi)$.

$\Theta_{K,J}$ denote  rational 1-forms in terms of $\theta_L\,(L : {\rm admissible})$.

Although it is generally very complicated,
at least in principle we can find out the variation formula $\nabla_B F_J$
by the use of Theorem 2, Theorem 4 and Theorem A.

The basic equation for it is related to (2) and can be stated as:
\be
d_B\,T_{-\varepsilon_J} {\cal J}_\lambda(\varphi) = T_{-\varepsilon_J}
\,d_B{\cal J}_\lambda(\varphi),
\ee
Indeed, Theorem 2 and Theorem 4 show
\be
&&\nabla_B F_J =  \nabla_B\,\mathfrak{M}_{f_J^{-1}}\varpi
\,\sim 
\mathfrak{M}_{f_J^{-1}}T_{-\varepsilon_J}\nabla_B \varpi\quad(\varepsilon_J=
\sum_{j\in J}\varepsilon_j )\\
&& \sim 
\sum_{p=1}^{n+1}\, \sum_{|K| = p; K : admissible} \,T_{-\varepsilon_J}\left(\frac{
\lambda_J}
{\prod_{\nu = 1}^{p-1} (\lambda_\infty + n - \nu )}\right)\,
\mathfrak{M}_{f_J^{-1}}W_0(K)\varpi\,\theta_{K}. \\
\ee

On the other hand, by seeing that
\be
&&\mathfrak{M}_{f_J^{-1}}T_{-\varepsilon_J}\left(g(\lambda)F_K\right)
= T_{-\varepsilon_J}(g(\lambda))\,\mathfrak{M}_{f_J^{-1}}F_K,\\
&& \mathfrak{M}_{f_J^{-1}f_K^{-1}}T_{-\varepsilon_J-\varepsilon_K}= \mathfrak{M}_{f_J^{-1}}T_{-\varepsilon_J}\,
\mathfrak{M}_{f_K^{-1}}T_{-\varepsilon_K},
\ee
and by applying Theorem A, 
$\nabla_BF_J$
can be explicitly written by the use of of $\theta_K$ in a recursive way.

Assume first the case where $p=1$ and $J = \{j\}$. Then
\beq
&& \sum_{K : admissible} F_K \,\Theta_{K,j}
\sim\nabla_B F_j = T_{-\varepsilon_j}\mathfrak{M}_{f_j^{-1}}\nabla_B\varpi\nonumber\\
&&
= 
\sum_{p=1}^{n+1}\, \sum_{|L| = p} 
T_{-\varepsilon_j}\left(\frac{
\lambda_L}
{\prod_{\nu = 1}^{p-1} (\lambda_\infty + n - \nu)}\right)\,
\mathfrak{M}_{f_j^{-1}}W_0(L)\varpi\,\theta_{L}\nonumber\\
&&
= \sum_{p=1}^{n+1}\, \sum_{|L| = p, j \in L} 
\,\left(
\frac{
\lambda_{\partial_jL}
}
{\prod_{\nu = 1}^{p-1} (\lambda_\infty + n - \nu -1)
}
\right)
\,
W_0^{(j)}(L)\varpi\,\theta_{L}\nonumber\\
&&+
\sum_{p=j}^{n+1}\, \sum_{|L| = p, j\notin L } 
\left(\frac{
\lambda_L}
{\prod_{\nu = 1}^{p-1} (\lambda_\infty + n - \nu - 1 )}\right)\,
\mathfrak{M}_{f_j^{-1}}W_0(L)\varpi\,\theta_{L}.
\eeq

Since 
\be
&&W_0^{(j)}(L)\varpi \sim \sum_{K} \,\gamma_{K,L}^j\,F_K\,\,\quad(j\in L)\quad({\rm see}\,(A.2)),\\
&&\mathfrak{M}_{f_j^{-1}}W_0(L)\varpi = \sum_{K}\,\frac{\gamma_{K,L}^j}{\lambda_j - 1}\,F_K\\
&&{}\quad \quad  = - \sum_{k\in L} \,
B\!
\left(
\begin{array}{ccc}
 0 &\star   &\partial_kL   \\
  0&k   &\partial_kL      
\end{array}
\right)\,F_{j\,\partial_kL} 
+ B(0\star L)\,F_{j\,L}\quad(j\notin L),
\ee
we may put
\beq
&&\Theta_{K,j} = \sum_{p=1}^{n+1}\sum_{|L| = p, j\in L},
\,\left(
\frac{
\lambda_{\partial_jL}
}
{\prod_{\nu = 1}^{p-1} (\lambda_\infty + n - \nu - 1)
}
\right)\,\gamma_{K,L}^1\,\theta_L\nonumber\\
&&{}\quad \quad  + 
\sum_{p=1}^{n+1}\sum_{|L| = p, j\notin L}
\,\left(
\frac{
\lambda_{L}
}
{\prod_{\nu = 1}^{p-1} (\lambda_\infty + n - \nu - 1)
}
\right)\,\frac{\gamma_{K,L}^1}{\lambda_j - 1}\,\theta_L.
\eeq

\medskip

%prop44
\begin{prop}
Fix an  admissible and unordered subset $J \in N_m\,(|J| = p,\,2\le p\le n+1)$ and  fix $j \in J$. Then
$\Theta_{K,J}$ in {\rm(71)} can be expressed as
\beq
\Theta_{K,J } =\sum_{L : admissible} \frac{{\tilde{\gamma}}_{K,L}^j}{\lambda_j - 1}\, T_{-\varepsilon_j}(\Theta_{L, \partial_jJ}).
\eeq
\end{prop}

\begin{pf}
In fact,
\be
&&\nabla_B \,F_J = \nabla_B(\mathfrak{M}_{f_j^{-1}}\,F_{\partial_jJ})\\
&&{}\qquad \quad \sim \sum_{L:admissible}\,T_{-\varepsilon_j}(\Theta_{L,\partial_jJ})\,\mathfrak{M}_{f_j^{-1}}F_L\\
&&{}\qquad \quad = \sum_{K:admissible}\sum_{L : admissible} \, T_{-\varepsilon_j}(\Theta_{L, \partial_jJ})\,\frac{{\tilde{\gamma}}_{K,L}^j}{\lambda_j - 1}\,F_K.
\ee
This implies (74). \qed
\end{pf}

In this way, we conclude the following.

\bigskip

%\vspace{.5cm}
%Theorem 5

{\bf[Theorem 5]}

{\it With respect to the admissible system $\{F_J ( = \frac{\varpi}{f_J})\}$ of $H_\nabla^n(X, \Omega^\cdot(*S))$,
\be
d_B {\cal J}_\lambda \bigl(\frac{1}{f_J}\bigr) = \sum_{K: admissible} \, {\cal J}_\lambda \bigl(\frac{1}{f_K}\bigr) \,\Theta_{K,J}
\ee
preresents the Gauss-Manin connection for the integrals ${\cal J}_\lambda(\varphi)$. Each element of the matrix {\rm 1}-form
$\Theta_{K,J}$ is a rational {\rm1}-form in the parameters $r_j^2, \rho_{kl}^2$
having poles along the loci $B(0\,J) = 0$ and $B(0\star J) = 0.$
It is rational in $\lambda$ having poles only on the set 
\be
\lambda_\infty + n = 1,2,3,\ldots .
\ee
whose numerator is a polynomial of degree $p\,(|J| = p).$ }

\bigskip

%\vspace{.5cm}

Proposition 44 means the following formula. Namely
when $J = \{j,k\}$,
applying the operations $T_{-\varepsilon_k}\,(k\ne j)$ on both sides of (72),
we can evaluate $\Theta_{K,jk}$ in a successive way
and derive the following recurrence formula
\beq
\Theta_{K, jk} = \frac{1}{\lambda_j- 1}\, \sum_{L : admissible}
\,T_{-\varepsilon_j}(\Theta_{L,k})
\,\tilde{\gamma}_{K,L}^j,
\eeq
which is symmetric with respect to $j,\,k$.

Theorem 5 holds true for $1\le m\le n+1,$ as well as Theorem 1 and Theorem 2 do so. 
We give explicit formulae for $\Theta_{K,J}$  in two simplest non trivial cases.

%ex
%
\bigskip
%\vspace{.5cm}

(i){\bf [Case $n\ge 1,\; m = 2$ ]}

\medskip

In this case, the system of admissible elements consists of $\{F_1, F_2, F_{12}\}$ which
coincides with the NBC basis. 

The following Lemmas are a direct consequence from Theorem A (see (A.12) and (A.23)).

%llm45
\begin{lm}
For $\{j,k\} = \{1,2\},$
\be
{\rm (i)} 
&&\gamma_{j,j}^j = - ( \lambda_\infty + \lambda_j + n - 2),\; \gamma_{k,j}^j = - \lambda_k,\\
&&\gamma_{jk,j}^j = - \lambda_k\,B\!
\left(
\begin{array}{ccc}
  0&\star   &k   \\
  0& \star  &j   
\end{array}
\right),\\
{\rm(ii)}
&&\gamma_{jk,k}^j = (\lambda_j - 1)\,B(0\star k),\\
{\rm(iii)}
&&\gamma_{j,jk}^j = \lambda_\infty + n - 2 , \;\gamma_{k,jk}^j = - (\lambda_\infty + n - 2),\\
&&\gamma_{jk,jk}^j = (\lambda_\infty + n - 2)\,
\,B\!
\left(
\begin{array}{ccc}
  0&\star   &k   \\
  0&j  &k   
\end{array}
\right).
\ee
\end{lm}

%lm46
\begin{lm}
For $\{j,k\} = \{1,2\},$
\be
&&\tilde{\gamma}_{j,j}^j = - \frac{(\lambda_\infty + \lambda_j + n -2)}{B(0\star j)},\;
\tilde{\gamma}_{k,j}^j = - \frac{\lambda_k}{B(0\star j)},\\
&&\tilde{\gamma}_{jk,j}^j = - \lambda_k
 \frac{B\!
\left(
\begin{array}{ccc}
  0&\star   &k   \\
  0&\star  &j   
\end{array}
\right)}{B(0\star j)},\\
&&\tilde{\gamma}_{jk,k}^j = \lambda_j - 1,\\
&&\tilde{\gamma}_{j,jk}^j = (\lambda_\infty + n - 2)\,
\frac{B\!
\left(
\begin{array}{ccc}
  0&\star   &j   \\
  0&\star &k   
\end{array}
\right)}{B(0\star j\,k)\,B(0\star j)}
- \lambda_j\,\frac{
B\!
\left(
\begin{array}{ccc}
  0&\star   &j   \\
  0&k &j  
\end{array}
\right)}{B(0\star j\,k)\,B(0\star j)},\\
&&\tilde{\gamma}_{k,jk}^j = - \lambda_k\,\frac{B\!
\left(
\begin{array}{ccc}
  0&\star   &j   \\
  0&k &j   
\end{array}
\right)}{B(0\star j\,k)\,B(0\star j)} - \frac{(\lambda_\infty + n - 2)}{B(0\star j\,k)},\\
&& \tilde{\gamma}_{jk,jk}^j = - \lambda_k \,\frac{B\!
\left(
\begin{array}{ccc}
  0&\star   &k   \\
  0&\star &j  
\end{array}
\right)\, B\!
\left(
\begin{array}{ccc}
  0&\star   &j   \\
  0&k &j  
\end{array}
\right)}{B(0\star j\,k)\,B(0\star j)} + (\lambda_\infty + \lambda_j + n - 3)
\frac{ 
B\!
\left(
\begin{array}{ccc}
  0&\star   &k   \\
  0&j &k
\end{array}
\right)}{B(0\star j\,k)}.
\ee
\end{lm}

Because of Lemmas 45 and Lemma 46, (73) and (74) give the explicit formulae for
$\Theta_{K,J}.$

%prop47
\begin{prop}
For $\{j,k\} = \{1,2\},$
\be
&&\Theta_{j,j} = - (\lambda_\infty + n - 2 + \lambda_j)\,\theta_j + \lambda_k\,\theta_{jk},\\
&&\Theta_{k,j} = - \lambda_k \,(\theta_j + \theta_{jk}),\\
&&\Theta_{jk,j} =  \lambda_k\, \{- B\!
\left(
\begin{array}{ccc}
  0&\star   &k   \\
  0&\star &j
\end{array}
\right)\,\theta_j + B(0\star k)\,\theta_k 
+ B\!
\left(
\begin{array}{ccc}
  0&\star   &k   \\
  0&j &k
\end{array}
\right)\,\theta_{jk}\},\\
&&\Theta_{j,jk} = - \lambda_j\,\{
\frac{B\!
\left(
\begin{array}{ccc}
  0&\star   &j   \\
  0&k &j
\end{array}
\right)}{B(0\star j\,k)}\,\theta_j + 
\frac{B\!
\left(
\begin{array}{ccc}
  0&\star   &k   \\
  0&j &k
\end{array}
\right)}{B(0\star j\,k)}\,\theta_k +
\frac{B(0\,j\,k)}{B(0\star j\,k)}\,\theta_{jk}\}\\
&&{}\quad  + (\lambda_\infty + n - 2) \{\frac{B\!
\left(
\begin{array}{ccc}
  0&\star   &j   \\
  0&\star &k
\end{array}
\right)}{B(0\star j\,k)}\,\theta_j -
\frac{B(0\star k))}{B(0\star j\,k)}\,\theta_k -
\frac{B\!
\left(
\begin{array}{ccc}
  0&\star   &k   \\
  0&j &k
\end{array}
\right)}{B(0\star j\,k)}\,\theta_{jk}\},
\ee
\be
&&\Theta_{12,12} =\\
&& \{- \lambda_2 \,\frac{B\!
\left(
\begin{array}{ccc}
  0&\star   &2  \\
  0&\star &1
\end{array}
\right)\,B\!
\left(
\begin{array}{ccc}
  0&\star   &1  \\
  0&2 &1
\end{array}
\right)}{B(0\star 1\,2)} + (\lambda_\infty + \lambda_1 + n -3)\, 
\frac{B\!
\left(
\begin{array}{ccc}
  0&\star   &2   \\
  0&1 &2
\end{array}
\right)\,B(0\star 1)}{B(0\star 1\,2)}\}\,\theta_1 \\
&&+\{- \lambda_1 \,\frac{B\!
\left(
\begin{array}{ccc}
  0&\star   &1  \\
  0&\star &2
\end{array}
\right)\,B\!
\left(
\begin{array}{ccc}
  0&\star   &2  \\
  0&1 &2
\end{array}
\right)}{B(0\star 1\,2)} + (\lambda_\infty + \lambda_2 + n -3)\, 
\frac{B\!
\left(
\begin{array}{ccc}
  0&\star   &2   \\
  0&1 &2
\end{array}
\right)\,B(0\star 2)}{B(0\star 1\,2)}\}\,\theta_2 \\
&&{}\quad + \{- (\lambda_\infty - 1)\,\bigl(2\frac{B\!
\left(
\begin{array}{ccc}
  0&\star   &1  \\
  0&\star &2
\end{array}
\right)B(0\,1\,2)}{B(0\star 1\,2)} + 1\bigr) \\
&&{}\quad - (n-1) \,\bigl(\frac{B\!
\left(
\begin{array}{ccc}
  0&\star   &1  \\
  0&\star &2
\end{array}
\right)B(0\,1\,2)}{B(0\star 1\,2)} + 1\bigr) \}\,\theta_{12}.
\ee

\end{prop}

\begin{pf}
The first four identities directly follow from the Definition.
 As for $\Theta_{jk,jk}$,  (73) and (74) show 
 \be
&& \Theta_{jk,jk} = \frac{\tilde{\gamma}_{jk,j}^k}{\lambda_k - 1}\,\{- (\lambda_\infty + \lambda_j + n -3)\,\theta_j + (\lambda_k - 1)\,\theta_{jk}\}
- \tilde{\gamma}_{jk,k}^k\,\{\theta_j + \theta_{jk}\} \\
&&{}\qquad \quad + \tilde{\gamma}_{jk,jk}^k\{- B\!
\left(
\begin{array}{ccc}
  0&\star   &k   \\
  0&\star&j
\end{array}
\right)\,\theta_j + B(0\star k)\,\theta_k + 
B\!
\left(
\begin{array}{ccc}
  0&\star   &k   \\
  0&j &k
\end{array}
\right)\,\theta_{jk}\}.
 \ee

The RHS is rewritten as in Proposition 47
in view of the Jacobi identities:
\be
&&B(0\star j\,k) + B\!
\left(
\begin{array}{ccc}
  0&\star   &j   \\
  0&k   &j   
\end{array}
\right) \,B\!
\left(
\begin{array}{ccc}
  0&\star   &k  \\
  0&j  &k   \\ 
\end{array}
\right) +  B\!
\left(
\begin{array}{ccc}
  0&\star   &j   \\
  0&\star   &k  
\end{array}
\right) \,B(0\,j\,k)= 0,\\
&&B(0\star j\,k) + B\!
\left(
\begin{array}{ccc}
  0&\star   &k  \\
  0&\star   &j   
\end{array}
\right) \,B\!
\left(
\begin{array}{ccc}
  0&\star   &j  \\
  0&k  &j   \\ 
\end{array}
\right) +  B\!
\left(
\begin{array}{ccc}
  0&\star   &k  \\
  0&j   &k  
\end{array}
\right) \,B(0\star j)= 0.
\ee
\qed
\end{pf}

\medskip

%co48
\begin{co}
The trace of $\Theta$ defined by
\be
T_r(\Theta) = \Theta_{1,1} + \Theta_{2,2} + \Theta_{12,12}
\ee
 can be explicitly written as
 \be
 T_r(\Theta) = d_B\log {\bf W},
 \ee
 where $\bf W$ denotes the Wronskian
  \beq
 {\bf W} = \frac{
 B(0\star 1)^{\lambda_1+ \frac{n-2}{2}}\,B(0\star 2)^{\lambda_2 + \frac{n-2}{2}}
 \,
 B(0\star 1\,2)^{
 \lambda_\infty + \frac{n-3}{2}}
 }
 {
 {B(0\,1\,2)}^{\frac{n-2}{2}}
 }.
 \eeq

\end{co}

\begin{pf}
Indeed due to Proposition 47, the identities
\be
&& \theta_j = - \frac{1}{2}\,d\log r_j^2\;(j = 1,2),\; \theta_{12} = \frac{1}{2}\,d\log\,\rho_{12}^2.\\
&&dB(0\star 1\,2) = - 2B\!
\left(
\begin{array}{ccc}
  0&\star   &1   \\
  0& 2  &1    
\end{array}
\right)\,d\,r_1^2  - 2B\!
\left(
\begin{array}{ccc}
  0&\star   &2   \\
  0& 1  &2    
\end{array}
\right)\,d\,r_2^2  - 2B\!
\left(
\begin{array}{ccc}
  0&\star   &1   \\
  0& \star  &2    
\end{array}
\right)\,d\,\rho_{12}^2
\ee
lead to Corollary 48.
\qed
\end{pf}

\bigskip
%\vspace{.5cm}

(ii){\bf [Case $n =1, \;m\ge 2$]}

\medskip

We compute $\Theta_{K,J}$ in the case of $n = 1$  in an explicit way.
Admissible elements of $H_\nabla^1(X, \Omega^\cdot(*S))$ (of dimension $2m - 1)$ consists 
of $F_j \,(1\le j\le m)$ and $F_{jk}\,(1\le j<k\le m)$. 
The NBC basis is given by $F_j \,(1\le j\le m)$ and $F_{2j}\,(1\le j\le m, j\ne 2)$
(see the Remark in \S6).

In the case $m \ge 3,$ there are fundamental relations among admissible elements. These 
are stated as in (48) and (50):
\beq
&&F_{jkl} =\frac{ B\!
\left(
\begin{array}{cccc}
  0&\star   &k&l   \\
  0&j   &k&l    
\end{array}
\right)}{B(0\,j\,k\,l)}\,F_{kl}
+ \frac{B\!
\left(
\begin{array}{cccc}
  0&\star   &j&l   \\
  0&k   &j&l    
\end{array}
\right)}{B(0\,j\,k\,l)}\,F_{jl}
+ \frac{B\!
\left(
\begin{array}{cccc}
  0&\star   &j&k   \\
  0&l   &j&k    
\end{array}
\right)}{B(0\,j\,k\,l)}\,F_{jk},\nonumber\\
\\
&&W_0(kl)\,\varpi = \frac{B\!
\left(
\begin{array}{ccc}
 0 & j  &l   \\
  0& k  &l  \\
\end{array}
\right)
}{B(0\,j\,l)}\,W_0(jl)\,\varpi
+ \frac{B\!
\left(
\begin{array}{ccc}
 0 & j  &k  \\
  0& l  &k \\
\end{array}
\right)
}{B(0\,j\,k)}\,W_0(jk)\,\varpi,
\eeq
where
 \be
 W_0(jk)\varpi = - B\!
\left(
\begin{array}{ccc}
  0&\star   & k  \\
 0 &j   &k   
\end{array}
\right)\,F_k 
- B\!
\left(
\begin{array}{ccc}
  0&\star   & j  \\
 0 &k   &j   
\end{array}
\right)\,F_j + B(0\star j\,k)\,F_{jk}.
\ee
 
Suppose for simplicity $\alpha_{j1} >  \alpha_{k1}$
such that $\rho_{jk} = - \rho_{k,j} = \alpha_{j1} - \alpha_{k1} > 0$ for $j < k$.
Then these relations give  linear expressions for $F_{1l}\,(3\le l \le m)$
and $F_{kl}\,(3\le k < l \le m)$ by the NBC basis.
\beq
&&W_0(1l)\,\varpi = \frac{\rho_{1l}}{\rho_{2l}}\, W_0(2l)\varpi
 + \frac{\rho_{1l}}{\rho_{12}}\,W_0(12)\varpi\quad(3\le l \le m),\\
 &&W_0(kl)\varpi = \frac{\rho_{kl}}{\rho_{2l}}\, W_0(2l)\varpi
 - \frac{\rho_{kl}}{\rho_{2k}}\,W_0(2k)\varpi\quad(3\le k<l \le m).
\eeq

We put
\be
\Delta_{jkl} = 
\left|
\begin{array}{ccc}
1  &1   &1   \\
 \alpha_{j0} &\alpha_{k0}   & \alpha_{l0}  \\
\alpha_{j1}  & \alpha_{k1}  &\alpha_{l1}   
\end{array}
\right|.
\ee

Then
\beq
&&B\!
\left(
\begin{array}{cccc}
 0 & \star  &k&l   \\
 0 &j   &k&l   \\
\end{array}
\right) = - \rho_{kl}\,\Delta_{jkl},\; 
B(0\star j\,k\,l) = - 2 \Delta_{jkl}^2,\\
&& \Delta_{jkl} = \rho_{jk}\,\rho_{jl}\,\rho_{kl} - \rho_{kl}r_j^2 + \rho_{jl}r_k^2 - \rho_{jk}r_l^2.
\eeq

(78) can be rewritten in terms of $F_j, \,F_{kl}$ in skew symmetric form
\beq
v_j\,F_j + v_k\,F_k + v_l\,F_l + v_{jk}F_{jk} - v_{jl}F_{jl} + v_{kl}F_{kl} = 0
\eeq
with
\be 
v_j = \frac{\Delta_{jkl}}{\rho_{kj}\,\rho_{lj}},\, v_{kl} = \frac{B(0\star kl)}{\rho_{kl}}\quad {\rm etc.}
\ee

The following two Lemmas are a direct consequence from Theorem A
 (see also (A.12) and  (A.23))  and  the Remark next to it (see (A.10) and (A.11)).

%lm49
\begin{lm}
For different indices $j,k,l,$ we have
\be
{\rm(i)}&&\gamma_{j,j}^j = - (\lambda_\infty + \lambda_j - 1),\; \gamma_{k,j}^j = - \lambda_k,\\
&&\gamma_{jk,j}^j = - \lambda_k \,B\!
\left(
\begin{array}{ccc}
 0 &\star   & k  \\
 0 &\star& j     
\end{array}
\right),\\
%\ee
%
%\be
{\rm(ii)}&&\gamma_{jk,k}^j = (\lambda_j - 1) \,B(0\star k), \\
%\ee
%
%\be
{\rm(iii)}&&\gamma_{j,jk}^j = \lambda_\infty - 1,\;
\gamma_{k,jk}^j = - (\lambda_\infty - 1),\\
&&\gamma_{jk,jk}^j = (\lambda_\infty - 1)\,B\!
\left(
\begin{array}{ccc}
 0 &\star   &k   \\
  0&j   &k   
\end{array}
\right) - \sum_{\nu\ne j,k}\,\lambda_\nu\, \frac{
B\!
\left(
\begin{array}{cccc}
 0 &\star   &k&\nu   \\
  0&\star   &k&j  
\end{array}
\right)\,
B\!
\left(
\begin{array}{cccc}
 0 &\star   &j&k   \\
  0&\nu   &j&k   
\end{array}
\right)\,
}{B(0\star \nu\, j\,k)},\\
&&\gamma_{jl,jk}^j = - \lambda_l 
\frac{
B\!
\left(
\begin{array}{cccc}
 0 &\star&j&k      \\
  0&l   &j&k   \\   
\end{array}
\right)\,
B\!
\left(
\begin{array}{cccc}
 0 &\star&j&l      \\
  0&\star   &k&l   \\   
\end{array}
\right)\,
}{B(0\star l\, j\,k)},\\
&&\gamma_{lk,jk}^j =  \lambda_l 
\frac{
B\!
\left(
\begin{array}{cccc}
 0 &\star&j&k     \\
  0&l   &j&k \\   
\end{array}
\right)\,
B(0\star k\,l))\,
}{B(0\star l\, j\,k)}, \\
%\ee
%
%\be
{\rm(iv)}&&\gamma_{jk,kl}^j = - (\lambda_j - 1) 
\frac{
B\!
\left(
\begin{array}{cccc}
 0 &\star&k&j      \\
  0&\star   &k&l   \\   
\end{array}
\right)\,
B\!
\left(
\begin{array}{cccc}
 0 &\star&k&l      \\
  0&j   &k&l   \\   
\end{array}
\right)\,
}{B(0\star l\, j\,k)},\\
&&\gamma_{kl,kl}^j = (\lambda_j -1 ) \frac{
B\!
\left(
\begin{array}{cccc}
 0 &\star&k&l     \\
  0&j  &k&l \\   
\end{array}
\right)\,
B(0\star k\,l)\,
}{B(0\star j\, k\,l)}.
\ee
All other $\gamma_{K,J}^j = 0$.
\end{lm}

In the same way due to (A.3),

%lm50
\begin{lm}
For different indices $j,k,l,$ we have
\be
&&\tilde{\gamma}_{j,j}^j = - \frac{(\lambda_\infty + \lambda_j - 1)}{B(0\star j)},\; \tilde{\gamma}_{k,j}^j = - \frac{\lambda_k}{B(0\star j)},\\
&&\tilde{\gamma}_{jk,j}^j = - \lambda_k \,\frac{B\!
\left(
\begin{array}{ccc}
 0 &\star   & k  \\
 0 &\star& j     
\end{array}
\right)}{B(0\star j)},\\
&&\tilde{\gamma}_{jk,k}^j = \lambda_j - 1, \\
&&\tilde{\gamma}_{j,jk}^j = (\lambda_\infty - 1)\,\frac{B\!
\left(
\begin{array}{ccc}
 0 &\star   &j   \\
 0 &\star   & k  
\end{array}
\right)
}{B(0\star j\,k))\,B(0\star j)} - \lambda_j\,\frac{B\!
\left(
\begin{array}{ccc}
 0 &\star   &j   \\
 0 &k   & j  
\end{array}
\right)
}{B(0\star j\,k))\,B(0\star j)},\\
&&\tilde{\gamma}_{k,jk}^j = -\lambda_k\,\frac{
B\!
\left(
\begin{array}{ccc}
 0 &\star   &j   \\
 0 &k   & j  
\end{array}
\right)
}{B(0\star j\,k))\,B(0\star j)} -
(\lambda_\infty - 1)\,
\frac{1}{B(0\star j\,k)},\\
&&\tilde{\gamma}_{l,jk}^j = - \lambda_l \,\frac{B\!
\left(
\begin{array}{ccc}
 0 &\star   &j   \\
 0 &k   & j  
\end{array}
\right)
}{B(0\star j\,k))\,B(0\star j)}, \\
&&\tilde{\gamma}_{jk,jk}^j = - \lambda_k \frac{B\!
\left(
\begin{array}{ccc}
 0 &\star   &k  \\
 0 &\star   & j  
\end{array}
\right)\,
B\!
\left(
\begin{array}{ccc}
 0 &\star   &j  \\
 0 &k   & j  
\end{array}
\right)}{B(0\star j\,k)\,B(0\star j)} +
(\lambda_\infty + \lambda_j - 2)
\frac{B\!
\left(
\begin{array}{ccc}
 0 &\star   &k \\
 0 &j   & k  
\end{array}
\right)
}{B(0\star j\,k)}\\
&&{}\qquad \quad  - \sum_{\nu\ne j,k}\,\lambda_\nu\,
\frac{B\!
\left(
\begin{array}{cccc}
 0&\star &k&\nu     \\
  0&\star&k   &j     
\end{array}
\right)\,
B\!
\left(
\begin{array}{cccc}
 0&\star &j&k   \\
  0&\nu&j   &k    
\end{array}
\right)}{B(0\star \nu\,j\,k)\,B(0\star j\,k)},\\
&&\tilde{\gamma}_{jl,jk}^j = - \lambda_l\, \bigl\{\frac{B\!
\left(
\begin{array}{cccc}
 0&\star &l    \\
  0&\star&j     
\end{array}
\right)\,
B\!
\left(
\begin{array}{cccc}
 0&\star &j   \\
  0&k&j    
\end{array}
\right)}{B(0\star j\,k)\,B(0\star j)}
+ \frac{B\!
\left(
\begin{array}{cccc}
 0&\star &j&k   \\
  0&l&j&k     
\end{array}
\right)\,
B\!
\left(
\begin{array}{cccc}
 0&\star &j &l  \\
  0&\star&k&l   
\end{array}
\right)}{B(0\star j\,k\,l)\,B(0\star j\,k)}
\},\\
&&\tilde{\gamma}_{kl,jk}^j =  \lambda_l\,\frac{B(0\star k\,l)\,
B\!
\left(
\begin{array}{cccc}
 0 &\star   &j&k   \\
  0&l&j   &k   \\
\end{array}
\right)
}{B(0\star j\,k\,l)\,B(0\star j\,k)},\\
&&\tilde{\gamma}_{jl,kl}^j =(\lambda_j - 1) \, \frac{
B\!
\left(
\begin{array}{cccc}
  0&\star   &j&l   \\
  0&k   &j&l   
\end{array}
\right)
}{B(0\star j\,k\,l)},\\
&&\tilde{\gamma}_{kl,kl}^j = (\lambda_j - 1)\,\frac{
B\!
\left(
\begin{array}{cccc}
  0&\star   &k&l   \\
  0&j   &k&l   
\end{array}
\right)
}{B(0\star j\,k\,l)}.
\ee
All other $\tilde{\gamma}_{K,J}^j = 0.$
\end{lm}

\medskip

Denote for different indices $j,k,l,$ the differential $1$-forms
\be
&&\zeta_{jk} =  B\!
\left(
\begin{array}{ccc}
 0 &\star   &j   \\
 0 &k   &j   \\
\end{array}
\right)\,\theta_j +
B\!
\left(
\begin{array}{ccc}
 0 &\star   &k   \\
 0 &j   &k   \\
\end{array}
\right)\,\theta_k
+ B(0\,j\,k)\,\theta_{jk},\\
&&\zeta_{jk,j} =  -  B\!
\left(
\begin{array}{ccc}
 0 &\star   &k   \\
0  &\star   &j   
\end{array}
\right)\,\theta_j + B(0\star k)\,\theta_k +
B\!
\left(
\begin{array}{ccc}
 0 &\star   &k   \\
0  &j   &k  
\end{array}
\right)\,\theta_{jk},\\
&&\zeta_{j\,k\,l} = B\!\left(
\begin{array}{cccc}
  0&\star   &j&k   \\
 0 &l   &j&k   
\end{array}
\right)\,\theta_{jk}
+ B\!\left(
\begin{array}{cccc}
  0&\star   &j&l   \\
 0 &k   &j&l   
\end{array}
\right)\,\theta_{jl}
+ B\!\left(
\begin{array}{cccc}
  0&\star   &k&l   \\
 0 &j  &k&l   
\end{array}
\right)\,\theta_{kl}.
\ee
$\zeta_{jk}$ is symmetric with respect to $j,k$
and $\zeta_{jkl}$ is symmetric with respect to $j,k,l.$ 

The following identity holds true by definition:
\beq
\zeta_{jkl} = 0.
\eeq
Indeed, by (81),
\be
&&\zeta_{jkl} = - 2\,\rho_{jk}\Delta_{ljk}\,\theta_{jk} - 2\,\rho_{jl}\Delta_{kjl}\,\theta_{jl} - 2\rho_{kl}\Delta_{jkl}\,\theta_{kl}\\
&&{}\quad \  = 2\Delta_{jkl}\{- d\rho_{jk} + d\rho_{jl} - d\rho_{kl}\} = 0,
\ee
 since $\rho_{jk} - \rho_{jl} + \rho_{kl} = 0$ by definition.

(73) and (74) give

%prop51
\begin{prop}

For different indices $j,k,l,$
we have 
\be
&&\Theta_{j,j} = - (\lambda_\infty + \lambda_j - 1)\,\theta_j 
+ \sum_{\nu\ne j}\,\lambda_\nu\,\theta_{j\nu},\\
&&\Theta_{k,j} = - \lambda_k\,(\theta_j + \theta_{jk}),\\
&&\Theta_{jk,j} = \lambda_k\, \{\zeta_{jk,j}
- \sum_{\nu \ne j,k}\,\frac{\lambda_\nu}{\lambda_\infty - 1}\,
\frac{B\!
\left(
\begin{array}{cccc}
  0&\star   &k&\nu   \\
  0&\star   &k&j   
\end{array}
\right)
}{B(0\star j\,k\,\nu)}\,\zeta_{j\,k\,\nu}
\} = \lambda_k \,\zeta_{jk,j},\\
&&\Theta_{kl,j} = \frac{\lambda_k\,\lambda_l}{\lambda_\infty - 1}\,
\frac{B(0\star k\,l)}{B(0\star j\,k\,l)}\,\zeta_{j\,k\,l} = 0.
\ee
\end{prop}

To obtain explicit formulae for $\Theta_{K,jk},$ the  identities (75) are used.
 For different indices $j,k,l,i,$
\beq
&&\Theta_{j,jk} = T_{-\varepsilon_k}(\Theta_{k,j})\,\frac{\tilde{\gamma}_{j,k}^k}{\lambda_k - 1}
+ T_{-\varepsilon_k}(\Theta_{jk,j})\,\frac{\tilde{\gamma}_{j,jk}^k}{\lambda_k - 1},\\
&&\Theta_{l,jk} = T_{-\varepsilon_k}(\Theta_{k,j})\,\frac{\tilde{\gamma}_{l,k}^k}{\lambda_k - 1}
 + T_{-\varepsilon_k}(\Theta_{jk,j})\,\frac{\tilde{\gamma}_{l,jk}^k}{\lambda_k - 1},
\eeq
which is symmetric with respect to $j,k$,
\beq
&&\Theta_{jk,jk} = T_{-\varepsilon_k}(\Theta_{j,j})\,\frac{\tilde{\gamma}_{jk,j}^k}{\lambda_k - 1}
+ T_{-\varepsilon_k}(\Theta_{k,j})\,\frac{\tilde{\gamma}_{jk,k}^k}{\lambda_k - 1}
+ T_{-\varepsilon_k}(\Theta_{jk,j})\,\frac{\tilde{\gamma}_{jk,jk}^k}{\lambda_k - 1}\nonumber\\
&&{}\qquad \quad + \sum_{\nu\ne j,k} T_{-\varepsilon_k}(\Theta_{j\nu,j})\,\frac{\tilde{\gamma}_{jk,j\nu}^k}{\lambda_k - 1},
\eeq
which is symmetric with respect to $j,k,$
\beq
\Theta_{jl,jk} = T_{-\varepsilon_k}(\Theta_{jk,j})\,\frac{\tilde{\gamma}_{jl,jk}^k}{\lambda_k - 1}
+ T_{-\varepsilon_k}(\Theta_{jl,j})\,\frac{\tilde{\gamma}_{jl,jl}^k}{\lambda_k - 1}.
\eeq
All other $\Theta_{K,J} = 0.$ Hence

%prop52
\begin{prop}
For different indices $j,k,l \in \{1,2,..., m\},$ we have
\be
&&\Theta_{j,jk} = - \lambda_j\,\frac{\zeta_{jk}}{B(0\star j\,k)}
- (\lambda_\infty - 1)\, \frac{\zeta_{jk,j}}{B(0\star j\,k)},\\
&&\Theta_{l,jk} = - \lambda_l\,\frac{\zeta_{jk}}
{B(0\star j\,k)},\\
%\ee
%
%\be
&&\Theta_{jk,jk} = \lambda_k \theta_j + 2\,(\lambda_j + \lambda_k - 1) \,\theta_j\,\frac{B(0\star j)\,B\!
\left(
\begin{array}{ccc}
  0&\star   &k   \\
  0&j   &k      
\end{array}
\right)
}{B(0\star j\,k)} \\
&&{}\quad + \lambda_j\,\theta_k + 2\,(\lambda_j + \lambda_k - 1)\,\theta_k \, \frac{B(0\star k)\,B\!
\left(
\begin{array}{ccc}
  0&\star   &j   \\
  0&k  &j      
\end{array}
\right)
}{B(0\star j\,k)}\\
&&{}\quad -
(\lambda_j + \lambda_k - 1)\,\theta_{jk}\,\{1 + 2\,\frac{B(0\,j\,k)\,B\!
\left(
\begin{array}{ccc}
  0&\star   &j   \\
  0&\star   &k    
\end{array}
\right)
}{B(0\star\,j\,k)}\} + \sum_{\nu\ne j,k}\,\lambda_\nu\,\theta_{j\nu}\\
&&{}\quad  + \sum_{\nu\ne j,k}\lambda_\nu\,\bigl\{
\frac{
B(0\star j)\,B\!
\left(
\begin{array}{cccc}
 0 &\star&k   &\nu   \\
 0 & j  &k&\nu   
\end{array}
\right)
}{B(0\star\,j\,k\,\nu)}\,\theta_j + 
\frac{
B(0\star k)\,B\!
\left(
\begin{array}{cccc}
 0 &\star&j   &\nu   \\
 0 & k  &j&\nu   
\end{array}
\right)
}{B(0\star\,j\,k\,\nu)}\,\theta_k\\
&&{}\quad +
\frac{
B(0\star \nu)\,B\!
\left(
\begin{array}{cccc}
 0 &\star&j   &k   \\
 0 & \nu  &j&k  
\end{array}
\right)
}{B(0\star\,j\,k\,\nu)}\,\theta_\nu\\
\ee
%page change
\be
&&{}\quad  +
\frac{B\!
\left(
\begin{array}{cccc}
  0&\star   &k&\nu   \\
  0&j   &k&\nu   \\
\end{array}
\right)\,
B\!
\left(
\begin{array}{cccc}
  0&\star   &j&\nu   \\
  0&\star   &j&k   \\
\end{array}
\right)\,
B\!
\left(
\begin{array}{ccc}
  0&\star   &k  \\
  0  &j&k  \\
\end{array}
\right)
}{B(0\star j\,k\,\nu)\,B(0\star\,j\,k)}\theta_{k\nu}\\
&&{}\quad  +
\frac{B\!
\left(
\begin{array}{cccc}
  0&\star   &j&\nu   \\
  0&k   &j&\nu  \\
\end{array}
\right)\,
B\!
\left(
\begin{array}{cccc}
  0&\star   &k&\nu  \\
  0&\star   &k&j   \\
\end{array}
\right)\,
B\!
\left(
\begin{array}{ccc}
  0&\star   &j   \\
  0  &k&j \\
\end{array}
\right)
}{B(0\star j\,k\,\nu)\,B(0\star\,j\,k)}\theta_{j\nu}\\
&&{}\quad +  \frac{B\!
\left(
\begin{array}{ccc}
  0&\star   &k   \\
  0&j   &k     
\end{array}
\right)\,
B\!
\left(
\begin{array}{ccc}
  0&\star   &j   \\
  0&k   &j     
\end{array}
\right)
}{B(0\star j\,k)}
\,\theta_{jk}
\bigr\}, \\
%\ee
%\be
&&\Theta_{jl, jk} = \lambda_l\,\{
- \bigl(
\frac{B(0\star j\,l)\,
B\!
\left(
\begin{array}{cccc}
 0 &\star&j   &k   \\
  0&l&j   &k     
\end{array}
\right)\,
B\!
\left(
\begin{array}{cccc}
 0 &\star&k   \\
  0&\star&j        
\end{array}
\right)
}{B(0\star j\,k\,l)\,B(0\star j\,k)}\\
&&{}\quad \bigr.
+
\bigl.
\frac{B\!
\left(
\begin{array}{cccc}
 0 &\star&j&l   \\
  0&k&j&l        
\end{array}
\right)\,
B\!
\left(
\begin{array}{cccc}
 0 &\star&l   \\
  0&\star&j        
\end{array}
\right)}{B(0\star j\,k\,l)}\bigr)
\,\theta_j
+ \frac{B(0\star k)\, B(0\star j\,l)\,B\!
\left(
\begin{array}{cccc}
 0 &\star   &j&k   \\
 0 &l&j   & k  \\
\end{array}
\right)
}{B(0\star j\,k\,l)\, B(0\star j\,k)}\,\theta_k \\
&&{}\quad + 
\frac{B(0\star l)\,B\!
\left(
\begin{array}{cccc}
 0 &\star   &j&l   \\
 0 &k&j   & l  \\
\end{array}
\right)
}{B(0\star j\,k\,l)}\,\theta_l
 + \frac{B\!
\left(
\begin{array}{ccc}
  0&\star   &l   \\
 0 &j  &l   \\
\end{array}
\right)\,
B\!
\left(
\begin{array}{cccc}
 0 &\star   &j&l   \\
 0 &k&j   & l  \\
\end{array}
\right)
}{B(0\star j\,k\,l)}
\,\theta_{jl}\\
&&{}\quad + 
\frac{B\!
\left(
\begin{array}{ccc}
  0&\star   &k   \\
  0&j   &k    
\end{array}
\right)B(0\star j\,l)\,
B\!
\left(
\begin{array}{cccc}
  0&\star&j   &k   \\
  0&l&j   &k    
\end{array}
\right)
}{B(0\star j\,k\,l)\,B(0\star j\,k)} \theta_{jk}\}
\ee
for every triple $j,k,l$.
Hence we have the symmetry with respect to the transposition 
$\sigma_{jk}$ between the indices $j,k${\rm :}
\be
&&\sigma_{jk}(\Theta_{j,jk}) = \Theta_{k,jk},\,
\sigma_{jk}(\Theta_{l,jk}) = \Theta_{l,jk},\\
&& \sigma_{jk}(\Theta_{jk,jk}) = \Theta_{jk,jk},\,
\sigma_{jk}(\Theta_{jl,jk}) = \Theta_{kl,jk}. 
\ee
\end{prop}

%pf of Prop 52
\begin{pf}
Indeed, the first two identities follow from (85) and (86). 
As to the third one,  (87) and Lemma 51 imply
\be
&&\Theta_{jk,jk} = - (\lambda_\infty + \lambda_j - 2)\,\theta_j + (\lambda_k - 1)\,\theta_{jk}
+ \lambda_j\,(\theta_j + \theta_{jk})\,\frac{B\!
\left(
\begin{array}{ccc}
  0&\star   &j   \\
  0&\star   &k   
\end{array}
\right)
}{B(0\star k)}\\
&& + \,\zeta_{jk,j}
\{- \lambda_j \,\frac{B\!
\left(
\begin{array}{ccc}
  0&\star   &j   \\
  0&\star   &k   
\end{array}
\right)\,
B\!
\left(
\begin{array}{ccc}
  0&\star   &k   \\
  0&j   &k   
\end{array}
\right)}{B(0\star j)\,B(0\star k)} +
(\lambda_k + \lambda_\infty - 2)\,\frac{B\!
\left(
\begin{array}{ccc}
  0&\star   &j   \\
  0&k   &j   
\end{array}
\right)}{B(0\star j\,k)}\\
&&{}\qquad - \sum_{\nu\ne j,k}\,\lambda_\nu\,
\frac{B\!
\left(
\begin{array}{cccc}
  0&\star&j   &\nu   \\
  0&\star&j   &k   
\end{array}
\right)\,
B\!
\left(
\begin{array}{cccc}
  0&\star& j  &k   \\
  0&\nu&j   &k   
\end{array}
\right)}{B(0\star\,j\,k\,\nu)\,B(0\star\,j\,k)}\} + \sum_{\nu\ne j,k}\,\lambda_\nu\,\theta_{j\nu}\\
&&{}\qquad +
\sum_{\nu\ne j,k}\,\lambda_\nu \,\zeta_{j\nu,j}\frac{B\!
\left(
\begin{array}{cccc}
  0&\star&j   &k   \\
  0&\nu&j   &k   
\end{array}
\right)}{B(0\star j\,k\,\nu)}.
\ee

This equals the {\rm RHS} of the formula for $\Theta_{jk,jk}$ in Proposition 52
due to Jacobi identity.

Similarly (88) and Lemma 50  imply
\be
&&\Theta_{jl,jk} = \lambda_l\,\{\zeta_{jk,j}\,
 \frac{B(0\star j\,l)\,
B\!
\left(
\begin{array}{cccc}
 0 &\star&j   &k   \\
  0&l&j   &k     
\end{array}
\right)}
{B(0\star j\,k\,l)\,B(0\star j\,k)}
+ \zeta_{jl,j}\,\frac{B\!
\left(
\begin{array}{cccc}
 0 &\star   &j&l   \\
 0 &k&j   & l  \\
\end{array}
\right)
}{B(0\star j\,k\,l)}\}.
\ee

This equals the {\rm RHS} of the formula for $\Theta_{jl,jk}$
in Proposition 52.
\qed
\end{pf}

\begin{re}
As is seen in (83), the fundamental linear relations among admissible forms are written 
independently of $\lambda_j$'s .   This means that the expression of the {\rm RHS}  stated in Propositions 51 and 52
are unique as rational function of $\lambda$.
\end{re}

We can make a conjectural formula for the Wronskian generalizing {\rm(76)}
as follows.

\bigskip
%\vspace{.5cm}

{\bf [Conjecture]}

\medskip

{\it Suppose that $n,\;m$ are general under the condition 
$m\le n+1.$ $T_r(\Theta)$ defined by
\be
T_r{(\Theta)} = \sum_{J \in {\cal B}}\,\Theta_{J,J},
\ee
has the representation 
\be
&&T_r(\Theta) = d_B\log{\bf W},\\
&&{\bf W} = C\,\prod_{p=1}^{n+1}
\prod_{J\in {\cal B}, |J| = p}  \,B(0\star J)^{\lambda_J + \frac{n - p - 1}{2}},
\ee
where $\lambda_J = \sum_{j\in J}\,\lambda_j$
and $C$ does not depend on $\lambda.$}

$C$ has not been known to the authors.

\bigskip

\bigskip

%\vspace{5cm}

%\pagebreak

%appendix

\def\appendix{
\par
\setcounter{section}{0}
\setcounter{subsection}{0}
\def\thesection{\Alph{section}}}
\renewcommand{\theequation}{A.\arabic{equation}}
\setcounter{equation}{0}

\vspace{.5cm}
%\appendix{\hspace{5cm}\bf Appendix}
%\appendix{}
\appendix{\noindent \Large{\bf Appendix} \quad {\LARGE \bf Contiguity relation in negative direction}}
%\section{Contiguity relation in negative direction}
\author{}
\date{}
%\maketitle

\bigskip

\bigskip

%\subsection{Theorem of contiguity in negative direction}
\noindent{\bf \Large A.1 \quad Theorem of contiguity }

%\vspace{.5cm}

\bigskip

\bigskip

For admissible $F_J$ or  $W_0(J)\varpi\,(1\le |J|\le n+1),$
we have by definition
\be
&&F_J^{(j)} = (\lambda_j-1)\,\mathfrak{M}_{f_j^{-1}}\,F_J = 
\left\{
\begin{array}{ccc}
(\lambda_j-1)\, F_{j\, J}\quad(j\notin J), \\
  (\lambda_j-1) f_j^{-2}\,F_{\partial_jJ }\quad(j\in J) .
\end{array}
\right.
\ee
In the same way,
\be
&&W_0^{(j)}(k)\varpi  = 
\left\{
\begin{array}{ccc}
 (\lambda_j -1) \,B(0\star k) \,F_{jk}\quad(j\ne k)\\
 B(0\star j) \, F_j^{(j)}\quad(k = j)
\end{array}
\right.
\ee
and for $|J| \ge 2$
\be
&&W_0^{(j)}(J)\varpi =
\left\{
\begin{array}{ccccc}
(\lambda_j-1) 
\bigl\{ - \sum_{k\in J}
 B\!
 \left(
\begin{array}{ccccc}
 0 &\star& \partial _k J    \\
 0 & k  &  \partial_k J \\
\end{array}
\right)
F_{j\, \partial_k J} +
B(0\star L)\,F_{j J}\bigr\}\\
\quad(j\in J^c),\\ 
-  (\lambda_j-1)\,B\!
 \left(
\begin{array}{ccccc}
 0 &\star& \partial _j J    \\
 0 & j  &  \partial_j J \\
\end{array}
\right)
F_{J}
 - \sum_{k\in \partial_jJ} 
 B\!
 \left(
\begin{array}{ccccc}
 0 &\star& \partial _k J    \\
 0 & k  &  \partial_k J \\
\end{array}
\right)
F_{ \partial_k J}^{(j)}\\
 +
B(0\star \,J)\, F_J^{(j)}
\quad(j\in J) .
\end{array}
\right. 
\ee

Since $\mathfrak{M}_{f_j^{-1}}F_J,\; \mathfrak{M}_{f_J^{-1}}W_0(J)\varpi$ belong to $H_\nabla^n(X,\Omega^\cdot(*S))$
as cohomology class,
we can represent them as linear combinations of admissible $F_K:$
\beq
&&F_J^{(j)} = (\lambda_j -1)\,\mathfrak{M}_{f_j^{-1}}\,F_J  \sim \sum_{K : {\rm admissible}} \tilde{\gamma}_{K,J}^j \,F_K,\\
&&W_0^{(j)}(J)\varpi =(\lambda_j -1)\,\mathfrak{M}_{f_j^{-1}}\,W_0(J)\varpi\sim \sum_{K : {\rm admissible}} {\gamma}_{K,J}^j \,F_K,
\eeq
where
the two matrices $\gamma^j = (\gamma_{K,J}^j),\,\tilde{\gamma}^j= (\tilde{\gamma}_{K,J}^j) $
are related with each other through $\tilde{\beta}$ due to Lemma 10:
\beq
\tilde{\gamma}_{K,J}^j = \sum_{L \subset J ; |L| \ge 1} \gamma_{K,L}^j \tilde{\beta}_{L,J}.
\eeq

Concerning 
$W_0^{(j)}(J)\varpi\ (|J| = p,\,1\le p\le n+1 ),$ the following 
recurrence (so called contiguity) relations hold.

%\vspace{.5cm}
\bigskip

%thmA
[{\bf Theorem A}]

\medskip

{\it
%\hspace{1cm}

{\rm(i)}\;  Case where $j\notin J$,
\beq
&&W_0^{(j)}(k)\,\varpi = (\lambda_j - 1)\,B(0\star k)\,F_{jk}\quad(J = \{k\}, \,j\ne k),\\
&&W_0^{(j)}(J)\varpi 
 =  (\lambda_j -1)\,\bigl\{
- \sum_{k\in J} B\!
\left(
\begin{array}{ccc}
 0  &\star  &\partial_k J   \\
  0& k  &\partial_k J   \\
\end{array}
\right)\, F_{j\,\partial_k J}
+ B(0\star J) F_{j \,J}\bigr\}\,\,(p\ge 2),\nonumber\\
\eeq

{\rm (ii)}\;
Case where $j\in J$,
\beq
&&W_0^{(j)}{(j)}\varpi \sim -(\lambda_\infty + n - 2 + \lambda_j)\,F_j \nonumber\\
&& {}\qquad \qquad \ - 
\sum_{k=1,k\ne j}^{m} \lambda_k \, \{F_k + B\!
\left(
\begin{array}{ccc}
 0 &\star   &k   \\
 0 &\star   &j   \\
\end{array}
\right)\,F_{jk}
\}\quad(p=1),\\
&&W_0^{(j)}(J)\varpi \sim (\lambda_\infty + n - p) Z_{J}^j + \sum_{k\in J^c} \lambda_k Z_{J,k}^j\quad(p\ge 2),
\eeq
where
\beq
&& Z_J^j = B(0\,\partial_j J)\, F_{\partial_j J} - \sum_{\nu \in \partial_j J} 
B\!
\left(
\begin{array}{ccc}
 0 &j   &\partial_j\partial_\nu J   \\
 0 &\nu   &\partial_j\partial_\nu J   \\
\end{array}
\right) F_{\partial_\nu J } +
B\!
\left(
\begin{array}{ccccc}
 0 &\star   & \partial_j J   \\
  0&j   &\partial_j J   
\end{array}
\right)\, F_{J},\nonumber\\
\\
&&
Z_{J,k}^j = B\!
\left(
\begin{array}{cccc}
0  &\star   &\partial_j J   \\
 0 & k  & \partial_j J  \\  
\end{array}
\right)\,F_{k\, \partial_j J}
- \sum_{\nu\in \partial_j J}B\!
\left(
\begin{array}{ccccc}
  0&\star   &j &\partial_j\partial_\nu J   \\
 0 &k   &\nu & \partial_j\partial_\nu J   \\
\end{array}
\right)\,F_{k\,\partial_\nu J}\nonumber\\
&&{}\qquad - B\!
\left(
\begin{array}{ccccc}
 0 &\star   &k &\partial_j J   \\
 0 &\star   &j&\partial_j J   \\
\end{array}
\right)\,F_{k\,J} \quad(k \in J^c).
\eeq
}
%end theorem

%re
\begin{re}
When $ p = n+1,$
the terms $F_{j\,J}$ in {\rm(A.5)} and  $F_{k\,J}$ in {\rm(A.9)}
can be described by linear combinations of
admissible $F_K\, (|K| = n+1)$
due to {\rm(48)}. As a result, $W_0^{(j)}\,\varpi$ and $Z_{J,k}^j$
have the following reduced forms:
\beq
&&W_0^{(j)}(J)\,\varpi = (\lambda_j - 1)\, 
\frac{B\!
\left(
\begin{array}{ccc}
  0&\star   &J   \\
  0&j   &J   
\end{array}
\right)}{B(0\star j\,J)}\,
\bigl\{B(0\star J)\,F_J \nonumber\\
&&{}\quad \quad \quad \quad \quad -
\sum_{\nu\in J}
B\!
\left(
\begin{array}{ccccc}
  0&\star&j   &\partial_\nu J   \\
  0&\star   &\nu&\partial_\nu J     
\end{array}
\right)\,F_{j\,\partial_\nu J}
\bigr\}, \\
%\eeq
%\beq
&&Z_{J,k}^j = \frac{B\!
\left(
\begin{array}{ccc}
  0&\star   &J   \\
  0& k  & J  \\
\end{array}
\right)
}
{B(0\star k\, J)}
\bigl\{
B(0\star k\,\partial_jJ)\,F_{k\,\partial_jJ}
- B\!
\left(
\begin{array}{ccccc}
 0 &\star   &k& \partial_jJ   \\
  0&\star&j&   \partial_jJ  
\end{array}
\right)\,F_{J}\nonumber\\
&&{}\qquad - \sum_{\nu\in \partial_jJ}
B\!
\left(
\begin{array}{cccccc}
  0&\star   &k&j&\partial_j\partial_\nu J   \\
  0&\star   &k&\nu&\partial_j\partial_\nu J   \\
\end{array}
\right)\,F_{k\,\partial_\nu J}\bigr\}.
\eeq

It seems remarkable that
$F_K$ appearing in the {\rm RHS} of {\rm (A.5) and (A.6)} are adjacent to $F_J$,
i.e., satisfy the condition
\be
|J - J\cap K | \le 1, |K - J\cap K| \le 1.
\ee

As an immediate consequence from Theorem A, we see that $\gamma_{K,J}^j,\,\tilde{\gamma}_{K,J}^j$ 
in (A.1) and (A.2)
are all inhomogeneous and linear in $\lambda$. 
\end{re}

\bigskip

\bigskip

%secA.2
\noindent{\bf \Large A.2 \quad Proof of Theorem A}

%\vspace{.5cm}
\bigskip

\bigskip

To prove Theorem A, we use contiguity relations in positive direction
stated in  Theorems 1 and 3 , Propositions 21 and 23.

\bigskip

%\vspace{.5cm}
%propA.1
{\bf Proposition A.1}
{\it
Fix $J \subset N, \,|J| = p,\,J^c = N_m - J$ and 
$j\in J$.
$F_J^{(j)}$ can be represented by linear combinations of admissible $F_K$
in recurrence form:
\beq
&&B(0\star j) \,F_j^{(j)} = W_0^{(j)}(j)\varpi 
\sim - (\lambda_\infty + n - 2 + \lambda_j) F_j \nonumber\\
&&{}\qquad \qquad \quad \quad  - 
\sum_{k=1, k\ne j}^{m} \lambda_k\{F_k + B\!
\left(
\begin{array}{ccc}
 0 &\star   &k   \\
  0&\star   &j   \\
\end{array}
\right)\,F_{jk}\}\quad(p=1) ,\\
 &&B(0\star J)\,F_{J}^{(j)} \sim U_0 + (\lambda_\infty + n- p-1)\,U_\infty 
+ \sum_{k=1}^{m} \lambda_k U_k\quad(2\le p\le n+1), \nonumber\\
\eeq
where
\be
&&U_0 = B(0\star \partial_{j}J) \sum_{\nu\in \partial_{j}\!J}\,F_{\partial_{j} J }^{(\nu)}
- \sum_{\mu \in \partial_{j}J}
B\!
\left(
\begin{array}{ccccc}
 0 &\star   &\mu& \partial_j\partial_\mu J  \\
  0& \star  &j& \partial_j\partial_\mu J  \\
\end{array}
\right)
\sum_{\nu \in \partial_\mu J}
\,F_{\partial_\mu J}^{(\nu)},\\
&&
U_\infty = B\!
\left(
\begin{array}{ccccc}
 0 &j   &\partial_{j} J  \\
  0& \star  &\partial_{j}J  \\
\end{array}
\right)\,F_{J},\\
&&U_{\mu} = - B\!
\left(
\begin{array}{ccccc}
 0 &\star   &\mu& \partial_j\partial_\mu J  \\
  0& \star  &j& \partial_j\partial_\mu J  \\
\end{array}
\right)
F_{J}\quad(\mu\in \partial_jJ), \\
&& U_{j} = B(0\star \partial_{j}J) F_{J},\\
&&
U_k = B(0\star \partial_{j}J)F_{k\,\partial_{j}J}
- \sum_{\mu \in \partial_{j}J} 
B\!
\left(
\begin{array}{ccccc}
 0 &\star   &\mu& \partial_j\partial_\mu J  \\
  0& \star  &j& \partial_j\partial_\mu J  \\
\end{array}
\right)\,F_{k\, \partial_\mu J}\\
&&{}\qquad - B\!
\left(
\begin{array}{ccccc}
 0 &\star   &k& \partial_{j}J  \\
  0& \star  &j& \partial_{j}J  \\
\end{array}
\right)\,F_{k J}\quad(k\in J^c).
\ee
}

Because of symmetry, it is sufficient to prove  Proposition A.1 in the special
case where $J = \{n-p+2,\ldots,n+1\}$ and $j= n+1.$
Suppose that $f_1,f_2,\ldots, f_{n+1}$ have a normalized form (7).
Then we have
\be
df_j\wedge *dQ = \{f_j + f_{n+1} + B\!
\left(
\begin{array}{ccc}
 0 &\star   &j   \\
 0 &\star   & n+1  \\
\end{array}
\right)
\}\,\varpi \quad(1\le j\le m).
\ee
 The following Lemma  is an immediate application of this formula:
 
\bigskip

%\vspace{.5cm}
%lmA.2
{\bf Lemma A.2}
{\it
The Stokes identity implies
\beq
&&0 \sim \nabla \frac{*dQ}{f_J}
= \sum_{\nu\in J}\,
B\!
\left(
\begin{array}{ccc}
 0 &\star   &\nu   \\
  0&\star   &n+1   \\
\end{array}
\right)\, F_{J}^{(\nu)} + R_{n+1},
\eeq
where
\beq
&&
R_{n+1} =  \sum_{\nu\in \partial_{n+1}J}\, F_{\partial_{n+1}J}^{(\nu)}
+ ( n - p -1 + \sum_{\nu\in \partial_{n+1}J} \lambda_\nu + 2\lambda_{n+1})\,F_{J}\nonumber\\
&&{}\qquad \quad + \sum_{k \in J^c}\,\lambda_k
\{B\!
\left(
\begin{array}{ccc}
 0 &\star   &k   \\
  0& \star  &n+1  
\end{array}
\right)\, F_{k J}
+\, F_{k\,\partial_{n+1}J} + F_{J}\}.
\eeq
}

Likewise by the transposition  $\sigma_{\mu\,n+1},$
we have similar identitiy to (A.14):
\beq
0 \sim \sum_{\nu\in J}\,
B\!
\left(
\begin{array}{ccc}
 0 &\star   &\nu   \\
  0&\star   &\mu   \\
\end{array}
\right)\, F_{J}^{(\nu)} + R_{\mu}\quad(\mu\in J),
\eeq
where
\be
R_\mu = \sigma_{\mu, n+1}(R_{n+1}).
\ee

\bigskip

%\vspace{.5cm}
%proof of PropA.1
[{\bf Proof of Proposition A.1}]

Consider $p\times p$ matrix $(B\!
\left(
\begin{array}{ccc}
 0 &\star   &\nu   \\
 0 & \star  &\mu   \\
\end{array}
\right))_{\mu,\nu \in J}
$. Its determinant and the cofactor of $\mu,\nu$th component 
can be evaluated through
Sylvester's determinant identity as follows:
$$
(-1)^{p-1}B(0\star J),\, (-1)^{\mu+\nu + p - 2}B\!
\left(
\begin{array}{ccc}
0  &\star   & \partial_\mu J  \\
 0 & \star  & \partial_\nu J \\
\end{array}
\right)
$$
respectively.
Hence  the system of linear equations (A,16) can be solved with respect to
$F_j^{(\nu)}$ as follows
\be
F_{J}^{(\nu)} \sim  \frac{1}{B(0\star J)}\, \bigl\{
R_\nu \,B(0\star \partial_\nu J) - \sum_{\mu\in J,\, \mu\ne \nu} R_\mu \, B\!
\left(
\begin{array}{cccccc}
 0 &\star   & \nu & \partial_\mu\partial_\nu J \\
 0 &\star   & \mu& \partial_\mu\partial_\nu J  \\
\end{array}
\right)\bigr\},
\ee
which is nothing else than (A.12) and (A.13). 
\qed

\bigskip

We can see by induction that the $F_K$ appearing in the RHS of (A.12) and (A.13) 
satisfy the condition
$ |K - J\cap K| \le 1$.

 Proposition A.1  enables us to evaluate  $F_J^{(j)}$
by recurrence with respect to $p.$

Comparing the part of weight $p+1$ in both sides of (A.13),
we have 
\beq
B(0\star J)\,\tilde{\gamma}_{k J, J}^j = - \lambda_k\,
B\!
\left(
\begin{array}{ccccc}
 0 & \star  &k &\partial_j J   \\
 0 & \star  &j& \partial_j J    \\
\end{array}
\right)\quad(j\in J).
\eeq
On the other hand, obviously from the definition
\beq
\tilde{\gamma}_{j\,J, J}^j = \lambda_j - 1\quad(j\in J^c)
\eeq
and $\tilde{\gamma}_{K, J}^j = 0$
for any other $ K$  such that $|K| = p+1$.

Furthermore, as to the part of weight $p,$
the similar identities hold true for  $F_{\partial_j J}^{(\nu)}\,(\nu \in \partial_j J)$
and $F_{\partial_\mu J}^{(\nu)}\,(\nu \in \partial_\mu J)$. 
Hence the equality of (A.13) restricted to the part of weight $p$ implies

\bigskip

%coA.3
%\vspace{.5cm}
{\bf Corollary A.3}

{\it
Suppose $j \in J$. Then

\beq
&&B(0\star J)\,\tilde{\gamma}_{J, J}^j = (\lambda_\infty + \lambda_j + n - p -1) B\!
\left(
\begin{array}{ccc}
 0 &j   &\partial_j J   \\
 0 &\star& \partial_j J     \\
\end{array}
\right)\nonumber\\
&&{}\qquad \qquad  -
\sum_{\mu\in \partial_j J} \lambda_\mu \,
\frac{B\!
\left(
\begin{array}{ccccc}
 0 &\star &\partial_\mu J   \\
 0 &\mu   & \partial_\mu J   \\
\end{array}
\right)\,
B\!
\left(
\begin{array}{ccccc}
 0 &\star   &\mu &\partial_j\partial_\mu J   \\
  0&\star   &j&\partial_j\partial_\mu J   \\
\end{array}
\right)
}
{B(0\star \partial_\mu J)}
,\nonumber\\
\\
&&B(0\star \partial_\mu J)\,\tilde{\gamma}_{J, \partial_\mu J}^j 
= - \lambda_\mu\,
B\!
\left(
\begin{array}{ccccc}
  0&\star   &\mu& \partial_j\partial_\mu J   \\
  0&\star   &j & \partial_j\partial_\mu J    \\
\end{array}
\right)\quad(\mu \in \partial_j J),\\
&&B(0\star J)\,\tilde{\gamma}_{k\,\partial_j J, J}^j
= \lambda_k\, B\!
\left(
\begin{array}{ccc}
 0 &\star   &\partial_j J   \\
  0&k   &\partial_j J   \\
\end{array}
\right)\quad(k \in J^c),\\
%(6)
&&B(0\star J)\,\tilde{\gamma}_{k\,\partial_\mu J, J}^j = 
- \lambda_k\, 
\frac{B\!
\left(
\begin{array}{ccccc}
 0 &\star  &\partial_\mu J   \\
 0 &k   & \partial_\mu J   \\
\end{array}
\right)\,
B\!
\left(
\begin{array}{ccccc}
 0 &\star   &\mu &\partial_j\partial_\mu J   \\
  0&\star   &j&\partial_j\partial_\mu J   \\
\end{array}
\right)
}
{B(0\star \partial_\mu J)}\nonumber\\
&&{}\qquad \qquad \qquad \qquad \qquad \qquad \qquad \qquad 
(k\in J^c,\,\mu\in \partial_j J).
\eeq
}

Remark that $\tilde{\gamma}_{K,J}^j = 0$ for other $K$ such that $|K| = p$.

To prove the identities (A.19) and (A.22), we have only to notify the equalities
\be
&&
B
\left(
\begin{array}{ccc}
  0&\star&K     \\
  0&j   &K  \\ 
\end{array}
\right)
= B(0\star K) - \sum_{\nu\in K}\,B
\left(
\begin{array}{ccccc}
0  &\star   &\nu&\partial_\nu K   \\
 0 &\star   &j&\partial_\nu K 
 \end{array} 
\right)\quad(j \in K^c),\\
&&B\!
\left(
\begin{array}{ccccc}
0& \star & j & K  \\
0 &\nu& j& K\\
\end{array}
\right)
B\!
\left(
\begin{array}{ccccc}
0& \star & k & K \\
0 &\star& j & K\\
\end{array}
\right)
-
B\!
\left(
\begin{array}{ccccc}
0& \star& j  & K \\
0& k& j & K\\
\end{array}
\right)
B\!
\left(
\begin{array}{ccccc}
0& \star & \nu & K \\
0& \star& j & K\\
\end{array}
\right)\\
&&{}\qquad \qquad \qquad = B(0\star j\,K)) \, 
B\!\left(
\begin{array}{ccccc}
 0 &k   &\nu &K   \\
 0 & j  &\star&K   
\end{array}
\right)\quad(j,k,\nu\in K^c).
\ee

In the special cases where $p= 1, 2\;(n\,{\rm general})$, i.e., for $J = \{j\}, \,\{j,k\},$
more explicit formulae are obtained from (A.12) and (A.13):
\be
&&F_k^{(j)} =(\lambda_j - 1)\, F_{jk}\quad(j\ne k),\\
&&B(0\star j)\,F_j^{(j)} \sim  - (\lambda_\infty + \lambda_j+ n-2)\,F_j
- \sum_{k=1, k\ne j}^{m} \,\lambda_k\{F_k + 
B\!\left(
\begin{array}{ccc}
 0 &\star   &k   \\
 0 &\star   &j   \\ 
\end{array}
\right)\,F_{kj},
\ee
namely,
\be
&&\tilde{\gamma}_{j,k}^j = \lambda_j -1,\; B(0\star j)\,\tilde{\gamma}_{j,j}^j = - (\lambda_\infty + \lambda_j + n - 2),\\
&&B(0\star j) \,\tilde{\gamma}_{k,j}^j = - \lambda_k\;(k\ne j),\;
B(0\star j)\,\tilde{\gamma}_{jk,j}^j = - \lambda_k \,B\!
\left(
\begin{array}{ccc}
  0& \star  &k   \\
 0 &\star   &j   
\end{array}
\right).
\ee

Likewise
\be
&&F_{kl}^{(j)} = (\lambda_j -1)  \,F_{jkl}\quad(j\ne k,l),\\
&&F_{jk}^{(j)} \sim \tilde{\gamma}_{j,jk}^j\,F_j + \tilde{\gamma}_{k,jk}^jF_{k}
+ \tilde{\gamma}_{jk,jk}^j\,F_{jk} \nonumber\\
&&{}\quad \quad + \sum_{l=1;\, l\ne j,k}^{m}\{\tilde{\gamma}_{l,jk}^j \,F_l + \tilde{\gamma}_{lj,jk}^j\,F_{lj}
+\tilde{\gamma}_{lk,jk}^j\,F_{lk}
+ \tilde{\gamma}_{ljk,jk}^j\,F_{ljk}\},
\ee
where
\be
&& \tilde{\gamma}_{jkl,kl}^{(j)} = \lambda_j -1\quad(j\ne k,l),\\
&&B(0\star jk)\,\tilde{\gamma}_{j,jk}^j = - \lambda_j\,
\frac{B\!
\left(
\begin{array}{ccc}
 0 &\star   &j   \\
  0& k  &j
\end{array}
\right)
}{B(0\star j)} +
(\lambda_\infty + n -2) 
\frac{B\!
\left(
\begin{array}{ccc}
 0 &\star   &j   \\
  0& \star  &k
\end{array}
\right)
}{B(0\star j)} ,\\
&&
B(0\star jk)\,\tilde{\gamma}_{k,jk}^j = - \lambda_k
\frac{B\!
\left(
\begin{array}{ccc}
 0 &\star   &j   \\
  0& k  &j
\end{array}
\right)
}{B(0\star j)} 
- (\lambda_\infty + n-2),\\
&&
B(0\star jk)\,\tilde{\gamma}_{jk,jk}^j = - \lambda_k 
\frac{B\!
\left(
\begin{array}{ccc}
 0 &\star   &j   \\
  0& \star  &k
\end{array}
\right)\,
B\!
\left(
\begin{array}{ccc}
 0 &\star   &j  \\
  0& k  &j
\end{array}
\right)}{B(0\star j)}\\
&&{}\qquad \qquad \qquad \quad +\, (\lambda_\infty  + \lambda_j + n-3)B\!
\left(
\begin{array}{ccc}
 0 &\star   &k   \\
  0& j  &k
\end{array}
\right), \\
%\ee
%
%\be
&&B(0\star jk)\,\tilde{\gamma}_{l,jk}^j = - \lambda_l \frac{B\!
\left(
\begin{array}{ccc}
 0 &\star   &j   \\
  0& k  &j
\end{array}
\right)}{B(0\star j)}\quad(l \ne j,k),\\
&&B(0\star  jk)\,\tilde{\gamma}_{lj,jk}^j = - \lambda_l
\frac{B\!
\left(
\begin{array}{ccc}
 0 &\star   &j   \\
  0& \star  &k
\end{array}
\right)\,
B\!
\left(
\begin{array}{ccc}
 0 &\star   &j  \\
  0& l  &j
\end{array}
\right)}{B(0\star j)}\quad(l\ne j,k), \\
%\ee
%
%\be
&&B(0\star  jk)\,\tilde{\gamma}_{lk,jk}^j = \lambda_l 
B\!
\left(
\begin{array}{ccc}
 0 &\star   &k  \\
  0& l  &k
\end{array}
\right)\quad(l \ne j,k),\\
&&
B(0\star  jk)\,\tilde{\gamma}_{ljk,jk}^j = - \lambda_l\,
B\!
\left(
\begin{array}{ccccc}
 0 &\star&l   &k  \\
  0&\star& j  &k
\end{array}
\right)\quad(l \ne j,k).
\ee

The above two equations are equivalent to (A.12) and the following relation
relative to $W_0^{(j)}(J)\,\varpi:$
\beq
&&W_0^{(j)}(jk)\,\varpi \sim (\lambda_\infty + n - 2)\,\{F_j - F_k + 
B\!
\left(
\begin{array}{ccc}
 0 &\star   &k   \\
  0&j   & k  \\
\end{array}
\right)
\,F_{jk}\}\nonumber\\
&& +
\sum_{l=1; l\ne j,k}^{m}\, \lambda_l\,\{
B\!
\left(
\begin{array}{ccc}
 0& \star&k   \\
 0 & l  &k   \\  
\end{array}
\right)
\,F_{lk}
- B\!
\left(
\begin{array}{ccc}
 0 &\star&j   \\
 0 & l  &k   \\  
\end{array}
\right)\,F_{lj}
-  B\!
\left(
\begin{array}{ccccc}
 0 &\star&k&l   \\
 0 & \star  &k&j   \\  
\end{array}
\right)\,F_{ljk}\}.\nonumber\\
\eeq

The following remark is an important consequence from Theorem A. 

%re
\begin{re}
$F_J^{(j)}$ being a linear combination of $F_K$, the coefficients $\tilde{
\gamma}_{K,J}^{j}$ are all {\rm(}inhomogeneous{\rm)} linear functions of $\lambda$ .
Indeed  if $j\notin J,$ then $F_J^{(j)} = (\lambda_j-1) F_{j J}$
is obviously linear in
$\lambda$. 
On the other hand, suppose $j\in J$. 
 If $p=1,$ then from {\rm(A.12)} $F_J^{(j)}$ 
is linear in $\lambda$ with respect to $F_K$.
If $p\ge 2$, then $U_\infty, \,U_k$ in the {\rm RHS} of {\rm(A.13)} does not depend on $\lambda$.
$U_0$ is a linear combination of $F_{\partial_j J}, \,F_{\partial_\mu J}$.
From the hypothesis of recurrence for $|\partial_j J|, |\partial_\mu J| \le p-1$,
it also depends linearly on $\lambda$. Hence
$F_J^{(j)}$ is linear in $\lambda$ with respect to the basis $F_K$.
\end{re}

\medskip

  By the above Remark, $W_0^{(j)}(J)\varpi$ also is a linear combination of $F_K$ whose
coefficients $\gamma_{K,J}^j$ are linear functions of  $\lambda$.

Since the identity (A.4) and (A.5)  are almost obvious
(partial fraction decomposition), we are going to prove (A.6) and (A.7).
As we have already seen, $W_0^{(j)}(J)\varpi$ is expressed by a linear combination
of $F_K$ whose coefficients are linear in $\lambda$
 and the weight of $F_K$  is at most  $p+1$.  The highest weight $F_K$
occurs in case of  $|K| = p+1$ such that $K = \{k\}\cup J \,(k\in J^c)$.
In this case, we have from (A.3) and (A.17)
\beq
&&{\gamma}_{k J,J}^j =  B(0\star J) \tilde{\gamma}_{k J,J}^j\nonumber\\
&& {}\qquad = -\lambda_k\, B\!
\left(
\begin{array}{ccccc}
 0 &\star   &k& \partial_j J   \\
 0 & \star  & j&\partial_j J  \\
\end{array}
\right).
\eeq

Next, consider the case $|K| = p $.
$F_K$ appears when and only when $K = J$ or $K = \{k, \partial_\mu J\}\,(k\in J^c, \mu\in J)$.

\bigskip

%lmA.4p.57
%\vspace{.5cm}

{\bf Lemma A.4}
{\it 
When $|K| = p,$
$\gamma_{K,J}^j\,(j \in J)$
have respectively the following expressions.

For $p =1,$
\beq
\gamma_{j,j}^j = - (\lambda_\infty + n - 2 + \lambda_j), \, \gamma_{k,j}^j = - \lambda_k\,(k\ne j),
\eeq
for $p \ge 2,$
\beq
&&\gamma_{J,J}^j = (\lambda_\infty + n - p)\, B\!
\left(
\begin{array}{ccccc}
 0 &\star   &\partial_j   J\\
  0&j   &\partial_j J   \\   
\end{array}
\right) ,\\
&&\gamma_{k\,\partial_j J, J}^j = \lambda_k\,
B\!
\left(
\begin{array}{ccc}
  0&\star   &\partial_j J  \\
  0&k   &\partial_j J  \\
\end{array}
\right)\quad(k\in J^c),\\
&&\gamma_{k\,\partial_\mu J, J}^j =
- \lambda_k\,B\!
\left(
\begin{array}{cccc}
 0 & \star  &j &\partial_j\partial_\mu J   \\
  0& k &\mu &\partial_j\partial_\mu J   \\  
\end{array}
\right)\quad(k\in J^c,\, \mu\in \partial_j J).
\eeq
}

%proof
\begin{pf}
(A.25) is a direct consequence of (A.12).
As regards to (A.26), by definition of $W_0^{(j)}(J)\varpi$
(A.1) and (A.2) means 
\be
&&\gamma_{J, J}^j 
 =  B(0\star J)\, \tilde{\gamma}_{J,J}^j - \sum_{\nu\in \partial_jJ}
B\!
\left(
\begin{array}{ccc}
0  &\star   &\partial_\nu J   \\
 0 &\nu   &\partial_\nu J  \\
\end{array}
\right)
\, \tilde{\gamma}_{J, \,\partial_\nu J}^j
- (\lambda_j -1) B\!
\left(
\begin{array}{ccc}
 0 &\star   & \partial_j J  \\
 0 &j   &\partial_j J   \\ 
\end{array}
\right)
.
\ee
(A.19) and (A.20),  and the cofactor  and expansion
\be
B(0\star \partial_\nu J) - 
B\!
\left(
\begin{array}{ccc}
  0&\nu   &\partial_\nu   J\\
 0 &\star   &\partial_\nu J   \\
\end{array}
\right)
- \sum_{\mu\in \partial_\nu J} B\!
\left(
\begin{array}{ccccc}
 0 &\star   &\nu &\partial_\mu\partial_\nu   J\\
  0&\star   &\mu &\partial_\mu\partial_\nu J  \\
\end{array}
\right) = 0
\ee\!
derives  (A.26).

(A.27) is obtained from the definition together with (A.21) and 
the identity
\be
\gamma_{k\,\partial_j J, J}^j=
B(0\star J)\, \tilde{\gamma}_{k\,\partial_j J,J}^j .
\ee

Likewise from the definition it holds
\be
\gamma_{k\,\partial_\mu J,J}^j = B(0\star J)\, \tilde{\gamma}_{k\,\partial_\mu J, J}^j
- B\!
\left(
\begin{array}{ccc}
 0 &\star   &\partial_\mu J   \\
 0 &\mu   &\partial_\mu  J \\
\end{array}
\right)\,\tilde{\gamma}_{k\,\partial_\mu J, \,\partial_\mu J}^j.
\ee

(A.28) is a consequence of (A.20) and  (A.22)
and the following Pl\"ucker 's identity
\be
&& B\!
\left(
\begin{array}{cccccc}
0& \star &j   & \partial_j\partial_\mu J  \\
 0 &\star   &\mu&\partial_j\partial_\mu J   \\
\end{array}
\right)\,
B\!
\left(
\begin{array}{cccccc}
0& \star &j   & \partial_j\partial_\mu J  \\
 0 &j   &k&\partial_j\partial_\mu J   \\
\end{array}
\right)\\
&&
= B(0\star \partial_\mu J)\, B\!
\left(
\begin{array}{ccccc}
 0 &\star   &j  &\partial_j\partial_\mu J \\
  0&\mu   &k  & \partial_j\partial_\mu J \\ 
\end{array}
\right) +
 B\!
\left(
\begin{array}{cccccc}
0& \star &j   & \partial_j\partial_\mu J  \\
 0 &\star   &k&\partial_j\partial_\mu J   \\
\end{array}
\right)\,
B\!
\left(
\begin{array}{cccccc}
0& \star &j   & \partial_j\partial_\mu J  \\
 0 &j &\mu&\partial_j\partial_\mu J   \\
\end{array}
\right).\\
\ee
\qed
\end{pf}

As a result of Theorems 1 and 3, we can prove
that $W_0^{(j)}(J)\varpi\,(j \in J)$
has the expression.

%lmA.5
%\vspace{.5cm}
\bigskip

{\bf Lemma A.5}
{\it
When $j\in J$ and $|J|= p \ge 2,$
\beq
W_0^{(j)}(J)\varpi  \sim (\lambda_\infty + n - p) Z_J^j + \sum_{k\in J^c} \lambda_k\,Z_{J,k}^j,
\eeq
where $Z_J, \,Z_{J, k}$ does not depend on $\lambda$.
$Z_J^j$ depends only on $f_\nu\,(\nu\in J)$, 
while $Z_{J,k}^j$ depends 
only on $f_\nu\,(\nu\in J)$ and $f_k$.
}

%pf
\begin{pf}
According to Lemma 6, when the operator $\mathfrak{M}_{f_j^{-1}}\,T_{-\varepsilon_j}$ is made on
both sides of (68), (68) transfers  into the following
\beq
&&(2\lambda_\infty + n -2)\,F_j \sim - W_0^{(j)}(j)\varpi - \sum_{k=1, k\ne j}^m \lambda_k B(0\star k)\,F_{jk}\nonumber\\
&&+ \sum_{q=2}^{n} \,
\frac{(-1)^q}{\prod_{\nu=1}^{q-1} (\lambda_\infty + n - \nu - 1)}\,
\{
\sum_{K\subset N_m, |K|=q, j \in K} \prod_{\nu\in \partial_j K}\,\lambda_\nu\,\,
W_0^{(j)}(K)\varpi \,\nonumber\\
&&{}\qquad + \,
\sum_{K\subset N_m, |K|=q, j\in K^c }\prod_{\nu\in K}\,\lambda_\nu\,\,
\mathfrak{M}_{f_j^{-1}}\,W_0(K)\varpi 
\}\nonumber\\
&&
{}\qquad + \frac{(-1)^{n+1}}{\prod_{\nu=1}^{n} (\lambda_\infty  + n - \nu - 1)}\,
\prod_{\nu\in \partial_j N}\,\lambda_\nu\,\,
W_0^{(j)}(N)\varpi .
\eeq

In view of Remark at the end of \S A.1,
this identity can be regarded as an identity among rational functions of $\lambda.$

From now on, for simplicity we may assume 
$j=1,\, J =\{12 \ldots p\}$
without losing generality.

If we restrict $\lambda$ to the postulate
\be
\lambda_k = 0\quad(p+1\le k\le m),
\ee
then (A.30) becomes
\beq
&&(2\overline{\lambda}_\infty + n -2)\,F_1 \sim - W_0^{(1)}(1)\varpi - \sum_{2\le  k\le p} \lambda_k B(0\star k)\,F_{jk}\nonumber\\
&&+ \sum_{q=2}^{p-1} \,
\frac{(-1)^q}{\prod_{\nu=1}^{q-1} (\overline{\lambda}_\infty + n - \nu - 1)}\,
\{
\sum_{ K\subset J, 1\in K, |K|=q} \prod_{\nu\in K, 2\le \nu}\,\lambda_\nu\,\,
W_0^{(1)}(K)\varpi \,\nonumber\\
&&{}\qquad + \,
\sum_{K\subset J, |K|=q, 1\in K^c }\prod_{\nu\in K}\,\lambda_\nu\,\,
\mathfrak{M}_{f_j^{-1}}\,W_0(K)\varpi 
\}\nonumber\\
&&
{}\qquad + \frac{(-1)^{p}}{\prod_{\nu=1}^{p-1} (\overline{\lambda}_\infty  + n - \nu - 1)}\,
\prod_{2\le \nu\le p}\,\lambda_\nu\,\,
W_0^{(1)}(J)\varpi .
\eeq
where $\overline{\lambda}_\infty = \sum_{j=1}^p \,\lambda_j.$

This equality shows, if  $\lambda_k = 0\,(k\in J^c),$ then
$W_0^{(1)}(J)\varpi$ can be a multiple of $\overline{\lambda}_\infty + n - p,$
because the factor $\overline{\lambda}_\infty + n - p$
appears only in the denominator of the last term of the RHS..
In view of $\overline{\lambda}_\infty = \lambda_\infty - \sum_{k\in J^c}\,\lambda_k$
for general $\lambda,$
we can conclude that (A.29) holds. \qed
\end{pf}

By the repeated use of (A.13), we can see by induction with respect to $p$
that 
$Z_J^j$ depends only on $J$, while $Z_{J,k}^j$ depends only on $J$ and $k.$
After the substitution of (A.29)
 into (A.30), we can compare, as rational functions of $\lambda$,
the coefficients of the term 
$(-1)^p\,\frac{\prod_{\nu\in \partial_j J}\,\lambda_\nu}{\prod_{\nu = 1}^{p-2}\,(\lambda_\infty + n - 1- \nu)}$
in (A.30) (see Lemma 43)
\beq
&&Z_J^j = \mathfrak{M}_{f_j^{-1}}\,W_0(\partial_j J)\varpi + \sum_{\nu\in \partial_j J}\,Z_{\partial_\nu J, \nu}^j \quad(p\ge 3),\\
&&
Z_{jk}^j = \hat{Z}_{jk}^j = 
 F_j- F_k -
  B\!
\left(
\begin{array}{ccc}
  0&\star   &k   \\
 0 &j   & k  \\
\end{array}
\right)\, F_{jk}\quad(p = 2).
\eeq

In fact, suppose first that $J = \{j, k\}\,(k \ne j)$.  By taking $\lambda_l = 0$ for
all $l$ except for $j, k$ in (A.30), we have the identities:
\be
&&(2\overline{\lambda}_\infty + n -2) F_j \sim - W_0^{(j)}(j)\varpi - \lambda_k\,B(0\star k)\,F_{jk}
+ \frac{\lambda_k}{\overline{\lambda}_\infty + n-2 }\,W_0^{(j)}(jk)\varpi,\\
&& W_0^{(j)}(jk)\varpi \sim (\overline{\lambda}_\infty + n-2)\,Z_{jk}^j,\\
&&W_0^{(j)}(j)\varpi \sim - (\overline{\lambda}_\infty + \lambda_j + n-2)\,F_j 
- \lambda_k\{F_k + B\!
\left(
\begin{array}{ccc}
 0 &\star   &k   \\
  0&\star   & j  \\  
\end{array}
\right)
F_{jk}
\},
\ee
where $\overline{\lambda}_\infty = \lambda_j + \lambda_k.$
As a result, we get the equality
\be
&&(2\lambda_k + n -2)\,F_j = (\lambda_k + n - 2)\,F_j +
\lambda_k\{F_k + B\!
\left(
\begin{array}{ccc}
 0 &\star   &k   \\
 0 & \star  &j   \\
\end{array}
\right)
F_{jk}
\} \\
&&
{}\qquad \qquad \qquad \qquad - \lambda_k\,B(0\star k)\,F_{jk} + \lambda_k\,Z_{jk}^j,
\ee
which implies (A.33).
Suppose now that $ J= \{j,k_2,\ldots,k_p\}$ such that every $k_\nu \ne j$.
By taking $\lambda_k = 0$ for all $k \in J^c$ and $\lambda_j = 0,$
we have the identity with respect to $\lambda_{k_2},\ldots, \lambda_{k_p}$
\be
Z_{j\,k_2\ldots k_p}^j = \sum_{\nu=2}^p \, Z_{j\,k_2\ldots\widehat{k_\nu}\ldots k_p,\,k_\nu}^j
+ \mathfrak{M}_{f_j^{-1}}\,W_0(k_2\ldots k_p)\varpi,
\ee
which is nothing else than (A.32).

%lm A.6
%\vspace{.5cm}

\bigskip

{\bf Lemma A.6}\quad
{\it
$Z_J^j,\,Z_{J,k}^j$ in the {\rm RHS} of {\rm (A.29)} equal the {\rm RHS} of {\rm (A.8), (A.9)}
respectively, so that $F_K$ appearing in the expression of $Z_J^j, \,Z_{J,k}^j$
are of weight at least $p -1,  p$.
}

%pf
\begin{pf}
Suppose first $p = 1$.  Then (A.4) is an immediate consequence of the definition, 
while  (A.6) coincides with (A.10).
We want to prove Lemma A.6 by induction on $p = |J| \,(p\ge 2).$ 
For the moment, denote the RHS of {\rm(A.8), (A.9)} by $\hat{Z}_{J}^j, \hat{Z}_{J,k}^j$
respectively. Then according to Lemma A.5, 
(A.23) are rewritten in terms of  $W_0^{(j)}(j\,j_2)\varpi$
as follows:
\be
W_0^{(j)}(j\,j_2)\varpi \sim (\lambda_\infty + n - 2)\,\hat{Z}_{j\, j_2}^j + \sum_{k=1;k\ne j,j_2}^{m} \lambda_k\,\hat{Z}_{{j\,j_2},k}^j.
\ee
This implies that  $Z_J^j = \hat{Z}_J^j,\, Z_{J,k}^j = \hat{Z}_{J,k}^j$ in the case $J = \{j, j_2\},\,p = 2.$ 
 We prove it in the case $p \ge 3$.
 Suppose that Lemma A.6 holds true for $|J| \le  p-1.$ 
 We must prove in the case $|J| = p.$
 Owing to the reccurrence relation  (A.32), we may assume that 
 $Z_J^j = \hat{Z}_J^j.$
 
  On the other hand,
as is seen from  {\rm (A.24), (A.26) and (A.28)},
a linear combination of  $F_K$ such that $|K| = p$ or $p+1$ appearing in $Z_{J,k}^j$
is the same as the one in $\hat{Z}_{J,k}^j$.

In view of (A.13), we see by induction that $Z_{J,k}^j$
is written by a linear combination of $F_K$ as follows:
\be
Z_{J,k}^j = c_k\,F_k + \sum_{L \subset J, 1\le |L| \le p-2} c_{k\,L}F_{k\,L}
+ \hat{Z}_{J,k}^j,
\ee
such that $k \in L^c,$
where $c_k, c_{kL}$ are constants independent of $\lambda$.
Notice that ${\rm weight}(\hat{Z}_J^j) \ge  p - 1, \,{\rm weight}(\hat{Z}_{J,k}^j) \ge   p,$
while ${\rm weight}(F_k) = 1, \,{\rm weight}(F_{kL}) \le  p - 1.$
We must prove that all $c_k = c_{k\,L} = 0$.
By Lemma A.5, $c_k,\,c_{k\,L}$
depend only on the coefficients of $f_k$ and $f_\nu\,(\nu\in J)$.
Hence $c_k, \,c_{k\,L}$ do not change
under the restriction
\beq
\lambda_h = 0\quad(h\in \partial_k J^c) .
\eeq

From now on, we fix $k$ and assume (A.34).
Theorem 1 still holds provided
that $\lambda_\infty$ are replaced by 
$\overline{\lambda}_\infty = \sum_{\nu\in J} \lambda_\nu + \lambda_k$.

For simplicity, we may assume that $J = \{12\ldots p\}, k = p+1$
without losing generality.
Furthermore, $c_k, c_{k L}$
do not change under the restriction
\be
\lambda_l = 0 \quad(l\ge p+2).
\ee

%p60
Operating $\mathfrak{M}_{f_j}T_{\varepsilon_j}$ on both sides of (A.29),
we have
\beq
&&\lambda_j \,W_0(J)\varpi =\mathfrak{M}_{f_j}T_{\varepsilon_j} W_0^{(j)}(J)\varpi 
\sim \mathfrak{M}_{f_j}T_{\varepsilon_j}\{(\overline{\lambda}_\infty + n - p) \,\hat{Z}_{J}^j +
\lambda_k\, Z_{J,k}^j \}\nonumber\\
&& {}\quad \quad \quad \quad \ \ = (\overline{\lambda}_\infty + n + 1 - p)\,\mathfrak{M}_{f_j}\hat{Z}_{J}^j
+  \lambda_k\{\, c_k\,\mathfrak{M}_{f_j}F_k\nonumber\\
&&{}\qquad \qquad  \qquad +\sum_{L \subset J, 1\le |L| \le p-2} \,c_{k\,L}\,\mathfrak{M}_{f_j}F_{k\,L}
+ \,\mathfrak{M}_{f_j}\hat{Z}_{J,k}^j\}.
\eeq

%p
$W_0(K)\varpi\,(L \subset \{k, J\}, 1\le |K| \le p+1)$ are linearly independent and make a basis of $H_\nabla^n(X, \Omega^\cdot(*S))$. 
Owing to Theorem 1 and  Propositions 34 and 36,
$\mathfrak{M}_{f_j}F_K$ can be represented explicitly as a linear combination of 
$W_0(L)\varpi$ and hence of $F_L\, (J \subset \{k, J\}, 1\le |L|\le p+1)$. 
As a linear combination of $F_K,$ denote by $[\mathfrak{M}_{f_j}F_K]_{r}$ the part of weight $r$ of $\mathfrak{M}_{f_j}F_K$
for $ r = 1, 2, \ldots, p +1$ (see Remark in \S4 for the weight).

We have
\be
&&\bigl[\mathfrak{M}_{f_j}F_k\bigr]_{1} = \{ -\lambda_k\, B(0\star k)+ B\!
\left(
\begin{array}{ccc}
 0 &j  &k   \\
 0 & \star  &k   \\
\end{array}
\right)\}\,F_k,\\
&&\bigl[\mathfrak{M}_{f_j}F_{k_1k_2}\varpi\bigr]_{1} = \frac{B\!
\left(
\begin{array}{ccc}
  0&j   &k_2   \\
  0&\star   &k_2   \\
\end{array}
\right)\,B(0\star k_2)
}
{B(0\star k_1k_2)}\,F_{k_2}
+ \frac{B\!
\left(
\begin{array}{ccc}
  0&j   &k_1   \\
  0&\star   &k_1   \\
\end{array}
\right)\,B(0\star k_1)
}{B(0\star k_1k_2)}\, F_{k_1},\\
&&\bigl[\mathfrak{M}_{f_j}F_{k_1k_2}\bigr]_{1} = 0 \quad(k\ne k_1,k_2),\\
&&\bigl[\mathfrak{M}_{f_j}F_{K}\bigr]_{1} = 0 \quad(|K| \ge 3),
\ee
since $W_0(k)\varpi = B(0\star k) \,F_k.$

 For more general $F_{k L}\,(j\notin L,\, L \subset J,\, 1\le |L| = q -1 \le p-2),$
we have
 \be
 &&\bigl[\mathfrak{M}_{f_j}F_{k\,L}\bigr]_{ q -1} =
 \frac{B\!
\left(
\begin{array}{ccc}
 0 &k   &  L   \\
 0 &j   &  L  \\ 
\end{array}
\right)
}{B(0\,k\,L)} \,F_{ L}
+ \sum_{\nu \in  L} 
\,\frac{B\!
\left(
\begin{array}{ccccc}
  0&\nu   &k& \partial_\nu L  \\
  0&j  &k&\partial_\nu L   \\
\end{array}
\right)
}
{B(0\,k\, L)} \, F_{k\,\partial_\nu L}.
 \ee
On the other hand, if $j\in L \subset J, 1\le |L| = q-1\le p -1,$ 
then
\be
\mathfrak{M}_{f_j}F_{k\,L} = F_{k\,\partial_j L}.
\ee

Hence,
\be
&&0 = (\lambda_k + n+1-p) \bigl[\mathfrak{M}_{f_j}\hat{Z}_J^j\bigr]_{1} +
\lambda_k\,c_k
\bigl\{
- \lambda_k,B(0\star k) + B\!
\left(
\begin{array}{ccc}
 0 &j   &k   \\
 0 &\star   &k   \\
\end{array}
\right)
 \bigr\}\,F_k\nonumber\\
&&{}\qquad \qquad \qquad \qquad + \lambda_k \sum_{l=1,\,l\ne j,k}^{n+1}\,c_{kl} \frac{B\!
\left(
\begin{array}{ccc}
 0 & j  &k   \\
 0 & \star  &k   \\
\end{array}
\right)
}{B(0\star k l)}
\bigl[\mathfrak{M}_{f_j}F_{kl}\bigr]_{1}.
\ee

$\bigl[\mathfrak{M}_{f_j}\hat{Z}_J^j \bigr]_{1}$ does not depend on $\lambda_k,$
because  $ {\rm weight}\, \hat{Z}_J^j \ge  p - 1\ge  2$.
Actually $\bigl[\mathfrak{M}_{f_j}\hat{Z}_J^j \bigr]_{1} = 0$ for $p \ge 4.$
This implies $c_k = 0.$

At second step
 notice that ${\rm weight}\, \mathfrak{M}_{f_j}\hat{Z}_J^j \ge  p - 1,\,
{\rm weight}\, \mathfrak{M}_{f_j}\hat{Z}_{J,k}^j \ge p-1$
 (see (A.8) and (A.9) and Proposition 21).
Comparing the part of weight $q - 1$ on both sides of (A.35), we have
\be
&&\sum_{q=2}^{p-1} \bigl(\sum_{j\notin L\subset J, |L| = q-1\le p-2}\, c_{kL}\{\frac{B\!
\left(
\begin{array}{ccc}
 0 &k   &L   \\
  0&j   & L  \\  
\end{array}
\right)
}{B(0\,k\,L)}\,F_L + 
\sum_{\nu\in L}\frac{B\!
\left(
\begin{array}{ccccc}
  0&\nu   &k& \partial_\nu L  \\
  0&j   &k&\partial_\nu L  \\
\end{array}
\right)
}{B(0\,k\,L)}\,F_{k\,\partial_\nu L}\}\\
&&{}\qquad \qquad \qquad \qquad \qquad + 
\sum_{j\in L \subset J, |L| = q-1\le p-2}\,c_{kL}\, F_{k\,\partial_j L}\bigr)
 = 0.
\ee
since $\bigl[\mathfrak{M}_{f_j}\hat{Z}_J^j\bigr]_{q - 1} = \bigl[\mathfrak{M}_{f_j}\hat{Z}_{J,k}^j\bigr]_{q - 1} = 0$,

Comparing the coefficients of $F_L$ on both sides, wee see that $c_{kL} = 0$
for $j\notin L.$
Comparing then the coefficients of $F_{k\,\partial_jL}$ on both sides 
we see that $c_{kL}= 0$ for $j \in L.$

By induction, $Z_{J,k}^j = \hat{Z}_{J,k}^j$ has thus been proved for all $J\,(2\le p\le n+1)$
and $k \in J^c$.\qed

\end{pf}

\bigskip

[{\bf proof of Theorem A}]

\medskip

Lemma A.6 implies that (A.8) and (A.9) hold true,  proving  Theorem A. 

\qed

%%%%%%%%%%%%%%%%%%%%%%%%%%%%%%%%%%%%%%%%%%%%

%

\bigskip

\bigskip

\noindent Kazuhiko AOMOTO,\\
5-1307 Hara, Tenpaku-ku, Nagoya-shi, 468-0015, Japan.\\
e-mail: kazuhiko@aba.ne.jp

\bigskip

\noindent Yoshinori MACHIDA,\\
4-9-37 Tsuji, Shimizu-shi, Shizuoka-shi, 424-0806, Japan. \\
%Numazu College of Technology, 3600 Ooka, Numazu-shi, Shizuoka 410-8501, Japan.\\
e-mail: yomachi212@gmail.com

\end{document}